\documentclass[a4paper,reqno,11pt]{article}

\pdfoutput=1

\usepackage{amsthm,enumerate,amssymb,amsmath,nccmath,mathtools,pgfmath,tikz,genyoungtabtikz,enumitem,setspace}
\usepackage{hyperref}
\hypersetup{colorlinks=true,linkcolor=blue,citecolor=red,urlcolor=blue}
\usepackage{scalerel}
\usepackage{mathrsfs}
\usepackage{blkarray}
\usepackage{tikz,tikz-cd}
\usepackage[a4paper,margin=0.8in]{geometry}
\usepackage{cleveref}
\usepackage{comment}
\usepackage[super]{nth}
\usepackage{multicol}
\usepackage{tcolorbox}
\tcbuselibrary{theorems}

\newtheorem{thm}{Theorem}[section]
\newtheorem{thm*}{Theorem}
\newtheorem{lem}[thm]{Lemma}
\newtheorem{prop}[thm]{Proposition}
\newtheorem{cor}[thm]{Corollary}
\newtheorem{ex}[thm]{Example}
\newtheorem{defn}[thm]{Definition}
\newtheorem{rmk}[thm]{Remark}

\numberwithin{equation}{section}

\crefname{thm}{Theorem}{Theorems}
\crefname{prop}{Proposition}{Propositions}
\crefname{lem}{Lemma}{Lemmas}
\crefname{cor}{Corollary}{Corollaries}
\crefname{conj}{Conjecture}{Conjectures}
\crefname{defn}{Definition}{Definitions}
\crefname{rmk}{Remark}{Remarks}
\crefname{section}{Section}{Sections}
\crefname{subsection}{Subsection}{Subsections}
\crefname{ex}{Example}{Examples}

\Crefname{thm}{Theorem}{Theorems}
\Crefname{prop}{Proposition}{Propositions}
\Crefname{lem}{Lemma}{Lemmas}
\Crefname{cor}{Corollary}{Corollaries}
\Crefname{conj}{Conjecture}{Conjectures}
\Crefname{defn}{Definition}{Definitions}
\Crefname{rmk}{Remark}{Remarks}
\Crefname{section}{Section}{Sections}
\Crefname{subsection}{Subsection}{Subsections}
\Crefname{ex}{Example}{Examples}

\newcommand\tts{\mathtt{S}}
\newcommand\ttt{\mathtt{T}}

\DeclareMathOperator{\im}{im}
\DeclareMathOperator{\spn}{span}
\DeclareMathOperator{\res}{res}
\DeclareMathOperator{\std}{Std}

\DeclareMathOperator{\dimn}{dim}
\DeclareMathOperator{\grdim}{grdim}
\DeclareMathOperator{\degr}{deg}

\DeclareMathOperator{\rem}{rem}
\DeclareMathOperator{\add}{add}
\DeclareMathOperator{\nor}{nor}
\DeclareMathOperator{\conor}{conor}
\DeclareMathOperator{\rad}{rad}
\DeclareMathOperator{\maxdeg}{maxdeg}
\DeclareMathOperator{\mindeg}{mindeg}
\DeclareMathOperator{\ind}{ind}
\DeclareMathOperator{\mx}{max}
\DeclareMathOperator{\soc}{soc}
\DeclareMathOperator{\md}{-mod}
\DeclareMathOperator{\sgn}{sgn}

\newcommand\plus{+}
\newcommand\minus{-}

\newcommand{\Mod}[1]{\ (\mathrm{mod}\ #1)}

\newcommand\emset{\varnothing}
\newcommand\gr{\Yfillcolour{lightgray}}

\newcommand\wh{\Yfillcolour{white}}
\newcommand\wht{\Yfillcolour{white}}

\tikzstyle{dotted}=                  [dash pattern=on \pgflinewidth off 2pt]
\tikzstyle{dashed}=                  [dash pattern=on 3pt off 3pt]
\tikzstyle{dashdotted}=              [dash pattern=on 3pt off 2pt on \the\pgflinewidth off 2pt]
\tikzstyle{densely dashdotted}=      [dash pattern=on 3pt off 1pt on \the\pgflinewidth off 1pt]
\tikzstyle{loosely dashdotted}=      [dash pattern=on 3pt off 4pt on \the\pgflinewidth off 4pt]
\allowdisplaybreaks
\makeatletter
\newcommand\Label[1]{&\refstepcounter{equation}(\theequation)\ltx@label{#1}&}

\makeatother

\allowdisplaybreaks

\author{Louise Sutton\footnote{Present address: National University of Singapore, 10 Lower Kent Ridge Road, Singapore 119076.}\ \footnote{Email address: matloui@nus.edu.sg.}
	\\	\normalsize Queen Mary University of London\\ Mile End Road\\ London E1 4NS}

\newcommand\runninghead[1]{\gdef\@runninghead{#1}}

\newcommand\auth{Louise Sutton}

\newcommand\toptitle{\date{}\maketitle\markboth{\auth}{\@runninghead}\pagestyle{myheadings}}

\begin{document}

\title{Specht modules labelled by hook bipartitions II}

\runninghead{Specht modules labelled by hook bipartitions II}

\toptitle

\begin{abstract}
We continue the study of Specht modules labelled by hook bipartitions for the Iwahori--Hecke algebra of type $B$ with $e\in\{3,4,\dots\}$ via the cyclotomic Khovanov--Lauda--Rouquier algebra $\mathscr{H}_n^{\Lambda}$.
Over an arbitrary field, we explicitly determine the graded decomposition submatrices for $\mathscr{H}_n^{\Lambda}$ comprising rows corresponding to hook bipartitions.
\end{abstract}

\begin{center}
\end{center}

\begin{center}
\end{center}

\section{Introduction}
	
The study of the representations of the $\mathbb{Z}$-graded \textit{cyclotomic Khovanov--Lauda--Rouquier algebras} (alternatively the \textit{cyclotomic quiver Hecke algebras}), denoted $\mathscr{H}_n^{\Lambda}$, has been motivated by their connection with the well-studied complex reflection groups and their deformations via  Brundan and Kleshchev's Graded Isomorphism Theorem in~\cite{bkisom}. This allows us to consider the Ariki--Koike algebras associated to a complex reflection group of type $G(l,1,n)$ as graded algebras.

The most important open question in the representation theory of the Ariki-Koike algebras is the Decomposition Number Problem.  One aims to understand the graded composition multiplicity $[S_{\lambda}:D_{\mu}\langle k\rangle]_v$ of the irreducible module, $D_{\mu}$, as a composition factor of the Specht module, $S_{\lambda}$, for all multipartitions $\lambda$ and for all regular multipartitions $\mu$.
Throughout this paper, we will fix $l=2$ and study the graded representation theory of the corresponding Iwahori--Hecke algebra of type $B$ from the perspective of $\mathscr{H}_n^{\Lambda}$.
In particular, we continue the study from~\cite{Sutton17I} of the special family of Specht modules labelled by \emph{hook bipartitions}, namely $S_{((n-m),(1^m))}$, as $\mathscr{H}_n^{\Lambda}$-modules. For the first time, we determine the corresponding graded decomposition numbers, which we observe are independent of the characteristic of the ground field.

\newtcolorbox{mybox}[1]{colback=black!10!white,
	colframe=black,fonttitle=\bfseries,
	title=#1}
\begin{mybox}{\large{Main Result}}
	Let $\mathbb{F}$ be arbitrary, $e\in\{3,4,\dots\}$, and $\lambda=((n-m),(1^m))$ with $m\in\{0,\dots,n\}$. We completely determine the graded decomposition numbers $[S_{\lambda}:D_{\mu}\langle k\rangle]_v$ for all regular bipartitions $\mu$.
\end{mybox}

Over a field of characteristic zero, we note that there exist recursive algorithms for determining decomposition numbers for the Ariki--Koike algebras. We know from Ariki's Categorification Theorem in~\cite{ariki96}, together with recent work of Brundan and Kleshchev~\cite{bk09}, that the graded decomposition numbers for the cyclotomic Khovanov--Lauda--Rouquier algebras can be determined from the canonical basis elements of the quantum affine algebra $U_q(\widehat{\mathfrak{sl}_e})$ via the LLT algorithm~\cite{LLT} in level one, and via an analogous algorithm~\cite{Fay10} in higher levels. 
While these are major breakthroughs in the field, the recursive nature of these algorithms means that explicit computations for all but sufficiently small $n$ are impossible, and the Decomposition Number Problem remains unsolved. 

In positive characteristic, we obtain the decomposition matrices for the Ariki--Koike algebras from the decomposition matrices in characteristic zero by post-multiplying them by certain adjustment matrices. However, there exists no analogue of the LLT algorithm for determining these adjustment matrices in positive characteristic, and moreover, we have very few explicit examples to hand. One of the most fundamental problems is to determine when the decomposition numbers in characteristic zero and positive characteristic coincide, and hence when the adjustment matrices are trivial. Due to Williamson's counterexamples for the symmetric groups~\cite{w13}, we now know that the long-standing James Conjecture~\cite[\S4]{j90} can no longer hope to provide a partial solution to this problem.
Except in a few cases, it is completely unknown when the adjustment matrices are trivial.
In level two, Brundan and Stroppel~\cite{bs11} and Hu and Mathas~\cite[Corollary B.5]{hm15} show that the decomposition numbers of the Iwahori--Hecke algebra of type $B$ do not depend on the characteristic of the ground field when $e$ is either infinite or sufficiently large. 
In level three, Lyle and Ruff~\cite{lr16} study certain blocks of the Ariki--Koike algebras, and determine that their corresponding adjustment matrices are, in fact, trivial for all quantum characteristics. This paper works in finite quantum characteristic $e\in\{3,4,\dots\}$, and adds to these recent developments by providing a special family of Specht modules for the Iwahori--Hecke algebra of type $B$ whose corresponding decomposition numbers are independent of the characteristic of the ground field.

In this paper, we add to this literature by studying the structure of Specht modules labelled by hook bipartitions.
We first recall from~\cite{Sutton17I} that the explicit presentation of $\mathscr{H}_n^{\Lambda}$ was used to determine composition series of $S_{((n-m),(1^m))}$ in which we defined its composition factors in terms of quotients either of the kernels or of the images of certain Specht module homomorphisms. In this way, we can write down explicit spanning sets for these composition factors in terms of standard basis elements of $S_{((n-m),(1^m))}$, which will later help us to determine certain properties of their gradings.
In general, we note that it is a non-trivial task to explicitly determine which of the regular multipartitions label the irreducible modules arising in the composition series of Specht modules. However, since we know that every composition factor of a Specht module arises as the head of a Specht module labelled by a regular multipartition, we are able to use Brundan and Kleshchev's $i$-restriction and $i$-induction functors~\cite{bk03} to find isomorphisms between the composition factors of $S_{((n-m),(1^m))}$ as presented in~\cite{Sutton17I} and the irreducible heads $D_{\mu}$ of certain Specht modules labelled by regular bipartitions.
We thus determine characteristic-free ungraded multiplicities $[S_{((n-m),(1^m))}:D_{\mu}]$ for all regular bipartitions $\mu$, and hence observe that the corresponding submatrices of the adjustment matrices are trivial. Furthermore, we completely determine the analogous graded composition multiplicities $[S_{((n-m),(1^m))}:D_{\mu}]_v$ by exploiting the combinatorial grading on Specht modules as defined in~\cite{bkw11}. 

We remark that one can alternatively keep track of the grading shifts throughout the preceding paper~\cite{Sutton17I} so that we immediately arrive at the graded results and hence implicitly recover the ungraded ones, however this method would give us little to no advantage since the resulting computations would be similar to those presented in this article. We instead only enter into the graded world in this paper to provide a distinction between the combinatorial calculations we now perform and those given in~\cite{Sutton17I} of the action of the $\mathscr{H}_n^{\Lambda}$-generators on standard basis elements of Specht modules.

The structure of this paper is as follows. In \cref{sec:back}, we present necessary background details of the graded representation theory of the cyclotomic Khovanov--Lauda--Rouquier algebras, and in particular we provide a brief overview of Specht modules labelled by hook bipartitions. In \cref{sec:onedims,sec:labelsgen}, we determine the composition factors of these Specht modules in terms of irreducible heads $D_{\mu}$ of Specht modules for certain regular bipartitions $\mu$. In doing so, it follows from~\cite{Sutton17I} that we completely determine the ungraded decomposition matrices of $\mathscr{H}_n^{\Lambda}$ corresponding to hook bipartitions; we present these results in \cref{sec:UDN}. Furthermore, by obtaining results in \cref{sec:GDS,sec:GDI} on the graded dimensions of Specht modules labelled by hook bipartitions and of their composition factors, we present the explicit graded decomposition numbers of $\mathscr{H}_n^{\Lambda}$ corresponding to hook bipartitions in \cref{sec:GDNum}.

\section{Background}\label{sec:back}

Throughout this paper, we let $\mathbb{F}$ be an arbitrary field and let $\mathfrak{S}_n$ be the symmetric group on $n$ letters. Let $q\in\mathbb{F}^{\times}$ be a cyclotomic $e$th root of unity such that $e \in \{3,4,\dots\}$; we call $e$ the \emph{quantum characteristic}. We set $I:=\mathbb{Z}/e\mathbb{Z}$ and identify $I$ with the set $\{0,1,\dots,e-1\}$. Recall that for a fixed \emph{level}, $l$, we let the \emph{$e$-multicharge of $l$} be the ordered $l$-tuple $\kappa=(\kappa_1,\kappa_2,\dots,\kappa_l)\in I^l$, with associated domaninant weight $\Lambda=\Lambda_{\kappa_1}+\dots +\Lambda_{\kappa_l}$ of level $l$. We refer the reader to~\cite[\S 2.2]{Sutton17I} for further details on the corresponding Lie-theoretic notation.

\subsection{Graded algebras and graded modules}

We familiarise the reader with the fundamental theory of graded algebras and modules; \cite{NO} provides a superb guide to graded representation theory. 

An $\mathbb{F}$-algebra $A$ is called \emph{graded}, more precisely $\mathbb{Z}$-graded, if there exists a direct sum decomposition $\textstyle{A=\bigoplus_{i\in\mathbb{Z}}A_i}$ such that $A_iA_j\subseteq{A_{i+j}}$ for all $i,j\in\mathbb{Z}$. An element in the summand $A_i$ is said to be \emph{homogeneous} of \emph{degree} $i$. For $a_i\in{A_i}$, we write $\degr(a_i)=i$.

Given a graded $\mathbb{F}$-algebra $A$, we say that the (left) $A$-module $M$ is \emph{$\mathbb{Z}$-graded} if there exists a direct sum decomposition $\textstyle{M=\bigoplus_{i\in\mathbb{Z}}M_i}$ such that $A_iM_j\subseteq{M_{i+j}}$ for all $i,j\in\mathbb{Z}$. We denote the abelian category of all finitely generated graded (left) $A$-modules by $A\md$.
If $M\in A\md$, then we obtain the module $M\langle {k} \rangle$ by shifting the grading on $M$ upwards by $k\in\mathbb{Z}$. For an indeterminate $v$, we set $M\langle {k} \rangle =v^kM$, so that the grading on $M\langle{k}\rangle$ is defined by $\left( M\langle{k}\rangle\right)_i=\left(v^kM\right)_i=M_{i-k}$. The \emph{graded dimension} of $M$ is defined to be the Laurent polynomial
\[\grdim(M)=\sum_{i\in\mathbb{Z}}\dimn(M_i)v^i \in\mathbb{N}[v,v^{-1}]. \]
Suppose that $A$ has a homogeneous anti-involution $\ast:A\rightarrow A$, and write $a^{\ast}$ for the image of $a\in A$ under this map. Then we define the \emph{dual} of $M$ to be the $\mathbb{Z}$-graded $A$-module
\[
M^{\circledast}=\bigoplus_{k\in\mathbb{Z}}\operatorname{Hom}_{\mathbb{F}}(M\langle k\rangle,\mathbb{F}),
\]
where the $A$-action is given by $(af)(m)=f(a^{\ast}m)$ for all $a\in A$, $f\in M^{\circledast}$ and $m\in M$.

We say that a \emph{graded composition series} for $M\in A\md$ is a filtration of graded submodules $\{0\}=M_0 \subset M_1 \subset \dots \subset M_{n-1} \subset M_n=M$ such that the quotients $M_{i}/M_{i-1}$ are irreducible for all $i\in\{1,\dots,n\}$, which we refer to as the \emph{graded composition factors of $M$}. The Jordan--H{\"o}lder theorem yields an analogous graded version, and thus it makes sense to study \emph{graded decomposition numbers} $[M:L]_v$ of $M$, where $L$ is a graded irreducible $A$-module. The graded multiplicity of $L$ as a composition factor of $M$ is defined to be the Laurent polynomial

\[[M:L]_v=\sum_{i\in\mathbb{Z}}[M:L\langle{i}\rangle]v^i \in\mathbb{N}[v,v^{-1}].\]

Note that by setting $v=1$ into the above definitions, we recover the ungraded analogues.

\subsection{Multipartitions, Young diagrams and tableaux}

We recall from \cite[\S 3]{Sutton17I} basic combinatorial notation and definitions in this section.

We write $\mathscr{P}_n^l$ for the set of all $l$-multipartitions of $n$, and in particular, we write $\varnothing$ for the empty multipartition. Let $\lambda=(\lambda^{(1)},\dots,\lambda^{(l)})\in\mathscr{P}_n^l$. We define the \emph{Young diagram} of $\lambda$ to be
\[ [\lambda]:= \left\{ (i,j,m)\in \mathbb{N} \times \mathbb{N} \times \{1,\dots,l\}\ \middle|\  1\leqslant{j}\leqslant{\lambda_i^{(m)}} \right\}. \]
We draw the $i$th component $\lambda^{(i)}$ of $[\lambda]$ above its $(i+1)$th component $\lambda^{(i+1)}$ for all $i\in\mathbb{N}$.
Each element $(i,j,m)\in[\lambda]$ is called a \emph{node} of $\lambda$, and in particular, an $(i,j)$-node of the $m$th component $\lambda^{(m)}$.
We say that the node $(i_1,j_1,m_1)\in[\lambda]$ lies \emph{strictly above} the node $(i_2,j_2,m_2)\in[\lambda]$ if either $i_1<i_2$ and $m_1=m_2$ or $m_1< m_2$.

We say that $A\in[\lambda]$ is a \emph{removable node} for $\lambda$ if ${{[\lambda]\backslash\{A\}}}$ is a Young diagram of an $l$-multipartition of $n-1$. Similarly, we say that $A\not\in[\lambda]$ is an \emph{addable node} for $\lambda$ if ${{[\lambda]\cup\{A\}}}$ is a Young diagram of an $l$-multipartition of $n+1$.

A \emph{$\lambda$-tableau} $\ttt$ is a bijection $\ttt:[\lambda]\rightarrow\{1,\dots,n\}$. We call $\ttt$ \emph{standard} if the entries in each row increase from left to right along the rows of each component, and the entries in each column increase from top to bottom down the columns of each component. We denote the set of all standard $\lambda$-tableaux by $\std(\lambda)$. The \emph{column-initial tableau} $\ttt_{\lambda}$ is the $\lambda$-tableau whose entries $1,\dots,n$ appear in order down consecutive columns, working from left to right in components $l,l-1,\dots,1$, in turn.

\subsection{Residues and degrees}

We fix an $e$-multicharge $\kappa=(\kappa_1,\dots,\kappa_l)\in{I^l}$.
The \emph{$e$-residue} of a node $A=(i,j,m)$ lying in the space $\mathbb{N}\times\mathbb{N}\times\{1,\dots,l\}$ is defined to be
\[
\res A := \kappa_m+j-i \Mod{e}.
\]
We say that an \emph{$i$-node} is a node of residue $i$.

Let $\ttt$ be a $\lambda$-tableau. We write $r=\ttt(i,j,m)$ to denote that the integer entry $r$ lies in node $(i,j,m)\in[\lambda]$, and set $\res_{\ttt}(r)=\res (i,j,m)$. The \emph{residue sequence} of $\ttt$ is defined to be
\[\mathbf{i}_{\ttt}=(\res_{\ttt}(1),\dots,\res_{\ttt}(n)).\]

We define the \emph{degree} of an addable $i$-node $A$ of $\lambda\in\mathscr{P}_n^l$ to be
\begin{align*}
d^A(\lambda):=\ 
&\#\left\{
\text{addable $i$-nodes of $\lambda$ strictly above $A$}
\right\}
-\#\left\{
\text{removable $i$-nodes of $\lambda$ strictly above $A$}
\right\}.
\end{align*}

Let $\ttt\in\std(\lambda)$ be such that $n$ lies in node $A$ of $\lambda$. We set $\deg (\varnothing):=0$, and define the \emph{degree} of $\ttt$ recursively via
\[
\deg (\ttt):=d^A(\lambda)+\deg (\ttt_{\leqslant{n-1}}),
\]
where $\ttt_{\leqslant{n-1}}$ is the standard tableau obtained by removing node $A$ from $\ttt$.

\begin{ex}\label{ex:deg}
	Let $e=3$ and $\kappa=(0,0)$. There are five standard $((1),(1^4))$-tableaux, namely
	\[
	\ttt_1=\young(1,,2,3,4,5),\:
	\ttt_2=\young(2,,1,3,4,5),\:
	\ttt_3=\young(3,,1,2,4,5),\:
	\ttt_4=\young(4,,1,2,3,5),\:
	\ttt_5=\young(5,,1,2,3,4).
	\]
	We find the degree of $\ttt_1$ as follows. We note that the degree of any node in the first row of the first component is $0$, so $d^{(1,1,1)}=0$ and hence $\degr(\ttt_{\leqslant{1}})=0$. Observe that 
	\[
	\ttt_{\leqslant 2}=\young(1,,2),
	\text{which has $3$-residues }
	\young(!\gr0,,!\wht0).
	\]
	Thus $(1,1,2)$ has removable $0$-node $(1,1,1)$ (shaded above), and hence $d^{(1,1,2)}=-1$. Now observe
	\[
	\ttt_{\leqslant 3}=\young(1,,2,3),
	\text{which has $3$-residues }
	\gyoung(;0,!<\Ylinethick{2pt}>\:,,!<\Ylinethick{0.5pt}>0,;2).
	\]
	Thus $(2,1,2)$ has addable $2$-node $(2,1,1)$ (outlined above), and hence $d^{(2,1,2)}=1$. Now observe
	\[
	\ttt_{\leqslant 4}=\young(1,,2,3,4),
	\text{which has $3$-residues }
	\gyoung(;0!<\Ylinethick{2pt}>\:,,!<\Ylinethick{0.5pt}>0!<\Ylinethick{2pt}>\:,!<\Ylinethick{0.5pt}>2,;1).
	\]
	Thus $(3,1,2)$ has addable $1$-nodes $(1,2,1)$ and $(1,2,2)$ (outlined above), and hence $d^{(3,1,2)}=2$. We finally observe that
	\[
	\ttt_1=\young(1,,2,3,4,5),
	\text{which has $3$-residues }
	\young(!\gr0,,!\wht0,2,1,0).
	\]
	Thus $(4,1,2)$ has removable $0$-node $(1,1,1)$ (shaded above), and hence $d^{(4,1,2)}=-1$. Hence
	\[
	\deg (\ttt_1)=d^{(1,1,1)}+d^{(1,1,2)}+d^{(2,1,2)}+d^{(3,1,2)}+d^{(4,1,2)}=1.
	\]
	
	Similarly, one can find that $\deg (\ttt_2)=\deg (\ttt_5)=3$, $\deg (\ttt_3)=2$ and $\deg (\ttt_4)=1$.
\end{ex}

\subsection{Cylotomic Khovanov--Lauda--Rouquier algebras and Specht modules}

The presentation of the cyclotomic Khovanov--Lauda--Rouquier algebra, $\mathscr{H}_{n}^{\Lambda}$, introduced independently by Khovanov and Lauda in~\cite{kl09} and Rouquier in~\cite{rouq08} endows $\mathscr{H}_{n}^{\Lambda}$ with a canonical $\mathbb{Z}$-grading. We know from Brundan and Kleshchev's \emph{Graded Isomorphism Theorem} \cite[Main Theorem]{bkisom} that $\mathscr{H}_{n}^{\Lambda}$ is isomorphic to a cyclotomic Hecke algebra (of type $A$).

We refer the reader to \cite{kmr} for the construction of \emph{Specht modules}, $S_{\lambda}$, over the cyclotomic Khovanov--Lauda--Rouquier algebras, which are indexed by multipartitions $\lambda$ and generated by the element $z_{\lambda}$ as an $\mathscr{H}_n^{\Lambda}$-module. We study (column) Specht modules as given in \cite{kmr}, which are dual to those given in \cite{bkw11} and consistent with James' classical construction of Specht modules over $\mathbb{F}\mathfrak{S}_n$.
For a $\lambda$-tableau $\ttt$, we recall from \cite[\S4]{Sutton17I} that for a reduced expression of $w_{\ttt}\in\mathfrak{S}_n$ such that $w_{\ttt}\ttt_{\lambda}=\ttt$, we can define the vector $v_{\ttt}=\psi_{w_{\ttt}}z_{\lambda}$ for some element $\psi_{w_{\ttt}}\in\mathscr{H}_n^{\Lambda}$ associated to the reduced expression of $w_{\ttt}$.
In general, we note that the vector $v_{\ttt}$ depends on the choice of a reduced expression of $w_{\ttt}$.
We observe that the existence of these vectors ensures that Specht modules naturally inherit a $\mathbb{Z}$-grading from $\mathscr{H}_n^{\Lambda}$.

\begin{thm}\cite[Corollary 4.6]{bkw11} and \cite[Proposition 7.14]{kmr}
	\label{gradedmod}
	Let $\lambda\in\mathscr{P}_n^l$. Then the set of vectors \[\left\{v_{\ttt}\ \middle |\  \ttt\in\std(\lambda)\right\}\]
	is a homogeneous $\mathbb{F}$-basis of $S_{\lambda}$ of degree determined by $\degr(v_{\ttt})=\degr(\ttt)$.
\end{thm}
Recall that this basis is called the \emph{standard homogeneous basis of $S_{\lambda}$}.
We can now define the \emph{graded dimensions} of Specht modules using the degree function on standard tableaux as follows.

\begin{defn}\label{defn:grdim}
	Let $\lambda\in\mathscr{P}_n^l$. 
	Then the \emph{graded dimension} of $S_{\lambda}$ is defined to be
	\[
	\grdim \left(S_{\lambda}\right):=\sum_{\ttt\in\std(\lambda)}v^{\deg (\ttt)}.
	\]
\end{defn}
We thus note that the graded dimensions of Specht modules depend only on the quantum characteristic $e$ and not directly on the ground field $\mathbb{F}$.

\begin{ex}
	Let $e=3$ and $\kappa=(0,0)$. Following \cref{ex:deg}, we know that 
	\[\grdim \left(S_{((1),(1^4))}\right)
	= v^{\degr(\ttt_1)}+v^{\degr(\ttt_2)}+v^{\degr(\ttt_3)}+v^{\degr(\ttt_4)}+v^{\degr(\ttt_5)}
	= 2v^3+v^2+2v.\]
\end{ex}

\subsection{Regular multipartitions}

We introduce numerous combinatorial definitions following~\cite{kbranch3}, most of which date back to~\cite{MiMi}, and we adopt notation introduced by Fayers in~\cite{Fay16}.

Let $\lambda\in\mathscr{P}_n^l$. We denote the total number of removable $i$-nodes of $\lambda$ by $\rem _i(\lambda)$, and we denote the total number of addable $i$-nodes of $\lambda$ by $\add _i(\lambda)$. We write the $l$-multipartition obtained by removing all of the removable $i$-nodes from $\lambda$ as ${\lambda^{\triangledown{i}}}$, and we write the $l$-multipartition obtained by adding all of the addable $i$-nodes to $\lambda$ as ${\lambda^{\vartriangle{i}}}$.

We define the \emph{$i$-signature of $\lambda\in\mathscr{P}_n^l$} by reading the Young digram $[\lambda]$ from the top of the first component down to the bottom of the last component, writing a $+$ for each addable $i$-node and writing a $-$ for each removable $i$-node, where the leftmost $+$ corresponds to the highest addable $i$-node of $\lambda$. We obtain the \emph{reduced $i$-signature of $\lambda$} by successively deleting all adjacent pairs $+-$ from the $i$-signature of $\lambda$, always of the form $-\dots-+\dots+$.

\begin{ex}
Let $e=3$, $\kappa=(0,0)$ and $\lambda=((7,4^2),(4))$. The $3$-residues of $\lambda$, as well as the $0$-addable and $0$-removable nodes of $\lambda$ are labelled, respectively, as follows
\begin{equation}\label{eq:sig}
\young(0120120,2012,1201,,,0120),
\qquad\qquad
\gyoung(;;;;;;;\minus,;;;;:\plus,;;;;,:\plus,,;;;;\minus).
\end{equation}
Thus, by removing all of the removable $0$-nodes from $\lambda$ (corresponding to the outlined nodes below), and respectively, adding all of the addable $0$-nodes of $\lambda$ (corresponding to the shaded nodes below) we have the following Young diagrams of multipartitions
\[
\left[\lambda^{\triangledown_0}\right]=
\gyoung(;;;;;;!<\Ylinethick{1.5pt}>;,!<\Ylinethick{0.5pt}>;;;;,;;;;,,;;;!<\Ylinethick{1.5pt}>;),
\qquad\qquad
\left[\lambda^{\vartriangle_0}\right]=
\gyoung(;;;;;;;,;;;;!\gr;,!\wh;;;;,!\gr;,,!\wh;;;;).
\]
Referring to \eqref{eq:sig}, the $0$-signature of $\lambda$ is $-++-$ (corresponding to the $-$ and $+$ labels from top to bottom in the diagram), and the reduced $0$-signature is $-+$ (corresponding to the nodes $(1,7,1)$ and $(2,5,1)$, respectively).
\end{ex}

The removable $i$-nodes corresponding to the $-$ signs in the reduced $i$-signature of $\lambda$ are called the \emph{normal} $i$-nodes of $\lambda$, and similarly, the addable $i$-nodes corresponding to the $+$ signs in the reduced $i$-signature of $\lambda$ are called the \emph{conormal} $i$-nodes of $\lambda$. We denote the total number of normal $i$-nodes of $\lambda$ by $\nor_i(\lambda)$ and the total number of conormal $i$-nodes of $\lambda$ by $\conor_i(\lambda)$. The lowest normal $i$-node of $[\lambda]$, if there is one, is called the \emph{good} $i$-node of $\lambda$, which corresponds to the last $-$ sign in the $i$-signature of $\lambda$. Similarly, the highest conormal $i$-node of $[\lambda]$, if there is one, is called the \emph{cogood} $i$-node of $\lambda$, which corresponds to the first $+$ sign in the $i$-signature of $\lambda$. 

For $r\in\{0,\dots,\nor_i(\lambda)\}$, we denote the multipartition obtained from $\lambda$ by removing the $r$ lowest normal $i$-nodes of $\lambda$ by $\lambda\downarrow_i^r$, and for $r\in\{0,\dots,\conor_i(\lambda)\}$, we denote the multipartition obtained from $\lambda$ by adding the $r$ highest conormal $i$-nodes of $\lambda$ by $\lambda\uparrow_i^r$. We set $\uparrow_i:=\uparrow_i^1$ when adding the cogood $i$-node of $\lambda$ and we set $\downarrow_i:=\downarrow_i^1$ when removing the good $i$-node of $\lambda$. It is easy to see that $A$ is a cogood $i$-node of $\lambda\in\mathscr{P}_n^l$ if and only if $A$ is a good $i$-node of $\lambda\cup\{A\}$. The operators $\uparrow_i^r$ and $\downarrow_i^r$ act inversely on a multipartition $\lambda\in\mathscr{P}_n^l$ in the following sense:
\begin{equation}\label{eq:arrows}
\lambda\downarrow_i^r\uparrow_i^r=\lambda
\quad
\left(0\leqslant{r}\leqslant{\nor_i(\lambda)}\right);
\qquad
\lambda\uparrow_i^s\downarrow_i^s=\lambda
\quad
\left(0\leqslant{s}\leqslant{\conor_i(\lambda)}\right).
\end{equation}

We define the set of all \emph{regular $l$-multipartitions of $n$} to be
\begin{equation}\label{def:regular}
\mathscr{RP}_n^l:=\left\{
\varnothing\uparrow_{i_1}\dots\uparrow_{i_n}\ \middle| \ i_1,\dots,i_n\in{I}
\right\}.
\end{equation}
If a multipartition $\lambda$ lies in $\mathscr{RP}_n^l$, then $\lambda$ is called \emph{regular}.
Hence $\lambda\in\mathscr{P}_n^l$ is \emph{regular} if and only if $[\lambda]$ is obtained by successively adding cogood nodes to $\varnothing$.

\subsection{Graded irreducible $\mathscr{H}_n^{\Lambda}$-modules}

In this section, we review a classification of the graded irreducible $\mathscr{H}_n^{\Lambda}$-modules. It is well known that the Specht module $S_{\lambda}$ has the quotient $D_{\lambda}:=S_{\lambda}/\rad S_{\lambda}$ for each $\lambda\in\mathscr{P}_n^l$, where the radical of $S_{\lambda}$ is defined from a homogeneous symmetric bilinear form on $S_{\lambda}$ of degree zero (see \cite[\S{2}]{hm10} for details). We know that each $D_{\lambda}$ is either absolutely irreducible or zero by \cite[Lemma 2.9]{hm10}, and moreover, $D_{\lambda}$ is absolutely irreducible if and only if $\lambda\in\mathscr{RP}_n^l$ \cite[Corollary 5.11]{hm10}.

\begin{thm}\cite[Theorem 4.11]{bk09} and \cite[Proposition 2.18]{hm10}\label{thm:selfdual}
\begin{enumerate}
\item{
$\left\{
D_{\lambda}\langle{i}\rangle\ \middle|\ \lambda\in\mathscr{RP}_n^l,i\in\mathbb{Z}
\right\}$
is a complete set of pairwise non-isomorphic irreducible graded $\mathscr{H}_{n}^{\Lambda}$-modules.
}
\item{
For all $\lambda\in\mathscr{RP}_n^l$, $D_{\lambda}\cong{D_{\lambda}^{\circledast}}$ as graded $\mathscr{H}_n^{\Lambda}$-modules.
}
\end{enumerate}
\end{thm}

\subsection{Graded decomposition numbers}\label{back:decomp}

Decomposition numbers record information about the structure of Specht modules. We denote the \emph{ungraded decomposition number} by $d_{\lambda,\mu}=
[S_{\lambda}:D_{\mu}]$ where $\lambda\in\mathscr{P}_n^l$ and $\mu\in\mathscr{RP}_n^l$, which is the multiplicity of $D_{\mu}$ appearing as a composition factor of $S_{\lambda}$.

We denote the ungraded decomposition matrix for $\mathscr{H}_n^{\Lambda}$ by $(d_{\lambda,\mu})$, and we write $(d_{\lambda,\mu}^{\mathbb{F}})$ when we want to emphasise the ground field. It is well known that we can compute the ungraded decomposition matrices, $(d_{\lambda,\mu}^{\mathbb{C}})$, for $\mathscr{H}_n^{\Lambda}$ via the generalised LLT algorithm given by Fayers in~\cite{Fay10}, whereas determining decomposition numbers in positive characteristics is an open problem.
We know from~\cite{bk09} that there exists an \emph{adjustment matrix} $(a_{\nu,\mu}^{\mathbb{F}})$ such that $(d_{\lambda,\nu}^{\mathbb{F}})=(d_{\lambda,\nu}^{\mathbb{C}})(a_{\nu,\mu}^{\mathbb{F}})$ where $\nu,\mu\in\mathscr{RP}_n^l$, but there exists no algorithm for determining the entries in this matrix.

We know from \Cref{gradedmod} that we can endow Specht modules with a $\mathbb{Z}$-grading, and since there exists a graded version of the Jordan--H\"{o}lder theorem, we can study their graded composition factors. We define the \emph{graded decomposition number} to be
\[
d_{\lambda,\mu}(v)=
\left[S_{\lambda}:D_{\mu}\right]_v:=
\sum_{i\in\mathbb{Z}} \left[S_{\lambda}:D_{\mu}\langle i\rangle\right]v^i
\in\mathbb{N}\left[v,v^{-1}\right],
\]
where $\lambda\in\mathscr{P}_n^l$ and $\mu\in\mathscr{RP}_n^l$.
We record these graded multiplicities in a \emph{graded decomposition matrix}, $(d_{\lambda\mu}(v))$, where its rows are indexed by multipartitions and its columns are indexed by regular multipartitions. 

The following result for $\mathscr{H}_n^{\Lambda}$ is a more general version of \cite[Corollary 12.2]{James} for $\mathbb{F}\mathfrak{S}_n$.

\begin{thm}\cite[Theorem 3.9 and Corollary 5.15]{bk09}
Let $\lambda\in\mathscr{P}_n^l$ and $\mu\in\mathscr{RP}_n^l$. Then
\begin{enumerate}
\item{ $d_{\mu,\mu}(v)=1$; }
\item{ $d_{\lambda,\mu}(v)\neq{0}$ only if $\mu\unrhd\lambda$.
}
\end{enumerate}
Moreover, if $\mathbb{F}=\mathbb{C}$ then $d_{\lambda,\mu}(v)\in v\mathbb{N}(v)$ whenever $\mu\rhd\lambda$.
\end{thm}

We denote the \emph{graded adjustment number} by $a_{\nu,\mu}^{\mathbb{F}}(v)$.

\begin{thm}\cite[Theorem 5.17]{bk09}
Let $\lambda\in\mathscr{P}_n^l$ and $\mu\in\mathscr{RP}^l_n$. Then
\[
d_{\lambda,\mu}^{\mathbb{F}}(v)=\sum_{\nu\in\mathscr{RP}_n^l}d_{\lambda,\nu}^{\mathbb{C}}(v)a_{\nu,\mu}^{\mathbb{F}}(v),
\]
for some $a_{\nu,\mu}^{\mathbb{F}}(v)\in\mathbb{N}[v,v^{-1}]$ with $a_{\nu,\mu}^{\mathbb{F}}(v)=a_{\nu,\mu}^{\mathbb{F}}(v^{-1})$.
\end{thm}

\subsection{Induction and restriction of $\mathscr{H}_n^{\Lambda}$-modules}\label{indres}

The \emph{Decomposition Number Problem} for $\mathscr{H}_n^{\Lambda}$ of determining the multiplicities $[S_{\lambda}:D_{\mu}]$ for all $\lambda\in\mathscr{P}_n^l$ and for all $\mu\in\mathscr{RP}_n^l$ is equivalent to the \emph{Branching Problem} of determining the multiplicities \[\left[\res_{\mathscr{H}_{n-1}^{\Lambda}}^{\mathscr{H}_n^{\Lambda}}D_{\lambda}:D_{\mu}\right]\]
for all $\lambda\in\mathscr{P}_n^l$ and for all $\mu\in\mathscr{RP}_n^l$.
The restriction of the ordinary representations of the symmetric group and their composition factors are well understood via the \emph{Classical Branching Rule for $\mathbb{F}\mathfrak{S}_n$} (for example, see~\cite[Theorem 9.2]{James}), which was first extended to the Ariki--Koike algebras (or the cyclotomic Hecke algebras) by Ariki--Koike~\cite[Corollary 3.12]{ak94}, and which has recently been extended to the cyclotomic Khovanov--Lauda--Rouquier algebras by Mathas~\cite{Mathas17}.

We first introduce Brundan and Kleshchev's $i$-restriction and $i$-induction functors, $e_i$ and $f_i$ respectively, acting on $\mathbb{F}\mathfrak{S}_n$-modules, as given in Section $2.2$ of~\cite{bk03}. These functors are exact, and originate from Robinson~\cite{Rob}; we extend these functors to act on $\mathscr{H}_n^{\Lambda}$-modules.

Let $M$ be an $\mathscr{H}_n^{\Lambda}$-module. For $i\in\mathbb{Z}/e\mathbb{Z}$, there are \emph{$i$-restriction functors} $e_i:\mathscr{H}_n^{\Lambda}\md\rightarrow\mathscr{H}_{n-1}^{\Lambda}\md$, 
and \emph{$i$-induction functors}
$f_i:\mathscr{H}_n^{\Lambda}\md\rightarrow
\mathscr{H}_{n+1}^{\Lambda}\md$, such that~\cite[Lemma 2.5]{bk03}
\begin{equation}\label{eq:indres}
\res^{\mathscr{H}_n^{\Lambda}}_{\mathscr{H}_{n-1}^{\Lambda}}M
\cong
\bigoplus_{i\in\mathbb{Z}/e\mathbb{Z}}
e_i M
\quad\text{and}\quad
\ind_{\mathscr{H}_n^{\Lambda}}^{\mathscr{H}_{n+1}^{\Lambda}}M
\cong\bigoplus_{i\in\mathbb{Z}/e\mathbb{Z}}f_iM.
\end{equation}

For $i\in\mathbb{Z}/e\mathbb{Z}$ and $r\geqslant{0}$, there exist the \emph{divided power $i$-restriction functors}
$e_i^{(r)}:\mathscr{H}_n^{\Lambda}\md\rightarrow\mathscr{H}_{n-r}^{\Lambda}\md$ and the \emph{divided power induction $i$-functors}
$f_i^{(r)}:\mathscr{H}_n^{\Lambda}\md\rightarrow
\mathscr{H}_{n+r}^{\Lambda}\md$,
which satisfy~\cite[Lemma 2.6]{bk03}
\[e_i^rM\cong\bigoplus_{k=1}^{r!}
e_i^{(r)}M
\quad\text{and}\quad
f_i^rM\cong\bigoplus_{k=1}^{r!}
f_i^{(r)}M.
\]

For a non-zero $\mathscr{H}_n^{\Lambda}$-module $M$, we define
\begin{equation}
\epsilon_i(M)=\mx\left\{r\geqslant{0}\ \middle|\ e_i^{(r)}M\neq{0}\right\}
\quad\text{and}\quad
\varphi_i(M)=\mx\left\{r\geqslant{0}\ \middle|\ f_i^{(r)}M\neq{0}\right\}.
\end{equation}
We now set 
\[e_i^{(\mx)}M=e_i^{(\epsilon_iM)}M
\qquad\text{and}\qquad
f_i^{(\mx)}M=f_i^{(\varphi_iM)}M.\]

By refining the Branching Rule for $\mathscr{H}_n^{\Lambda}$-modules, we obtain~\cite[Lemma 4.1]{Fay16} and its analogue.

\begin{lem}\label{Spechtmax}
Let $i\in\mathbb{Z}/e\mathbb{Z}$ and $\lambda\in\mathscr{P}_n^l$.
\begin{enumerate}
\item{
Then $\epsilon_i\left(S_{\lambda}\right)=\rem_i(\lambda)$ and
$
e_i^{(\mx)}S_{\lambda}
\cong
S_{\lambda^{\triangledown{i}}}
$.
}
\item{
Then $\varphi_i\left(S_{\lambda}\right)=\add_i(\lambda)$ and
$
f_i^{(\mx)}S_{\lambda}\cong
S_{\lambda^{\vartriangle{i}}}
$.
}
\end{enumerate}
\end{lem}

\subsection{Modular branching rules for $\mathscr{H}_n^{\Lambda}$-modules}

Kleshchev developed the analogous theory for restricting the modular representations of the symmetric group~\cite{kbranch1,kbranch2,kbranch3}, which Brundan extended to Hecke algebras of type $A$~\cite{brundan98}. These modular branching rules were generalised for cyclotomic Hecke algebras, proven by Ariki in the proof of~\cite[Theorem 6.1]{ariki06}. Thus modular branching rules for the cyclotomic Khovanov--Lauda--Rouquier algebras make sense, which we note here.

\begin{thm}\cite[\S{2.6}]{bk03}\label{thm:modbranch}
Let $i\in\mathbb{Z}/e\mathbb{Z}$ and $\lambda\in\mathscr{RP}_n^l$.
\begin{enumerate}
	\item Then $\epsilon_i\left({D_{\lambda}}\right)=\nor_i(\lambda)$ and $e_i^{(\mx)}D_{\lambda}\cong D_{\lambda\downarrow^{\nor_i(\lambda)}_i} $.
	\item Then $\varphi_i\left({D_{\lambda}}\right)=\conor_i(\lambda)$ and $f_i^{(\mx)}D_{\lambda}\cong D_{\lambda\uparrow^{\conor_i(\lambda)}_i} $.
\end{enumerate}
\end{thm}

\begin{ex}\label{ex:br}
	Let $e=3$, $\kappa=(0,2)$ and $\lambda=((9,6,2^2,1),(4,3,2))$. Since we can obtain $\lambda$ from $(\varnothing,\varnothing)$ by adding certain conormal nodes as follows
	\[\lambda=(\varnothing,\varnothing)\uparrow_2\uparrow_1\uparrow_0\uparrow_2\uparrow_0\uparrow_1^2\uparrow_2^2\uparrow_0^4\uparrow_1^4\uparrow_2\uparrow_0\uparrow_2^4\uparrow_0^2\uparrow_1^3\uparrow_2,\]
	we know from (\ref{def:regular}) that $\lambda$ is a regular bipartition. 
	We observe from~\Cref{modbranchfig} the $3$-residues of $\lambda$, together with its addable nodes. Thus $\lambda$ has $2$-signature $-+--++$ and reduced $2$-signature $--++$. One can also observe that we have drawn the bipartitions obtained from $\lambda$ by: 1) removing all of the normal $2$-nodes of $\lambda$ (outlined in~\Cref{modbranchfig}), corresponding to the $-$ signs in the reduced $2$-signature of $\lambda$, and 2) adding all of the conormal $2$-nodes of $\lambda$ (shaded in~\Cref{modbranchfig}), corresponding to the $+$ signs in the reduced $2$-signature of $\lambda$.
	It thus follows from \cref{thm:modbranch} that
	\begin{align*}
	e_2^{(2)}D_{\lambda}\cong{D_{((8,6,2^2,1),(3^2,2))}};\qquad
	f_2^{(2)}D_{\lambda}\cong{D_{((9,6,2^2,1),(4,3^2,1))}}.
	\end{align*}
	\begin{figure}
	\[
	\begin{tikzpicture}
	\tgyoung(0cm,0cm,;0;1;2;0;1;2;0;1;2:0,;2;0;1;2;0;1:2,;1;2:0,;0;1,;2:0,:1,,;2;0;1;2:0,;1;2;0:1,;0;1:2,:2)
	\tgyoung(-5cm,-6cm,;;;;;;;;!<\Ylinethick{2pt}>;,!<\Ylinethick{0.5pt}>;;;;;;,;;,;;,;,,;;;!<\Ylinethick{2pt}>;,!<\Ylinethick{0.5pt}>;;;,;;)
	\tgyoung(5cm,-6cm,;;;;;;;;;,;;;;;;,;;,;;,;,,;;;;,;;;,;;!\gr;,;)
	\draw[->](-0.5,-3)
	to node[above]{$e_2^{(2)}$}
	(-3,-5);
	\draw[->](3,-3)
	to node[above]{\ \ $f_2^{(2)}$}
	(5.5,-5);
	\end{tikzpicture}
	\]
	\caption{\footnotesize{The $3$-residues of $((9,6,2^2,1),(4,3,2))$, $((9,6,2^2,1),(4,3,2))\downarrow^{\nor_2(\lambda)}_2$ and $((9,6,2^2,1),(4,3,2))\uparrow^{\conor_2(\lambda)}_2$.}}
	\label{modbranchfig}
	\end{figure}
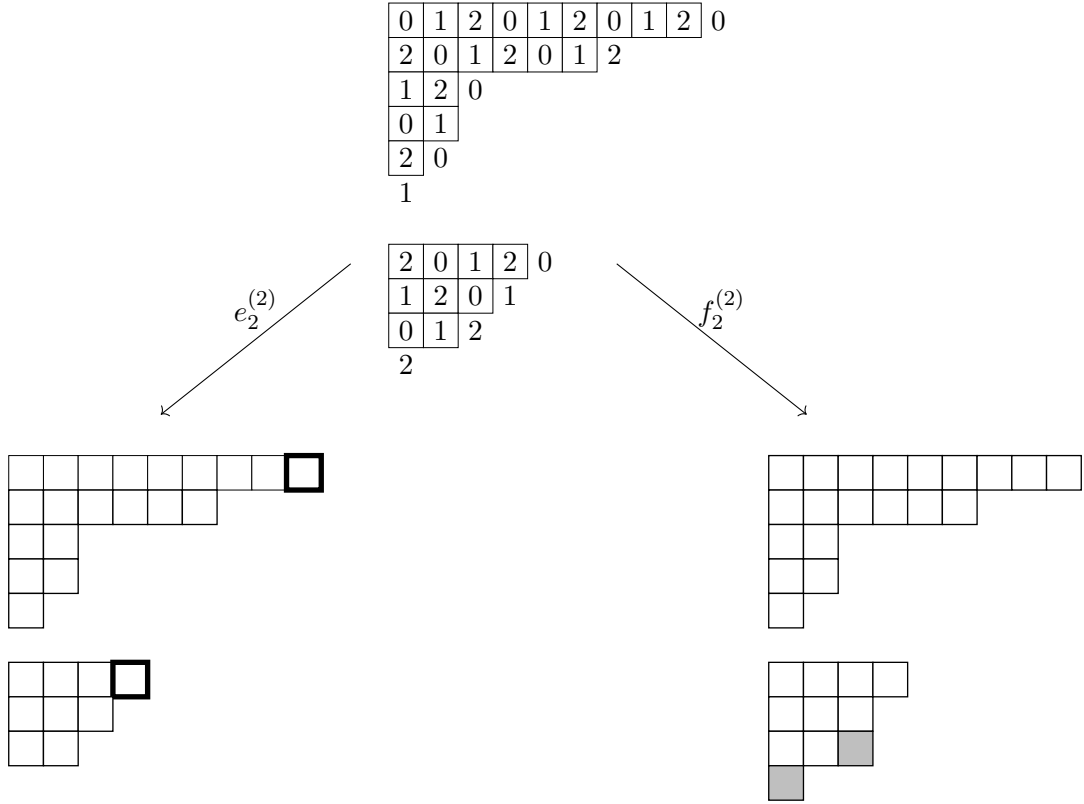
\end{ex}

For each $i\in\mathbb{Z}/e\mathbb{Z}$, there is at most one good $i$-node of $\lambda$, and hence at most $e$ good nodes of $\lambda$. It follows from~\cite[Theorem 0.5]{kbranch2} that the socle of the restriction of an irreducible $\mathscr{H}_n^{\Lambda}$-module $D_{\lambda}$ to an $\mathscr{H}_{n-1}^{\Lambda}$-module is a direct sum of at most $e$ indecomposable $\mathscr{H}_n^{\Lambda}$-summands.
Moreover, we also know from~\cite{kbranch2} that we can verify that the residue sequence of $\lambda\backslash\{A\}$ is distinct for each good node $A$ of $\lambda$, so that each summand $D_{\lambda\backslash\{A\}}$ belongs to a distinct block of $\mathscr{H}_n^{\Lambda}$. We generalise this result to ``divided powers'' as follows.

\begin{cor}\label{cor:modbr}
Let $i\in\mathbb{Z}/e\mathbb{Z}$ and $\lambda\in\mathscr{RP}_n^l$.
\begin{enumerate}
\item{If $r\leqslant{\nor_i(\lambda)}$, then $\soc\left(e_i^{(r)}D_{\lambda}\right)
	\cong
	D_{\lambda\downarrow_i^r}$.
}
\item{If $r\leqslant{\conor_i(\lambda)}$, then $\soc\left(
	f_i^{(r)}D_{\lambda}\right)\cong
	D_{\lambda\uparrow_i^r}$.
}
\end{enumerate}
\end{cor}

It follows that the modular branching rules for Specht modules of the cyclotomic Khovanov--Lauda--Rouquier algebras $\mathscr{H}_n^{\Lambda}$, together with the operators $\uparrow_i^r$ and $\downarrow_i^r$, provide a combinatorial algorithm for determining the labels of irreducible $\mathscr{H}_n^{\Lambda}$-modules.

\begin{prop}\label{prop:modbr}
Let $r\geqslant 0$ and $i\in\mathbb{Z}/e\mathbb{Z}$. 
If $D$ is an irreducible $\mathscr{H}_n^{\Lambda}$-module with $e_i^{(r)}D\cong{D_{\lambda}}$ for some $\lambda\in\mathscr{RP}_{n-r}^l$, then $D=D_{\lambda\uparrow_i^r}$.
\end{prop}

\begin{proof}
	Suppose that $D=D_{\mu}$ where $\mu\in\mathscr{RP}_n^l$, so that $e_i^{(r)}D=e_i^{(r)}D_{\mu}\cong{D_{\lambda}}$.
	We know that $r\leqslant\nor_i(\mu)$ since $e_i^{(r)}D\neq{0}$, then from the first part of \cref{cor:modbr} we have $\soc\left(e_i^{(r)}D_{\mu}\right)\cong{D_{\nu}}$ where $\nu=\mu\downarrow_i^r$.
	Since $e_i^{(r)}D_{\mu}\cong{D_{\lambda}}$, we have $\nu=\lambda$.
	Then, by (\ref{eq:arrows}), $\lambda\uparrow_i^r=\mu\downarrow_i^r\uparrow_i^r=\mu$, as required.
\end{proof}

Let $0\leqslant{r}\leqslant{\nor_i(\lambda)}$ with $e_i^{(r)}D_{\mu}\cong{D_{\lambda}}$ for some $\mu\in\mathscr{RP}_n^l$ and $\lambda\in\mathscr{RP}_{n-r}^l$. Then the normal $i$-nodes of $\mu$ and the conormal $i$-nodes of $\lambda$ coincide, and hence \[\soc\left(f_i^{(r)}\left(e_i^{(r)}D_{\mu}\right)\right)\cong{D_{\mu\downarrow_i^r\uparrow_i^r}}={D_{\mu}}.\]

For non-irreducible $\mathscr{H}_n^{\Lambda}$-modules, we can determine the labels of their composition factors by applying the same combinatorial algorithm using the following result.

\begin{cor}\label{irrcomp}
Let $r\geqslant{0}$ and $i\in\mathbb{Z}/e\mathbb{Z}$.
If $M$ is an $\mathscr{H}_n^{\Lambda}$-module with $e_i^{(r)}M\cong{D_{\mu}}$ for some $\mu\in\mathscr{RP}_{n-r}^l$, then one of the composition factors of $M$ is $D_{\mu\uparrow_i^r}$. Moreover, all of the other composition factors of $M$ are killed by $e_i^{(r)}$.
\end{cor}

\begin{ex}
Let $e=3$, $\kappa=(0,2)$ and $\lambda=((6),(1^3))$. We successively remove the maximum number of removable $i$-nodes (shaded below) from $\lambda$ as follows.
\[
\gyoung(01201!\gr2,,!\wh2,1,0)\xrightarrow{e_2}
\gyoung(0120!\gr1,,!\wh2,1,0)\xrightarrow{e_1}
\gyoung(012!\gr0,,!\wh2,1,!\gr0)\xrightarrow{e_0^{(2)}}
\young(01!\gr2,,!\wh2,1)\xrightarrow{e_2}
\young(0!\gr1,,!\wh2,!\gr1)\xrightarrow{e_1^{(2)}}
\young(0,,!\gr2)\xrightarrow{e_2}
\gyoung(!\gr0,,:\emset)\xrightarrow{e_0}
\gyoung(:\emset,,:\emset).
\]
Hence $e_0e_2e_1^{(2)}e_2e_0^{(2)}e_1e_2S_{\lambda}\cong S_{(\varnothing,\varnothing)}$, which we know is irreducible, and moreover, $(\varnothing,\varnothing)\uparrow_0\uparrow_2\uparrow_1^2\uparrow_2\uparrow_0^2\uparrow_1\uparrow_2=((6,2,1),\varnothing)$. It thus follows from \cref{irrcomp} that $D_{((6,2,1),\varnothing)}$ is a composition factor of $S_{\lambda}$.
\end{ex}

\subsection{Specht modules labelled by hook bipartitions}

We fix $l=2$ from now on and recall that $e\in\{3,4,\dots\}$. A \emph{hook bipartition of $n$} is defined to be a bipartition of the form $((n-m),(1^m))$ for some $m\in\{0,\dots,n\}$. We refer to the first component of a hook bipartition as its \emph{arm} and to its second component as its \emph{leg}. We call the node $(1,n-m,1)$ lying at the end of its arm its \emph{hand node}, and the node $(m,1,2)$ lying at the end of its leg its \emph{foot node}.

Let $\ttt$ be the standard $((n-m),(1^m))$-tableau with entries $a_1,\dots,a_m\in\{1,\dots,n\}$ lying in its leg, and recall from~\cite[\S{5.2}]{Sutton17I} that we define $v(a_1,\dots,a_m):=v_{\ttt}$ to be the corresponding standard basis element of $S_{((n-m),(1^m))}$.
We note that since $\ttt$ is completely defined by the strictly increasing entries $a_1,\dots,a_m$ that lie in a single column,  the corresponding vector $v_{\ttt}$ is independent of the choice of a reduced expression of $w_{\ttt}$ (which is generally not the case).

We remind the reader of some of the Specht module homomorphisms that were introduced in~\cite[Proposition 5.6]{Sutton17I}, which we will require later on.

\begin{prop}
\label{prop:homs}
We have the following non-zero homomorphisms of Specht modules.
\begin{itemize}
\item{
If $n\equiv{\kappa_2-\kappa_1+1}\Mod{e}$ and $0\leqslant{m}\leqslant{n-1}$, we have
	\begin{align*}
	\gamma_m:
	S_{((n-m),(1^m))}
	&\longrightarrow
	S_{((n-m-1),(1^{m+1}))},\ 
	\gamma_m(z_{((n-m),(1^m))})=v(1,\dots,m,n).
	\end{align*}
}
\item{
If $\kappa_2\equiv{\kappa_1-1}\Mod{e}$ and $1\leqslant{m}\leqslant{n-1}$, we have
	\begin{align*}
	\chi_m:
	S_{((n-m,1^m),\varnothing)}
	&\longrightarrow
	S_{((n-m),(1^m))},\ 
	\chi_m(z_{((n-m,1^m),\varnothing)})=v(2,3,\dots,m+1).
	\end{align*}
}
\item{
If $\kappa_2\equiv{\kappa_1-1}\Mod{e}$, $n\equiv{0}\Mod{e}$ and $1\leqslant{m}\leqslant{n-1}$, we have
	\begin{align*}
	\phi_m:
	S_{((n-m+1,1^{m-1}),\varnothing)}
	&\longrightarrow
	S_{((n-m),(1^m))},\ 
	\phi_m(z_{((n-m+1,1^{m-1}),\varnothing)})=v(2,3,\dots,m,n).
	\end{align*}
}
\end{itemize}
\end{prop}

\begin{comment}

The standard basis elements of the kernels and images of the above Specht modules homomorphisms are as follows.

\begin{lem}\cite[Lemma 5.7]{Sutton17I}
\label{imker}
\begin{enumerate}
\item{
If $n\equiv{\kappa_2-\kappa_1+1}\Mod{e}$, then
	\begin{enumerate}
	\item{
	$\im(\gamma_m)
	=\spn\left\{
	v_{\ttt}\ 
	\mid\ 
	\ttt\in\std\left((n-m-1),(1^{m+1})\right),
	\ttt(m+1,1,2)=n
	\right\}$;
	}
	\item{
	$\ker(\gamma_m)
	=\spn\left\{
	v_{\ttt}\ 
	\mid\ 
	\ttt\in\std\left((n-m),(1^m)\right),
	\ttt(m,1,2)=n
	\right\}$.
	}
	\end{enumerate}
}
\item{If $\kappa_2\equiv\kappa_1-1\Mod{e}$, then
	\begin{enumerate}
		\item{
			$\im(\chi_m)=\spn\left\{v_{\ttt}\ \mid\
			\ttt\in\std\left((n-m),(1^m)\right),\ttt(1,1,1)=1\right\}$; }
			\item{ $\ker(\chi_m)=0$; 
	}
	\end{enumerate}
}
\item{
If $\kappa_2\equiv\kappa_1-1\Mod{e}$ and $n\equiv{\kappa_2-\kappa_1+1}\Mod{e}$, then
		\begin{enumerate}
		\item If $m<n-1$, then
		\[\im(\phi_m)=\spn\left\{v_{\ttt}\mid
		\ttt\in\std\left((n-m),(1^m)\right),\ttt(1,1,1)=1,\ttt(m,1,2)=n\right\}.\]
		\item If $m=n-1$, then
		\[\im(\phi_m)=\spn\{
		v_{\ttt}\mid
		\ttt\in\std\left((1),(1^{n-1})\right),\ttt(1,1,1)=1\}.\]
		\end{enumerate}			
}
\end{enumerate}
\end{lem}
\end{comment}

\section{One-dimensional Specht modules}\label{sec:onedims}

We determine the labels of the irreducible $\mathscr{H}_n^{\Lambda}$-modules that are isomorphic to the one-dimensional Specht modules, namely $S_{((n),\varnothing)}$ and $S_{(\varnothing,(1^n))}$, which arise as composition factors of $S_{((n-m),(1^m))}$. We note that we work solely with \emph{ungraded} cyclotomic Khovanov--Lauda--Rouquier modules up to and including \Cref{sec:UDN}.
We let $l$ be the residue of $\kappa_2-\kappa_1$ modulo $e$ throughout, so that $l\in\{0,\dots,e-1\}$.

We know that $S_{((n),\varnothing)}=\left\{z_{((n),\varnothing)}\right\}$ and $S_{(\varnothing,(1^n))}=\left\{z_{(\varnothing,(1^n))}\right\}$ are both one-dimensional $\mathscr{H}_n^{\Lambda}$-modules, and hence are both irreducible. In fact, $S_{((n),\varnothing)}={D_{((n),\varnothing)}}$. We now introduce a \emph{$\sgn$-functor} to determine the bipartition $\mu\in\mathscr{RP}_n^2$ such that $S_{(\varnothing,(1^n))}\cong{D_{\mu}}$ as ungradwed $\mathscr{H}_n^{\Lambda}$-modules.

For $1\leqslant{r}\leqslant{n}$, $S_{(\varnothing,(1^r))}$ has only one removable node, namely $(r,1,2)$, that satisfies $\res(r,1,2)=\kappa_2+1-r\Mod{e}$. Thus the only restriction functor which acts non-trivially on $S_{(\varnothing,(1^r))}$ is $e_{\kappa_2+1-r}:\mathscr{H}_r^{\Lambda}\md\longrightarrow{\mathscr{H}_{r-1}^{\Lambda}\md}$, where $e_{\kappa_2+1-r}S_{(\varnothing,(1^r))}\cong{S_{(\varnothing,(1^{r-1}))}}$. For $r\geqslant{0}$, we now define the \emph{$\sgn$-restriction functor} to be the composition of restriction functors
\[
e_{\sgn}:=e_{\kappa_2}\circ{e_{\kappa_2-1}}\circ\dots\circ{e_{\kappa_2+1-r}}
:\mathscr{H}_{r}^{\Lambda}\md\longrightarrow{\mathscr{H}_{0}^{\Lambda}\md},
\]
with the property that
	\[
	e_{\sgn}S_{(\varnothing,(1^r))}\cong{S_{(\varnothing,\varnothing)}}.
	\]
With $r=n$, we observe that $e_{\sgn}$ is the only composition of $n$ $i$-restriction functors which acts non-trivially on $S_{(\varnothing,(1^n))}$. Analogously, we now define the \emph{$\sgn$-induction functor} to be
\[
f_{\sgn}:=f_{\kappa_2+1-r}\circ{f_{\kappa_2+2-r}}\circ\dots\circ{f_{\kappa_2}}
:\mathscr{H}_{0}^{\Lambda}\md\longrightarrow{\mathscr{H}_{r}^{\Lambda}\md}
\]
for $r\geqslant 0$. The $\sgn$-induction functor acts non-trivially on $S_{(\varnothing,\varnothing)}$; we now determine the socle of $f_{\sgn}S_{(\varnothing,\varnothing)}$.

\begin{defn}
For each $a\in\mathbb{N}\cup\{0\}$, we define the following weakly decreasing sequence of $e-1$ integers that sum to $a$ and that differ by at most one
\begin{align*}
\{a\}
:=
\left\lfloor\frac{a+e-2}{e-1}\right\rfloor,
\left\lfloor\frac{a+e-3}{e-1}\right\rfloor,
\dots,
\left\lfloor\frac{a}{e-1}\right\rfloor.
\end{align*}
\end{defn}
We now give an explicit description of the regular bipartition that labels the irreducible $\mathscr{H}_n^{\Lambda}$-module that is isomorphic to $S_{(\varnothing,(1^n))}$.

\begin{defn}
Let
$
{(\varnothing,(1^n))^R}:=\begin{cases}
{(\varnothing,(1^n))}&\text{if $n<l$,}\\
{\left(\left(\{n-l\}\right),\left(1^l\right)\right)}&\text{if $n\geqslant{l}$.}
\end{cases}
$
\end{defn}

\begin{lem}\label{simpsign}
Let $n\in\mathbb{N}\cup\{0\}$. Then $S_{(\varnothing,(1^n))}\cong D_{(\varnothing,(1^n))^R}$ as ungraded $\mathscr{H}_n^{\Lambda}$-modules.
\end{lem}

\begin{proof}
Let $1\leqslant{r}\leqslant{n}$, and suppose that $S_{(\varnothing,(1^n))}\cong{D_{\mu}}$ for some $\mu\in\mathscr{RP}_n^2$. It follows from (\ref{eq:indres}) that
	\[
	\res^{\mathscr{H}_n^{\Lambda}}_{\mathscr{H}_0^{\Lambda}}
	S_{(\varnothing,(1^n))}
	\cong
	e_{\sgn}
	S_{(\varnothing,(1^n))}.
	\]
For any $r>1$, there is only one removable $(\kappa_2+1-r)$-node of $[(\varnothing,(1^r))]$, so that $\epsilon_{\kappa_2+1-r}\left(S_{(\varnothing,(1^r))}\right)=1$.
It thus follows from~\cref{Spechtmax} that
	\[
	e_{\sgn}
	S_{(\varnothing,(1^n))}
	\cong
	S_{(\varnothing,(1^n))^{\triangledown{e_{\sgn}}}}
	=
	S_{{(\varnothing,(1^n))}
	^{{\triangledown{(\kappa_2+1-n)}}
	{\triangledown{(\kappa_2+2-n)}}
	{\cdots}
	{\triangledown{{\kappa_2}}}}}
	\cong
	S_{(\varnothing,\varnothing)}.
	\]
Define 
	$(\varnothing,\varnothing)\uparrow_{\sgn}^n
	:= (\varnothing,\varnothing)
	\uparrow_{\kappa_2}	\uparrow_{\kappa_2-1} \dots	\uparrow_{\kappa_2+1-n}$.
Since $S_{(\varnothing,(1^n))}$ is irreducible, we know from \cref{prop:modbr} that
	$S_{(\varnothing,(1^n))} \cong
	S_{(\varnothing,\varnothing)\uparrow_{\sgn}^n}
	= D_{(\varnothing,\varnothing)\uparrow_{\sgn}^n}$.
To calculate $(\varnothing,\varnothing)\uparrow_{\sgn}^n$, we successively add the highest conormal node of $e$-residue $\kappa_2,\kappa_2-1,\dots,\kappa_2+1-n$, respectively, to $[(\varnothing,\varnothing)]$.

Firstly, we successively add the highest $l$ conormal nodes of $e$-residue $\kappa_2,\kappa_2-1,\dots,\kappa_2-l+1$, respectively, to $[(\varnothing,\varnothing)]$. Since $\kappa_1={\kappa_2-l}\Mod{e}$, it is easy to see that $(\varnothing,(1^i))$ has $(\kappa_2-i)$-signature $+$, corresponding to node $(i+1,1,2)$ for each $i\in\{0,\dots,l-1\}$. Hence 
\[
(\varnothing,\varnothing)
\uparrow_{\kappa_2}
\uparrow_{\kappa_2-1}
\dots
\uparrow_{\kappa_2-l+1}
={(\varnothing,(1^l))}.
\]
If $n\leqslant{l}$, then we are done. Instead suppose that $n>l$.
We now successively add the highest $e$ conormal nodes to $[(\varnothing,(1^l))]$ of $e$-residue $\kappa_1,\kappa_1-1,\dots,\kappa_1+1$, respectively. Notice that $((1^i),(1^l))$ has $(\kappa_1-i)$-signature $+$ for each $i\in\{0,\dots,e-1\}$, corresponding to node $(i+1,1,1)$, except in the following cases.
\begin{itemize}
\item{
The $\kappa_1$-signature of $(\varnothing,(1^l))$ is $++$, corresponding to the nodes $(1,1,1)$ and $(l+1,1,2)$, respectively. Hence $(\varnothing,(1^l))\uparrow_{\kappa_1}=((1),(1^l))$.
}
\item{
Let $l>0$. Then the $(\kappa_1+l+1)$-signature of $((1^{e-l-1}),(1^l))$ is $++$, corresponding to the nodes $(e-l,1,1)$ and $(1,2,2)$, respectively. Hence $((1^{e-l-1}),(1^l))\uparrow_{\kappa_1+l+1}=((1^{e-l}),(1^l))$.
}
\item{
If $l>0$, then the $(\kappa_1+1)$-signature of $((1^{e-1}),(1^l))$ is $++-$, corresponding to the nodes $(1,2,1)$, $(e,1,1)$ and $(l,1,2)$, respectively. If $l=0$, then the $(\kappa_1+1)$-signature of $((1^{e-1}),\varnothing)$ is $++$, corresponding to the nodes $(1,2,1)$ and $(e,1,1)$, respectively. Hence $((1^{e-1}),(1^l))\uparrow_{\kappa_1+1}=((2,1^{e-2}),(1^l))$.
}
\end{itemize}
It thus follows that
\[
(\varnothing,(1^l))
\uparrow_{\kappa_1}
\uparrow_{\kappa_1-1}
\dots
\uparrow_{\kappa_1+1}
=((2,1^{e-2}),(1^l)),
\]
and so the first component of $(\varnothing,\varnothing)\uparrow_{\sgn}^n$ has $e-1$ non-empty rows.

Finally, we successively add the remaining nodes to the first component of $[((2,1^{e-2}),(1^l))]$, down each column from left to right.
Observe that there are $n-l-r+1$ nodes in 
\[
[(\varnothing,\varnothing)\uparrow_{\sgn}^n]
\backslash{\{(1,1,1),\dots,(r-1,1,1)\}\cup\{(1,1,2),\dots,(l,1,2)\}}
\]
for all $r\in\{1,\dots,e-1\}$. Since there are $e-1$ non-empty rows in the first component of $[((2,1^{e-2}),(1^l))]$, there are also $e-1$ non-empty rows in the first component of $\mu$, and moreover, we observe that there are $\left\lfloor\tfrac{n-l-r+e-1}{e-1}\right\rfloor$ nodes in the $r$th row of the first component of $[(\varnothing,\varnothing)\uparrow_{\sgn}^n]$.
\end{proof}

\section{Labelling the composition factors of $S_{((n-m),(1^m))}$}\label{sec:labelsgen}

In the preceding paper~\cite[\S{6}]{Sutton17I}, the composition factors of $S_{((n-m),(1^m))}$ were constructed as quotients either of the images or of the kernels of the Specht module homomorphisms given in \cref{prop:homs} --- both of which do not depend on the characteristic of $\mathbb{F}$.
Since each of these quotients is isomorphic (up to a grading shift) to a particular head of a Specht module, we now determine the regular bipartitions that label these irreducible $\mathscr{H}_n^{\Lambda}$-modules.
Recall that $l$ is the residue of $\kappa_2-\kappa_1$ modulo $e$.

\subsection{Labelling the composition factors of $S_{((n-m),(1^m))}$ with $\kappa_2\not\equiv{\kappa_1-1}\Mod{e}$}\label{sec:labels1}

We fix $\kappa_2\not\equiv{\kappa_1-1}\Mod{e}$ throughout this subsection.

When $n\not\equiv{l+1}\Mod{e}$, we recall from~\cite[Theorem 6.8]{Sutton17I} that $S_{((n-m),(1^m))}$ is an irreducible $\mathscr{H}_n^{\Lambda}$-module, that is, $S_{((n-m),(1^m))}\cong{D_{\lambda_m}}$ for some regular bipartition $\lambda_m\in\mathscr{RP}_n^2$.

\begin{defn}\label{def:lnot-1,nnot1}
Let $\kappa_2\not\equiv{\kappa_1-1}\Mod{e}$ and $n\not\equiv{l+1}\Mod{e}$. For $0\leqslant{m}<n$, we define
\begin{align*}
\mu_{n,m}:=
\begin{cases}
((n-m),(1^m))
&\text{ if $0\leqslant{m}<{l+1}$,}\\
((n-m,\{m-l-1\}),(1^{l+1}))
&\text{ if $l+1\leqslant{m}<{n-\frac{n}{e}}$,}\\
((\{m-l\},n-m-1),(1^{l+1}))
&\text{ if $n-\frac{n}{e}\leqslant{m}<n$.}
\end{cases}
\end{align*}
\end{defn}

In fact, we claim that $\lambda_m=\mu_{n,m}$ for all $m\in\{0,\dots,n-1\}$.

When $n\equiv{l+1}\Mod{e}$ and $1\leqslant{m}<n$, we recall from~\cite[Corollary 6.11]{Sutton17I} that $S_{((n-m),(1^m))}$ has two composition factors, namely $\im(\gamma_{m-1})$ and $\im(\gamma_m)$. Thus $\im(\gamma_{m-1})\cong{D_{\lambda_m}}$ and $\im(\gamma_m)\cong{D_{\mu_m}}$ for some regular bipartitions $\lambda_m,\mu_m\in\mathscr{RP}_n^2$.

\begin{defn}
Let $\kappa_2\not\equiv{\kappa_1-1}\Mod{e}$ and $n\equiv{l+1}\Mod{e}$. For $0\leqslant{m}<n$, we define
\begin{align*}
\mu_{n,m}:=
\begin{cases}
((n-m),(1^{m}))
&\text{if }0\leqslant{m}<{l+1},\\
\left(\left(n-m,\{m-l-1\}\right),\left(1^{l+1}\right)\right)
&\text{if }l+1\leqslant{m}<{n-\frac{n}{e}},\\
\left(\left(\{m-l+1\},n-m-2\right),\left(1^{l+1}\right)\right)
&\text{if }n-\frac{n}{e}\leqslant{m}\leqslant{n-2},\\
\left(\left(\{n-l\}\right),(1^l)\right)
&\text{if $m=n-1$}.
\end{cases}
\end{align*}
\end{defn}

Notice that $\mu_{n,m-1}$ and $\mu_{n,m}$ are distinct. We claim that the two labels $\lambda_m,\mu_m$ of the composition factors of $S_{((n-m),(1^m))}$ as heads of some Specht modules are, in fact, $\mu_{n,m-1}$ and $\mu_{n,m}$, respectively, and hence that the corresponding composition factors are non-isomorphic.

We require the following combinatorial result in order to confirm our claims.

\begin{lem}\label{lem:labels1}
Let $\kappa_2\not\equiv{\kappa_1-1}\Mod{e}$ and $0\leqslant{m}<n$.
\begin{enumerate}
\item{ If $n\equiv{l}\Mod{e}$, then \begin{equation}\label{leg1}
\mu_{n,m}\uparrow_{\kappa_2-m}=\mu_{n+1,m}.
\end{equation} }
\item{ If $n\not\equiv{l}\Mod{e}$, then \begin{equation}\label{leg2}
\mu_{n,m}\uparrow_{\kappa_2-m}=\mu_{n+1,m+1}.
\end{equation} }
\end{enumerate}
\end{lem}

\begin{proof}
\begin{enumerate}[label=(\roman*)]
\item{
Let $1\leqslant{m}<l+1$. 
Observe that $((n-m),(1^m))$ has addable $(\kappa_2-m)$-node $(m+1,1,2)$, as well as $(1,n-m+1,1)$ if $n\equiv{l}\Mod{e}$, and has removable $(\kappa_2-m)$-node $(1,n-m,1)$ if $n\equiv{l+1}\Mod{e}$. We note that the addable nodes $(2,1,1)$ and $(1,2,2)$ of $((n-m),(1^m))$ cannot have residue $\kappa_2-m$ since $l<e-2$.

If $n\equiv{l}\Mod{e}$ then $((n-m),(1^m))$ has $(\kappa_2-m)$-signature $++$, corresponding to the conormal nodes $(1,n-m+1,1)$ and $(m+1,1,2)$. Adding the higher of these conormal nodes, we have
\[
\mu_{n,m}\uparrow_{\kappa_2-m}
=((n-m),(1^m))\uparrow_{\kappa_2-m}
=((n-m+1),(1^m))
=\mu_{n+1,m}.
\]
Now suppose that $n\equiv{l+1}\Mod{e}$. Then $((n-m),(1^m))$ has $(\kappa_2-m)$-signature $-+$, and if $n-l\not\equiv{0,1}\Mod{e}$ then $((n-m),(1^m))$ has $(\kappa_2-m)$-signature $+$. The conormal node in each sequence is $(m+1,1,2)$, whereby adding this node gives us
\[
\mu_{n,m}\uparrow_{\kappa_2-m}
=((n-m),(1^m))\uparrow_{\kappa_2-m}
=((n-m),(1^{m+1}))
=\mu_{n+1,m+1}.
\]
}
\item{
Let $l+1\leqslant{m}<{n-\tfrac{n}{e}}$.
Observe that $((n-m,\{m-l-1\}),(1^{l+1}))$ has the following addable and removable $(\kappa_2-m)$-nodes
\begin{itemize}
	\item {addable node $(1,n-m+1,1)$ if $n\equiv{l}\Mod{e}$,}
	\item {removable node $(1,n-m,1)$ if $n\equiv{l+1}\Mod{e}$,}
	\item {addable node at the end of the $\lfloor(m+e-l-2)/(e-1)\rfloor$th column in the first component,}
	\item {addable node $(e+1,1,1)$ and removable node $(l+1,1,2)$ if $m\equiv{l}\Mod{e}$,}
	\item {addable node $(1,2,2)$ if $m\equiv{-1}\Mod{e}$,}
	\item {addable node $(l+2,1,2)$ if $m\equiv{l+1}\Mod{e}$.}
\end{itemize}

First suppose that $n\equiv{l}\Mod{e}$. Then $((n-m,\{m-l-1\}),(1^{l+1}))$ has $(\kappa_2-m)$-signature
\begin{itemize}
	\item {$+++-$ if $m\equiv{l}\Mod{e}$,}
	\item {$++++$ if $m\equiv{-1}\Mod{e}$ and $l=e-2$,}
	\item {$+++$ if $m\equiv{-1}\Mod{e}$ and $l\neq{e-2}$ or $m\not\equiv{-1}\Mod{e}$ and $m\equiv{l+1}\Mod{e}$,}
	\item {$++$ for all other cases.}
\end{itemize}
Adding the highest conormal $(\kappa_2-m)$-node in these sequences, $(1,n-m+1,1)$, we have
\begin{align*}
\mu_{n,m}\uparrow_{\kappa_2-m}
&=((n-m,\{m-l-1\}),(1^{l+1}))\uparrow_{\kappa_2-m}\\
&=((n-m+1,\{m-l-1\}),(1^{l+1}))\\
&=\mu_{n+1,m}.
\end{align*}
We now suppose that $n\not\equiv l\Mod{e}$.
If $n\equiv{l+1}\Mod{e}$, then $((n-m,\{m-l-1\}),(1^{l+1}))$ has $(\kappa_2-m)$-signature
\begin{itemize}
	\item {$-++-$ if $m\equiv{l}\Mod{e}$,}
	\item {$-+++$ if $m\equiv{-1}\Mod{e}$ and $l=e-2$,}
	\item {$-++$ if $m\equiv{-1}\Mod{e}$ and $l\neq{e-2}$ or $m\not\equiv{-1}\Mod{e}$ and $m\equiv{l+1}\Mod{e}$,}
	\item {$-+$ for all other cases.}
\end{itemize}
If $n-l\not\equiv{0,1}\Mod{e}$, then $((n-m,\{m-l-1\}),(1^{l+1}))$ has $(\kappa_2-m)$-signature
\begin{itemize}
	\item {$++-$ if $m\equiv{l}\Mod{e}$,}
	\item {$+++$ if $m\equiv{-1}\Mod{e}$ and $l=e-2$,}
	\item {$++$ if $m\equiv{-1}\Mod{e}$ and $l\neq{e-2}$ or $m\not\equiv{-1}\Mod{e}$ and $m\equiv{l+1}\Mod{e}$,}
	\item {$+$ for all other cases.}
\end{itemize}
Thus, for $n\not\equiv{l}\Mod{e}$, we observe that the highest conormal $(\kappa_2-m)$-node in each $(\kappa_2-m)$-signature of $((n-m,\{m-l-1\}),(1^{l+1}))$ is the addable node lying at the bottom of its $\lfloor(m+e-l-2)/(e-1)\rfloor$th column in the first component. Adding this node, we have
\begin{align*}
\mu_{n,m}\uparrow_{\kappa_2-m}
&=((n-m,\{m-l-1\}),(1^{l+1}))\uparrow_{\kappa_2-m}\\
&=((n-m,\{m-l\}),(1^{l+1}))\\
&=\mu_{n+1,m+1}.
\end{align*}
}
\item{
Let $m\geqslant{n-\tfrac{n}{e}}$.
Firstly, suppose that $n\not\equiv{l+1}\Mod{e}$.
If $n\equiv{l}\Mod{e}$, then we find that $((\{m-l\},n-m-1),(1^{l+1}))$ has $(\kappa_2-m)$-signature
\begin{itemize}
	\item {$+++-$ if $m\equiv{l}\Mod{e}$,}
	\item {$++++$ if $m\equiv{-1}\Mod{e}$ and $l=e-2$,}
	\item {$+++$ if $m\equiv{-1}\Mod{e}$ and $l\neq{e-2}$ or $m\not\equiv{-1}\Mod{e}$ and $m\equiv{l+1}\Mod{e}$,}
	\item {$++$ for all other cases.}
\end{itemize}
For $n-l\not\equiv{0,1}\Mod{e}$, $(e,n-m,1)$ is no longer an addable $(\kappa_2-m)$-node. Upon discounting this node, we observe that $((\{m-l\},n-m-1),(1^{l+1}))$ has $(\kappa_2-m)$-signatures $++-$, $+++$, $++$ and $+$ corresponding to the above cases, respectively. Thus, for $n\not\equiv{l+1}\Mod{e}$, the highest conormal $(\kappa_2-m)$-node in each $(\kappa_2-m)$-signature of $((\{m-l\},n-m-1),(1^{l+1}))$ is the addable node at the bottom of the $\lfloor(m+e-l-2)/(e-1)\rfloor$th column in the first component.
Hence
\begin{align*}
\mu_{n,m}\uparrow_{\kappa_2-m}
&=((\{m-l\},n-m-1),(1^{l+1}))\uparrow_{\kappa_2-m}\\
&=((\{m-l+1\},n-m-1),(1^{l+1}))\\
&=
\begin{cases}
\mu_{n+1,m}
&\text{if $n\equiv{l}\Mod{e}$,}\\
\mu_{n+1,m+1}
&\text{if $n-l\not\equiv{0,1}\Mod{e}$.}
\end{cases}
\end{align*}

Secondly, suppose that $n\equiv{l+1}\Mod{e}$.
Then we find that $((\{m-l+1\},n-m-2),(1^{l+1}))$ has $(\kappa_2-m)$-signature
\begin{itemize}
	\item {$-++-$ if $m\equiv{l}\Mod{e}$,}
	\item {$-+++$ if $m\equiv{-1}\Mod{e}$ and $l=e-2$,}
	\item {$-++$ if $m\equiv{-1}\Mod{e}$ and $l\neq{e-2}$ or $m\not\equiv{-1}\Mod{e}$ and $m\equiv{l+1}\Mod{e}$,}
	\item {$-+$ for all other cases.}
\end{itemize}
The highest conormal $(\kappa_2-m)$-node in each sequence is $(e,n-m-1,1)$, and adding this node we have
\begin{align*}
\mu_{n,m}\uparrow_{\kappa_2-m}
&=((\{m-l+1\},n-m-2),(1^{l+1}))\uparrow_{\kappa_2-m}\\
&=((\{m-l+1\},n-m-1),(1^{l+1}))\\
&=\mu_{n+1,m+1}.
\end{align*}
\qedhere
}
\end{enumerate}
\end{proof}

\begin{thm}\label{factors}
Suppose that $\kappa_2\not\equiv{\kappa_1-1}\Mod{e}$ and $n\not\equiv{l+1}\Mod{e}$.
Then $S_{((n-m),(1^m))}\cong{D_{\mu_{n,m}}}$ as ungraded $\mathscr{H}_n^{\Lambda}$-modules for all $m\in\{0,\dots,n-1\}$.
\end{thm}

\begin{proof}We proceed by induction on $m$.
\begin{enumerate}
\item{Suppose that $n-l\not\equiv{2}\Mod{e}$.
First observe that $e_{\kappa_2}S_{((n-1),(1))}\cong{S_{((n-1),\varnothing)}}$, and moreover, by applying \Cref{leg2} we have $((n-1),\varnothing)\uparrow_{\kappa_2}=((n-1),(1))$. It thus follows from \cref{irrcomp} that $S_{((n-1),(1))}\cong{D_{((n-1),(1))}}$.

Assuming that $S_{((n-m),(1^{m-1}))}\cong{D_{\mu_{n-1,m-1}}}$ for some $m>1$, then 
\[
e_{\kappa_2+1-m}S_{((n-m),(1^m))}\cong{S_{((n-m),(1^{m-1}))}}\cong{D_{\mu_{n-1,m-1}}}.
\]
Since $S_{((n-m),(1^m))}$ is an irreducible $\mathscr{H}_n^{\Lambda}$-module, we can apply \cref{irrcomp} and \Cref{leg2} to obtain
\begin{align*}
S_{((n-m),(1^m))}\cong{D_{\mu_{n-1,m-1}\uparrow_{\kappa_2-m+1}}}={D_{\mu_{n,m}}}.
\end{align*}
}
\item{
Suppose that $n-l\equiv{2}\Mod{e}$.
We obtain the regular bipartition that labels the irreducible module which $S_{((n-m),(1^m))}$ is isomorphic to, up to a grading shift, as follows. We first restrict $S_{((n-m),(1^m))}$ to an irreducible $\mathscr{H}_{n-2}^{\Lambda}$-module, say $D_{\mu}$ for some $\mu\in\mathscr{RP}_{n-2}^2$, by removing both the hand and foot node of residue $\kappa_2+1-m$ modulo $e$ from the hook bipartition $((n-m),(1^m))$, and then inducing $D_{\mu}$ to an irreducible $\mathscr{H}_{n}^{\Lambda}$-module by adding the two highest conormal $(\kappa_2+1-m)$-nodes to $\mu$.
We have
$e_{\kappa_2}^{(2)}S_{((n-1),(1))}
\cong{S_{((n-2),\varnothing)}}$.
By both \Cref{leg1,leg2},
$
((n-2),\varnothing)\uparrow_{\kappa_2}^2
=((n-1),\varnothing)\uparrow_{\kappa_2}
=((n-1),(1))$.
Hence $S_{((n-1),(1))}\cong{D_{((n-1),(1))}}$ by \cref{irrcomp}.

Assuming that ${S_{((n-m-1),(1^{m-1}))}}\cong{D_{\mu_{n-2,m-1}}}$ for some $m>1$, then
\[
{e_{\kappa_2-m+1}^{(2)}S_{((n-m),(1^m))}}\cong{S_{((n-m-1),(1^{m-1}))}}\cong{D_{\mu_{n-2,m-1}}}.
\]
Since $S_{((n-m),(1^m))}$ is an irreducible $\mathscr{H}_n^{\Lambda}$-module, we apply \cref{irrcomp} to obtain
\begin{align*}
{S_{((n-m),(1^m))}}\cong
{D_{\mu_{n-2,m-1}\uparrow_{\kappa_2-m+1}^2}}
&={D_{\mu_{n-1,m-1}\uparrow_{\kappa_2-m+1}}}
&&\text{(by \Cref{leg1})}\\
&={D_{\mu_{n,m}}}
&&\text{(by \Cref{leg2}).}
\end{align*}
\qedhere
}
\end{enumerate}
\end{proof}

We now use this result to give an explicit description of the composition factors of Specht modules labelled by hook bipartitions in the following case.

\begin{thm}\label{thm:labels case 2}
Suppose that $\kappa_2\not\equiv{\kappa_1-1}\Mod{e}$ and $n\equiv{l+1}\Mod{e}$. Then the composition factors of $S_{((n-m),(1^m))}$ are $D_{\mu_{n,m-1}}$ and $D_{\mu_{n,m}}$ for all $m\in\{1,\dots,n-1\}$. Moreover, ${D_{\mu_{n,m}}}\cong{\im(\gamma_m)}$ as ungraded $\mathscr{H}_n^{\Lambda}$-modules.
\end{thm}

\begin{proof}
We obtain the regular bipartitions that label the two composition factors of $S_{((n-m),(1^m))}$ as heads of some Specht modules by first restricting this Specht module to an irreducible $\mathscr{H}_{n-1}^{\Lambda}$-module, say $D_{\mu}$ for some $\mu\in\mathscr{RP}_{n-1}^2$, by either 1) removing the foot node of residue $\kappa_2+1-m$ from $((n-m),(1^m))$ or 2) by removing the hand node of residue $\kappa_2-m$ from $((n-m),(1^m))$. We then induce $D_{\mu}$ to an irreducible $\mathscr{H}_n^{\Lambda}$-module by adding the highest conormal node of residue $\kappa_2+1-m$ or $\kappa_2-m$, respectively, to $\mu$.
\begin{enumerate}
\item{
By removing the foot node of $[((n-m),(1^m))]$, we obtain
\[
e_{\kappa_2-m+1}S_{((n-m),(1^m))}\cong{S_{((n-m),(1^{m-1}))}}
\cong D_{\mu_{n-1,m-1}}
\qquad(\text{by \cref{factors}}).
\]
We now observe from \Cref{leg1} that $\mu_{{n-1,m-1}\uparrow_{\kappa_2+1-m}}=\mu_{n,m-1}$. It thus follows that $D_{\mu_{n,m-1}}$ is a composition factor of $S_{((n-m),(1^m))}$ by \cref{irrcomp}.
}
\item{
By removing the hand node of $[((n-m),(1^m))]$, we obtain \[e_{\kappa_2-m}S_{((n-m),(1^m))}\cong{S_{((n-m-1),(1^m))}}
\cong D_{\mu_{n-1,m}}
\qquad\qquad(\text{by \cref{factors}}).\]
We now observe from \Cref{leg1} that $\mu_{n-1,m}\uparrow_{\kappa_2-m}=\mu_{n,m}$. Then, by \cref{irrcomp}, $D_{\mu_{n,m}}$ is a composition factor of $S_{((n-m),(1^m))}$.
}
\end{enumerate}

Furthermore, we know from~\cite[Corollary 6.11]{Sutton17I} that $\operatorname{im}(\gamma_{m-1})$ and $\operatorname{im}(\gamma_m)$ are in bijection with $D_{\mu_{n,m-1}}$ and $D_{\mu_{n,m}}$, up to isomorphism and grading shift. We notice that $\operatorname{im}(\gamma_m)$ is a composition factor of both $S_{((n-m),(1^m))}$ and $S_{((n-m-1),(1^{m+1}))}$, and hence must be isomorphic to $D_{\mu_{n,m}}$, as required.
\end{proof}

\subsection{Labelling the composition factors of $S_{((n-m),(1^m))}$ with $\kappa_2\equiv{\kappa_1-1}\Mod{e}$}\label{sec:labels2}

We note that $\kappa_2\equiv{\kappa_1-1}\Mod{e}$ throughout this subsection.

For $n\not\equiv{0}\Mod{e}$ and $1\leqslant{m}\leqslant{n-1}$, we recall from~\cite[Proposition 6.13]{Sutton17I} that $S_{((n-m),(1^m))}$ has two composition factors, namely $\im(\chi_m)$ and $S_{((n-m),(1^m))}/\im(\chi_m)$. Thus $\im(\chi_m)\cong{D_{\lambda_m}}$ and $S_{((n-m),(1^m))}/\im(\chi_m)\cong{D_{\mu_m}}$ for some regular bipartitions $\lambda_m,\mu_m\in\mathscr{RP}_n^2$.

\begin{defn}\label{def:labels,case 3}
Let $\kappa_2\equiv{\kappa_1-1}\Mod{e}$ and $n\not\equiv{0}\Mod{e}$.
For $1\leqslant{m}\leqslant{n-1}$, we define
\begin{align*}
\mu_{n,2m}&:=
\begin{cases}
((n-m,\{m\}),\varnothing)
&\hspace{50pt}\text{ if $1\leqslant{m}<n-\frac{n}{e}$,}\\
((\{m+1\},n-1-m),\varnothing)
&\hspace{50pt}\text{ if $n-\frac{n}{e}\leqslant{m}\leqslant{n-1}$,}
\end{cases}\\
\\
\mu_{n,2m+1}&:=
\begin{cases}
((n-m),(1^m))
&\text{ if $1\leqslant{m}<e$,}\\
((n-m,\{m-e\}),(2,1^{e-2}))
&\text{ if $e\leqslant{m}<n-\frac{n}{e}$,}\\
((\{m-e+1\},n-1-m),(2,1^{e-2}))
&\text{ if $n-\frac{n}{e}\leqslant{m}\leqslant{n-1}$.}
\end{cases}
\end{align*}
\end{defn}

Notice that $\mu_{n,2m}$ and $\mu_{n,2m+1}$ are distinct. We claim that the labels $\lambda_m,\mu_m$ of the two composition factors of $S_{((n-m),(1^m))}$ are $\mu_{n,2m}$ and $\mu_{n,2m+1}$, respectively, and hence that the corresponding composition factors are non-isomorphic.

For $n\equiv{0}\Mod{e}$, we recall from~\cite[Theorem 6.16]{Sutton17I} that $S_{((n-m),(1^m))}$ has four composition factors $\im(\phi_m)$, $\im(\phi_{m+1})$, $\ker(\gamma_m)/\im(\phi_m)$ and $\ker(\gamma_{m+1})/\im(\phi_{m+1})$ for $m\in\{2,\dots,n-2\}$, and that $S_{((n-1),(1))}$ and $S_{((1),(1^{n-1}))}$ both have three composition factors. Thus, for $m\in\{2,\dots,n-2\}$, $\im(\phi_m)\cong{D_{\lambda_m}}$, $\im(\phi_{m+1})\cong{D_{\mu_m}}$, $\ker(\gamma_m)/\im(\phi_m)\cong{D_{\nu_m}}$ and $\ker(\gamma_{m+1})/\im(\phi_{m+1})\cong{D_{\eta_m}}$ for some regular bipartitions $\lambda_m,\mu_m,\nu_m,\eta_m\in\mathscr{RP}_n^2$.

\begin{defn}\label{def:labels,case 4}
Let $\kappa_2\equiv{\kappa_1-1}\Mod{e}$ and $n\equiv{0}\Mod{e}$.
For $2\leqslant{m}\leqslant{n-1}$, we define
\begin{align*}
\mu_{n,2m}&:=
	\begin{cases}
	((n-m+1,\{m-1\}),\varnothing)
	&\hspace{50pt}\text{ if $2\leqslant{m}\leqslant{n-\tfrac{n}{e}}$,}\\
	((\{m+1\},n-m-1),\varnothing)
	&\hspace{50pt}\text{ if $n-\tfrac{n}{e}<m\leqslant{n-1}$,}
	\end{cases}\\
\\
\mu_{n,2m+1}&:=
	\begin{cases}
	((n-m+1),(1^{m-1}))
	&\text{ if $2\leqslant{m}\leqslant{e}$,}\\
	((n-m+1,\{m-e-1\}),(2,1^{e-2}))
	&\text{ if $e<{m}\leqslant{n-\tfrac{n}{e}}$,}\\
	((\{m-e+1\},n-m-1),(2,1^{e-2}))
	&\text{ if $n-\tfrac{n}{e}<m\leqslant{n-1}$.}
	\end{cases}
\end{align*}
\end{defn}

We notice that $\mu_{n,2m}$, $\mu_{n,2m+1}$, $\mu_{n,2m+2}$ and $\mu_{n,2m+2}$ are distinct bipartitions. For $m\in\{2,\dots,n-2\}$, we claim that the labels $\lambda_m,\mu_m,\nu_m,\eta_m$ of the four composition factors of $S_{((n-m),(1^m))}$ are $\mu_{n,2m}$, $\mu_{n,2m+2}$, $\mu_{n,2m+1}$ and $\mu_{n,2m+3}$, respectively, and hence that the corresponding composition factors are non-isomorphic.

To confirm our claims above, we need the following combinatorial result, which is analogous to \cref{lem:labels1} and can be proved in a similar manner.

\begin{lem}\label{lem:labels2}
Suppose that $\kappa_2\equiv{\kappa_1-1}\Mod{e}$.
\begin{enumerate}
\item{
Let $m\in\{1,\dots,n-1\}$. If $n\not\equiv{0}\Mod{e}$, then
\begin{subequations}
\begin{align}
\mu_{n,2m}\uparrow_{\kappa_2-m}=\mu_{n+1,2m+2},
\label{leg3}
\\
\mu_{n,2m+1}\uparrow_{\kappa_2-m}=\mu_{n+1,2m+3}.
\label{leg4}
\end{align}
\end{subequations}
}
\item {
Let $m\in\{2,\dots,n-1\}$. If $n\equiv{0}\Mod{e}$, then
\begin{subequations}
\begin{align}
\mu_{n,2m}\uparrow_{\kappa_2+1-m}&=\mu_{n+1,2m},
\label{leg5} 
\\
\mu_{n,2m+1}\uparrow_{\kappa_2+1-m}&=\mu_{n+1,2m+1}.
\label{leg6}
\end{align}
\end{subequations}
}
\end{enumerate}
\end{lem}

\begin{thm}\label{labels:case 3}
Suppose that $\kappa_2\equiv{\kappa_1-1}\Mod{e}$ and $n\not\equiv{0}\Mod{e}$. Then the composition factors of $S_{((n-m),(1^m))}$ are $D_{\mu_{n,2m}}$ and $D_{\mu_{n,2m+1}}$ for all $m\in\{1,\dots,n-1\}$.
Moreover, ${D_{\mu_{n,2m}}}\cong{\im(\chi_m)}$ and ${D_{\mu_{n,2m+1}}}\cong{S_{((n-m),(1^m))}/\im(\chi_m)}$ as ungraded $\mathscr{H}_n^{\Lambda}$-modules.
\end{thm}

\begin{proof}
We first show that $D_{\mu_{n,3}}$ is a composition factor of $S_{((n-1),(1))}$. We have $f_{\kappa_2-1}^{(2)}S_{((n-1),(1))}\cong S_{((n),(1^2))} \text{ if } n\equiv{-1}\Mod{e}$, and $f_{\kappa_2-1}S_{((n-1),(1))}\cong S_{((n-1),(1^2))} \text{ if } n\not\equiv{-1}\Mod{e}$.

For $n\equiv{-1}\Mod{e}$, $D_{\mu_{n+2,5}}$ is a composition factor of $S_{((n),(1^2))}$ by downwards induction on $n$. Hence, by \cref{irrcomp}, $D_{\mu_{n+2,5}}\downarrow_{\kappa_2-1}^2$ is a composition factor of $S_{((n-1),(1))}$. We have
\begin{align*}
((n-1),(1))\uparrow_{\kappa_2-1}^2
=
\mu_{n,3}\uparrow_{\kappa_2-1}^2
&=
\mu_{n+1,5}\uparrow_{\kappa_2-1}
&&\text{(\Cref{leg4})}\\
&=\mu_{n+2,5}
&&\text{(\Cref{leg6})}\\
&=((n),(1^2)).
\end{align*}
Its inverse gives us $\mu_{n,3}=\mu_{n+2,5}\downarrow_{\kappa_2-1}^2$, and hence $D_{\mu_{n,3}}$ is a composition factor of $S_{((n-1),(1))}$.

Similarly, for $n\not\equiv{-1}\Mod{e}$, $D_{\mu_{n+1,5}}$ is a composition factor of $S_{((n-1),(1^2))}$. Thus, by \cref{irrcomp}, $D_{\mu_{n+1,5}}\downarrow_{\kappa_2-1}$ is a composition factor of $S_{((n-1),(1))}$. Observe that
\[
((n-1),(1))\uparrow_{\kappa_2-1}
=
\mu_{n,3}\uparrow_{\kappa_2-1}
=
\mu_{n+1,5}
=((n-1),(1^2))
\qquad\text{(\Cref{leg4})}.
\]
Its inverse gives us $\mu_{n,3}=\mu_{n+1,5}\downarrow_{\kappa_2-1}$, and hence $D_{\mu_{n,3}}$ is a composition factor of $S_{((n-1),(1))}$.

\begin{enumerate}
\item{Suppose that $n-l\not\equiv{2}\Mod{e}$. We have $e_{\kappa_2+1-m}S_{((n-m),(1^m))}\cong{S_{((n-m),(1^{m-1}))}}$, and by induction, $D_{\mu_{n-1,2m-2}}$ and $D_{\mu_{n-1,2m-1}}$ are composition factors of $S_{((n-m),(1^{m-1}))}$. It thus follows from \cref{irrcomp} that $D_{\mu_{n-1,2m-2}\uparrow_{\kappa_2+1-m}}$ and $D_{\mu_{n-1,2m-1}\uparrow_{\kappa_2+1-m}}$ are composition factors of $S_{((n-m),(1^m))}$.
We observe that $\mu_{n-1,2m-2}\uparrow_{\kappa_2+1-m}=\mu_{n,2m}$ by \Cref{leg3},
and $\mu_{n-1,2m-1}\uparrow_{\kappa_2+1-m}=\mu_{n,2m+1}$ by \Cref{leg4}.
Hence $D_{\mu_{n,2m}}$ and $D_{\mu_{n,2m+1}}$ are composition factors of $S_{((n-m),(1^m))}$.
}
\item{
Suppose that $n-l\equiv{2}\Mod{e}$. We have $e_{\kappa_2+1-m}^{(2)}S_{((n-m),(1^m))}\cong{S_{((n-m-1),(1^{m-1}))}}$, and by induction, $D_{\mu_{n-2,2m-2}}$ and $D_{\mu_{n-2,2m-1}}$ are composition factors of $S_{((n-m-1),(1^{m-1}))}$.
We observe that
\begin{align*}
\mu_{n-2,2m-2}\uparrow_{\kappa_2+1-m}^2
&=\mu_{n-1,2m-2}\uparrow_{\kappa_2+1-m}
&&\text{(\Cref{leg5})}\\
&=\mu_{n,2m}
&&\text{(\Cref{leg3})}.
\end{align*}
Thus, by \cref{irrcomp}, $D_{\mu_{n,2m}}$ is a composition factor of $S_{((n-m),(1^m))}$.

We also observe that
\begin{align*}
\mu_{n-2,2m-1}\uparrow_{\kappa_2+1-m}^2
&=\mu_{n-1,2m-1}\uparrow_{\kappa_2+1-m}
&&\text{(\Cref{leg6})}\\
&=\mu_{n,2m+1}
&&\text{(\Cref{leg4})}.
\end{align*}
Thus, by \cref{irrcomp}, $D_{\mu_{n,2m+1}}$ is composition factor of $S_{((n-m),(1^m))}$.
}
\end{enumerate}

Furthermore, we know from~\cite[Proposition 6.13]{Sutton17I} that the composition factors $D_{\mu_{n,2m}}$ and $D_{\mu_{n,2m+1}}$ of $S_{((n-m),(1^m))}$ are in bijection with $\operatorname{im}(\chi_m)$ and $S_{((n-m),(1^m))}/\operatorname{im}(\chi_m)$, up to isomorphism and grading shift.
By~\cite[Lemma 5.10]{Sutton17I},
\begin{itemize}
	\item $\operatorname{im}(\chi_m)=\operatorname{span}\left\{v_\ttt\ \mid \ \ttt\in\operatorname{Std}((n-m),(1^m)),\ttt(1,1,1)=1\right\}$;
	\item $S_{((n-m),(1^m))}/\operatorname{im}(\chi_m)=\operatorname{span}\left\{v_\ttt\ \mid \ \ttt\in\operatorname{Std}((n-m),(1^m)),\ttt(1,1,2)=1\right\}$.
\end{itemize}
Now let $\tts,\ttt\in\operatorname{Std}((n-m),(1^m))$ be such that $1$ lies in the arm of $\ttt$ and $1$ lies in the leg of $\tts$.
Then every tableau $\ttt$ has residue sequence $(\kappa_1,i_2,\dots,i_n)$ where $i_r\in\{0,\dots,e-1\}$, and every tableau $\tts$ has residue sequence $(\kappa_2,j_2,\dots,j_n)$ where $j_r\in\{0,\dots,e-1\}$.
The only non-empty component of $\mu_{n,2m}$ is its first component, whereas both of the components of $\mu_{n,2m+1}$ are non-empty. Thus, only the residue sequence of $\mu_{n,2m+1}$ can begin with residue $\kappa_2$, and hence ${D_{\mu_{n,2m}}}\cong{\operatorname{im}(\chi_m)}$, as required.
\end{proof}

Similarly to the results in \cref{sec:labels1}, we use this result to describe the composition factors of Specht modules labelled by hook bipartitions in the following case.

\begin{thm}\label{labels:case 4}
Suppose that $\kappa_2\equiv{\kappa_1-1}\Mod{e}$, $n\equiv{0}\Mod{e}$ and let $m\in\{1,\dots,n-1\}$. Then $S_{((n-m),(1^m))}$ has composition factors

\begin{enumerate}
	\item $S_{((n),\varnothing)}$, $D_{\mu_{n,4}}$ and $D_{\mu_{n,5}}$ if $m=1$;
	\item $D_{\mu_{n,2m}}$, $D_{\mu_{n,2m+1}}$, $D_{\mu_{n,2m+2}}$ and $D_{\mu_{n,2m+3}}$ if $m\in\{2,\dots,n-2\}$;
	\item $S_{(\varnothing,(1^n))}$, $D_{\mu_{n,2n-2}}$ and $D_{\mu_{n,2n-1}}$ if $m=n-1$.
\end{enumerate}

Moreover, ${D_{\mu_{n,2m}}}\cong{\im(\phi_m)}$ and ${D_{\mu_{n,2m+1}}}\cong{\ker(\gamma_m)/\im(\phi_m)}$ as ungraded $\mathscr{H}_n^{\Lambda}$-modules.
\end{thm}

\begin{proof}
\begin{enumerate}[label=(\roman*)]
\item{Firstly, by removing the foot node of ${((n-1),(1))}$, we have \[e_{\kappa_1}S_{((n-1),(1))}\cong{S_{((n-1),\varnothing)}}\cong{D_{((n-1),\varnothing)}}.\] The $\kappa_2$-signature of $((n-1),\varnothing)$ is $++$, corresponding to the conormal nodes $(1,n,1)$ and $(1,1,2)$. Adding the higher of these nodes, $((n-1),\varnothing)\uparrow_{\kappa_2}=((n),\varnothing)$, and by \cref{irrcomp}, $D_{((n),\varnothing)}$ is a composition factor of $S_{((n-1),(1))}$.

Now suppose that $2\leqslant{m}\leqslant{n-1}$.
By removing the foot node of $((n-m),(1^m))$, we have
\[
e_{\kappa_2+1-m}S_{((n-m),(1^m))}\cong{S_{((n-m),(1^{m-1}))}}.
\]
It follows from \cref{labels:case 3} that $D_{\mu_{n-1,2m-2}}$ and $D_{\mu_{n-1,2m-1}}$ are composition factors of $S_{((n-m),(1^{m-1}))}$. Observe that $\mu_{n-1,2m-2}\uparrow_{\kappa_2+1-m}=\mu_{n,2m}$ by \Cref{leg3}, and that $\mu_{n-1,2m-1}\uparrow_{\kappa_2+1-m}=\mu_{n,2m+1}$ by \Cref{leg4}. Thus, by \cref{irrcomp}, both $D_{\mu_{n,2m}}$ and $D_{\mu_{n,2m+1}}$ are composition factors of $S_{((n-m),(1^m))}$.
}
\item{First suppose that $1\leqslant{m}\leqslant{n-2}$.
By removing the hand node of $((n-m),(1^m))$, we have
\[
e_{\kappa_2-m}S_{((n-m),(1^m))}\cong{S_{((n-m-1),(1^m))}}.
\]
By \cref{labels:case 3}, $D_{\mu_{n-1,2m}}$ and $D_{\mu_{n-1,2m+1}}$ are composition factors of  $S_{((n-m-1),(1^m))}$. Observe that $\mu_{n-1,2m}\uparrow_{\kappa_2-m}$ by \Cref{leg3}, and that $\mu_{n-1,2m+1}\uparrow_{\kappa_2-m}=\mu_{n,2m+3}$ by \Cref{leg4}. Thus, $D_{\mu_{n,2m+2}}$ and $D_{\mu_{n,2m+3}}$ are also composition factors of $S_{((n-m),(1^m))}$ by \cref{irrcomp}.

Secondly, suppose that $m=n-1$. By removing the hand node of ${((1),(1^{n-1}))}$, we have
\[
{e_{\kappa_1}S_{((1),(1^{n-1}))}}\cong
{S_{(\varnothing,(1^{n-1}))}}\cong
{D_{(\{n-e\},(1^{e-1}))}}
\qquad\text{(by \cref{simpsign})}.
\]
The $\kappa_1$-signature of $(\{n-e\},(1^{e-1}))$ is $+++$, corresponding to the conormal nodes $(1,\lfloor(n-2)/(e-1)\rfloor+1,1)$, $(1,2,2)$ and $(e,1,2)$. Adding the highest of these nodes, we have $(\{n-e\},(1^{e-1}))\uparrow_{\kappa_1}=(\{n-e+1\},(1^{e-1}))$. By \cref{simpsign}, ${D_{(\{n-e+1\},(1^{e-1}))}}\cong{S_{(\varnothing,(1^n))}}$, and hence $S_{(\varnothing,(1^n))}$ is a composition factor of $S_{((1),(1^{n-1}))}$ by \cref{irrcomp}.
}
\end{enumerate}

Furthermore, for all $m\in\{2,\dots,n-2\}$, we know from~\cite[Theorem 6.16]{Sutton17I} that the composition factors $D_{\mu_{n,2m}}$, $D_{\mu_{n,2m+1}}$, $D_{\mu_{n,2m+2}}$ and $D_{\mu_{n,2m+3}}$ of $S_{((n-m),(1^m))}$ are in bijection with $\operatorname{im}(\phi_m)$, $\operatorname{im}(\phi_{m+1})$, $\operatorname{ker}(\gamma_m)/\operatorname{im}(\phi_m)$ and $\operatorname{ker}(\gamma_{m+1})/\operatorname{im}(\phi_{m+1})$, up to isomorphism and grading shift. Moreover, $\operatorname{im}(\phi_{m+1})$ and $\operatorname{ker}(\gamma_{m+1})/$ $\operatorname{im}(\phi_{m+1})$ are composition factors of both $S_{((n-m),(1^m))}$ and $S_{((n-m-1),(1^{m+1})}$, and hence are in bijection with $D_{\mu_{n,2m+2}}$ and $D_{\mu_{n,2m+3}}$, up to isomorphism and grading shift.

Let $\mathscr{T}=\operatorname{Std}((n-m),(1^m))$. Then, by~\cite[Lemmas 5.9 and 5.10]{Sutton17I}, we have that
\begin{itemize}
	\item $\operatorname{im}(\phi_{m+1})
	\cong\operatorname{span}\left\{v_{\ttt}\ |\ \ttt\in\mathscr{T},\ttt(1,1,1)=1,\ttt(1,n-m,1)=n\right\}$;
	\item $\operatorname{ker}(\gamma_{m+1})/\operatorname{im}(\phi_{m+1})\cong
	\operatorname{span}\left\{v_{\ttt}\ |\ \ttt\in\mathscr{T},\ttt(1,1,2)=1,\ttt(1,n-m,1)=n\right\}$.
\end{itemize}
It follows, together with~\cite[Lemma 5.10]{Sutton17I}, that
\begin{itemize}
\item $\ttt(1,1,1)=1$ if $v_{\ttt}$ lies in either $\operatorname{im}(\phi_m)$ or $\operatorname{im}(\phi_{m+1})$;
\item $\ttt(1,1,2)=1$ if $v_{\ttt}$ lies in either $\operatorname{ker}(\gamma_m)/\operatorname{im}(\phi_m)$ or $\operatorname{ker}(\gamma_{m+1})/\operatorname{im}(\phi_{m+1})$.
\end{itemize}
We now observe that only the first component of $\mu_{n,2m}$ is non-empty, whereas both components of $\mu_{n,2m+1}$ are non-empty. It follows that $1$ can only lie in the leg of $\ttt$ if $v_{\ttt}$ lies in $D_{\mu_{n,2m+1}}$ or $D_{\mu_{n,2m+3}}$, and hence ${D_{\mu_{n,2m}}}\cong{\operatorname{im}(\phi_m)}$ and ${D_{\mu_{n,2m+1}}}\cong{\operatorname{ker}(\gamma_m)/\operatorname{im}(\phi_m)}$, as required.
\end{proof}

\section{Ungraded decomposition numbers corresponding to $S_{((n-m),(1^m))}$}\label{sec:UDN}

We remind the reader that we found the characteristic-free composition series of Specht modules labelled by hook bipartitions in terms of the basis vectors of $S_{((n-m),(1^m))}$ in~\cite[\S{6}]{Sutton17I}, and furthermore, in \Cref{sec:labels1,sec:labels2} we established the regular bipartitions that label these composition factors. We can thus determine the \emph{ungraded} multiplicities $[{S_{((n-m),(1^m))}}:{D_{\mu}}]$ for all regular bipartitions $\mu\in\mathscr{RP}_n^2$. Recall that $l\equiv\kappa_2-\kappa_1\Mod{e}$.

\subsection[Case I: $\kappa_2\not\equiv{\kappa_1-1}\Mod{e}$ and $n\not\equiv{l+1}\Mod{e}$]{Case I: $\kappa_2\not\equiv{\kappa_1-1}\Mod{e}$ and $n\not\equiv{l+1}\Mod{e}$}

We recall from~\cite[Theorem 6.8]{Sutton17I} that $S_{((n-m),(1^m))}$ is irreducible for all $m\in\{0,\dots,n\}$. Moreover, we know from \cref{simpsign} that $S_{(\varnothing,(1^n))}\cong D_{(\varnothing,(1^n))^R}$ and from \cref{factors} that $S_{((n-m),(1^m))}\cong{D_{\mu_{n,m}}}$ for all $m<n$, leading us to the following result.

\begin{thm}
Let $\kappa_2\not\equiv{\kappa_1-1}\Mod{e}$ and $n\not\equiv{l+1}\Mod{e}$. Then the decomposition submatrix $\left(d_{((n-m),(1^m))},\mu\right)$ of $\mathscr{H}_n^{\Lambda}$, under a specific ordering on its columns, is
\[
\begin{blockarray}{ccccccccc}
\begin{block}{c(ccccc|ccc)}
S_{((n),\varnothing)}
	& 1 &&&&
	&&&\\
S_{((n-1),(1))}
	& & 1 && \text{\huge$0$} &
	&&&\\
S_{((n-2),(1^2))}
	&&& 1 &&
	&& \text{\Huge$0$} &\\
\vdots
	&& \text{\huge$0$} && \ddots & 
	&&& \\
S_{(\varnothing,(1^n))}
	&&&&& 1 
	&&&\\
\end{block}
\end{blockarray}
\]
for all regular bipartitions $\mu\in\mathscr{RP}_n^2$.
\end{thm}

\subsection[Case II: $\kappa_2\not\equiv{\kappa_1-1}\Mod{e}$ and $n\equiv{l+1}\Mod{e}$]{Case II: $\kappa_2\not\equiv{\kappa_1-1}\Mod{e}$ and $n\equiv{l+1}\Mod{e}$}

We know from \cref{thm:labels case 2} that the composition factors of $S_{((n-m),(1^m))}$ are $D_{\mu_{n,m-1}}$ and $D_{\mu_{n,m}}$ for all $m\in\{1,\dots,n-1\}$. Hence $D_{\mu_{n,m}}$ is a composition factor of both $S_{((n-m),(1^m))}$ and $S_{((n-m-1),(1^{m+1}))}$ whenever $m\in\{1,\dots,n-2\}$. We also note that $D_{\mu_{n,0}}=S_{((n),\varnothing)}$ and $D_{\mu_{n,n-1}}=D_{(\varnothing,(1^n))^R}$. Furthermore, since the bipartitions $\mu_{n,0},\mu_{n,1},\dots,\mu_{n,n-1}$ are distinct, the irreducible modules $D_{\mu_{n,0}},D_{\mu_{n,1}},\dots,{D_{\mu_{n,n-1}}}$ are non-isomorphic.

\begin{thm}
Let $\kappa_2\not\equiv{\kappa_1-1}\Mod{e}$ and $n\equiv{\kappa_2-\kappa_1+1}\Mod{e}$. 
Then the decomposition submatrix $\left(d_{((n-m),(1^m))},\mu\right)$ of $\mathscr{H}_n^{\Lambda}$, under a specific ordering on its columns, is
\[
\begin{blockarray}{cccccccccc}
\begin{block}{c(cccccc|ccc)}
 S_{((n),\varnothing)}& 1 &  &   &   &    & 									&&\\
 S_{((n-1),(1))}&  1 & 1 &  & &   \text{\huge$0$}  &                         && \\
 S_{((n-2),(1^2))} &  & 1  & 1 &  &   &                                      &&     \\
 S_{((n-3),(1^3))} & & & 1 & 1 & &											&&\text{\Huge$0$} \\
 \vdots& &     & & \ddots & \ddots &                                        &&   \\
 S_{((1),(1^{n-1}))}  & & \text{\huge$0$}  & & &  1  & 1 &                  &&  \\
 S_{(\varnothing,(1^n))}   &   & & &  & & 1                                   & &\\
\end{block}
\end{blockarray}
 \]
 for all regular bipartitions $\mu\in\mathscr{RP}_n^2$.
\end{thm}

\subsection[Case III: $\kappa_2\equiv{\kappa_1-1}\Mod{e}$ and $n\not\equiv{0}\Mod{e}$]{Case III: $\kappa_2\equiv{\kappa_1-1}\Mod{e}$ and $n\not\equiv{0}\Mod{e}$}

We know from \cref{labels:case 3} that the composition factors of $S_{((n-m),(1^m))}$ are $D_{\mu_{n,2m}}$ and $D_{\mu_{n,2m+1}}$ for all $m\in\{1,\dots,n-1\}$. Furthermore, since the bipartitions $((n),\varnothing)$, $\mu_{n,2}$, $\mu_{n,3},\dots,\mu_{2n-1}$, $(\varnothing,(1^n))^R$ are distinct, we know that the irreducible modules $S_{((n),\varnothing)}$, $D_{\mu_{n,2}}$, $D_{\mu_{n,3}},\dots,{D_{\mu_{n,2n-1}}}$, $D_{(\varnothing,(1^n))^R}$ are non-isomorphic.

\begin{thm}
Let $\kappa_2\equiv{\kappa_1-1}\Mod{e}$ and $n\not\equiv{0}\Mod{e}$.
Then the decomposition submatrix $\left(d_{((n-m),(1^m))},\mu\right)$ of $\mathscr{H}_n^{\Lambda}$, under a specific ordering on its columns, is
\[
\begin{blockarray}{cccccccccccccc}
\begin{block}{c(cccccccccc|ccc)}
S_{((n),\varnothing)} & 1 &&&&&&						&&&	 &  &    \\
S_{((n-1),(1))} && 1 & 1 &&&&&& \text{\huge$0$} 			&	 &  &    \\
S_{((n-2),(1^2))} &&&& 1 & 1 &&&						&&&	&&& \\
S_{((n-3),(1^3))} &&&&&& 1 & 1 &						&&&&\text{\Huge$0$}& \\
\vdots &&&&&&& \ddots & \ddots  								&&	 &  &  \\
S_{((1),(1^{n-1}))} && \text{\huge$0$} &&&&&& 1 & 1 & 		&&& \\
S_{(\varnothing,(1^n))} &&&&&&&&&& 1						 &&    \\ 
\end{block}
\end{blockarray}
\]
for all regular bipartitions $\mu\in\mathscr{RP}_n^2$.
\end{thm}

\subsection[Case IV: $\kappa_2\equiv{\kappa_1-1}\Mod{e}$ and $n\equiv{0}\Mod{e}$]{Case IV: $\kappa_2\equiv{\kappa_1-1}\Mod{e}$ and $n\equiv{0}\Mod{e}$}

We recall from \cref{labels:case 4} that the composition factors of $S_{((n-m),(1^m))}$ are: $S_{((n),\varnothing)}$, $D_{\mu_{n,4}}$ and $D_{\mu_{n,5}}$ if $m=1$; $D_{\mu_{n,2m}}$, $D_{\mu_{n,2m+1}}$, $D_{\mu_{n,2m+2}}$ and $D_{\mu_{n,2m+3}}$ if $m\in\{2,\dots,n-2\}$; $D_{\mu_{n,2n-2}}$, $D_{\mu_{n,2n-1}}$ and $D_{(\varnothing,(1^n))^R}$ if $m=n-1$.
\begin{comment}
\begin{itemize}
\item{ $D_{((n),\varnothing)}$, $D_{\mu_{n,4}}$ and $D_{\mu_{n,5}}$ if $m=1$; }
\item{ $D_{\mu_{n,2m}}$, $D_{\mu_{n,2m+1}}$, $D_{\mu_{n,2m+2}}$ and $D_{\mu_{n,2m+3}}$ if $m\in\{2,\dots,n-2\}$; }
\item{ $D_{\mu_{n,2n-2}}$, $D_{\mu_{n,2n-1}}$ and $D_{(\varnothing,(1^n))^R}$ if $m=n-1$. }
\end{itemize}
\end{comment}
Thus, for all $m\in\{1,\dots,n-2\}$, $D_{\mu_{n,2m+2}}$ and $D_{\mu_{n,2m+3}}$ are composition factors of both $S_{((n-m),(1^m))}$ and $S_{((n-m-1),(1^{m+1}))}$. Furthermore, since the bipartitions $((n),\varnothing)$, $\mu_{n,4},\dots,\mu_{n,2n-1}$, $(\varnothing,(1^n))^R$ are distinct, the irreducible modules $S_{((n),\varnothing)}$, $D_{\mu_{n,4}},\dots,D_{\mu_{n,2n-1}}$, $D_{(\varnothing,(1^n))^R}$ are non-isomorphic.

\begin{thm}
Let $\kappa_2\equiv{\kappa_1-1}\Mod{e}$ and $n\equiv{0}\Mod{e}$.
Then the decomposition submatrix $\left(d_{((n-m),(1^m))},\mu\right)$ of $\mathscr{H}_n^{\Lambda}$, under a specific ordering on its columns, is
\[
\begin{blockarray}{ccccccccccccccccc}
\begin{block}{c(cccccccccccc|cccc)}
S_{((n),\varnothing)} &  1&&&&&&&&&							&&	 &  & &  \\
S_{((n-1),(1))} &  1&1&1&&&&& &&\text{\huge$0$} 			&&	 &  &  &  \\
S_{((n-2),(1^2))} &  &1&1&1&1&&&&&							&&	&&& \\
S_{((n-3),(1^3))} &  &&&1&1&1&1&&&							&&	&&& \\
S_{((n-4),(1^4))} &&&&&&1&1&1&1&							&&&&\text{\Huge$0$}&\\
\vdots &  &&&&&& \ddots & \ddots & \ddots & \ddots 				&&	 &  && 	\\
S_{((2),(1^{n-2}))} &  &&&&&&&1&1&1&1& 						&	&&	\\
S_{((1),(1^{n-1}))} &  && \text{\huge$0$} &&&&&&&1&1&1		&	&&	\\
S_{(\varnothing,(1^n))} &  &&&&&&&&&&&1						 &  & & 	\\
\end{block}
\end{blockarray}
\]
for all regular bipartitions $\mu\in\mathscr{RP}_n^2$.
\end{thm}

\begin{rmk}
Notice that the above decomposition submatrices of $\mathscr{H}_n^{\Lambda}$ are independent of the characteristic of the ground field, and thus the corresponding adjustment submatrices are trivial.
\end{rmk}

\section{Graded dimensions of $S_{((n-m),(1^m))}$}\label{sec:GDS}

From now on, we study \emph{graded} Specht modules labelled by hook bipartitions, using the combinatorial $\mathbb{Z}$-grading defined on these $\mathscr{H}_n^{\Lambda}$-modules to determine their graded dimensions.

We first determine the removable and addable $i$-nodes of hook bipartitions as follows.

\begin{lem}\label{lem:entries}
Let $1\leqslant{i}\leqslant{k}$. Then $((k-i),(1^i))$ has neither an addable nor a removable $(\kappa_2+1-i)$-node in the first row of the first component, except in the following cases.
\begin{enumerate}[label=(\roman*)]
\item{
If $k\equiv{l+1}\Mod{e}$, then $(1,k-i+1,1)$ is an addable $(\kappa_2+1-i)$-node of $((k-i),(1^i))$.
}
\item{If $k\equiv{l+2}\Mod{e}$ and $k>i$, then $(1,k-i,1)$ is a removable $(\kappa_2+1-i)$-node of $((k-i),(1^i))$.}
\end{enumerate}
\end{lem}

\begin{proof}Let $\ttt\in\std((n-m),(1^m))$ be such that $\ttt(i,1,2)=k$.
\begin{enumerate}
\item
Suppose that $\ttt(i,1,2)=l+1+\alpha{e}$ for some $\alpha\in\mathbb{N}\cup\{0\}$. Then $1,\dots,l+\alpha{e}$ must lie in the set of nodes 
$\{(1,1,2),\dots,(i-1,1,2)\}
\cup
\{(1,1,1),\dots,(1,j,1)\}$, where $j=l+\alpha{e}-i+1$.
There are $j$ and $i-1$ entries strictly smaller than $l+\alpha{e}+1$ in the arm and the leg of $\ttt$, respectively. We now observe that	$\res(1,j+1,1) = \kappa_1+j = \kappa_2-i+1 = \res(i,1,2) \Mod{e}$, and since $\ttt(i,1,2)>\ttt(1,j,1)$, it follows that $(1,j+1,1)=(1,k-i+1,1)$ is an addable $(\kappa_2+1-i)$-node for $((k-i),(1^i))$.
\item
Suppose that $\ttt(i,1,2)=l+k+\alpha{e}$ for some $\alpha\in\mathbb{N}\cup\{0\}$ such that $k\in\{2,\dots,e\}$, and prove this is in a similar fashion to the first part, treating the cases $k=2$ and $k>2$ separately.
\qedhere
\end{enumerate}
\end{proof}

For any $\ttt\in\std((n-m),(1^m))$, we define
\[
a_{\ttt}:=\#\{
i\mid
\ttt(i,1,2)\equiv{l+1}\Mod{e}
\}-\#\{i\mid
\ttt(i,1,2)\equiv{l+2}\Mod{e}\}.
\]
We are now able to obtain the degree of an arbitrary standard $((n-m),(1^m))$-tableau.

\begin{lem}\label{lem:deg}
Let $\ttt\in\std((n-m),(1^m))$ and $1\leqslant{i}\leqslant{m}<n$. Then
\begin{align*}
\deg(\ttt)=\left\lfloor\tfrac{m+e-l-2}{e}\right\rfloor
+\lfloor\tfrac{l+1}{e}\rfloor
+\left\lfloor\tfrac{m}{e}\right\rfloor+a_{\ttt}.
\end{align*}
\end{lem}

\begin{proof}
Suppose that $\ttt(i,1,2)=k$ for some $k\in\{i,\dots,n\}$, so that $\ttt_{\leqslant k}$ is a standard $((k-i),(1^i))$-tableau. Applying \cref{lem:entries}, we have
\begin{align*}
\deg(\ttt) = &\
\#\{i\mid
(i,1,2)\text{ has addable $(\kappa_1-1)$-node }(2,1,1)\}\\
&+\#\{i\mid(i,1,2)\text{ has addable $(\kappa_2+1)$-node }(1,2,2)\}\\
&+\#\{i\mid
(i,1,2)\text{ has addable $(\kappa_2+1-i)$-node in the first row of $\ttt$}\}\\
&-\#\{i\mid
(i,1,2)\text{ has removable $(\kappa_2+1-i)$-node in the first row of $\ttt$}\}\\
= &\ \#\{i\mid
i\equiv{l+2}\Mod{e},k>i\}\\
&+\#\{i\mid
i\equiv{0}\Mod{e}\}\\
&+\#\{i\mid
k\equiv{l+1}\Mod{e}\}\\
&-\#\{i\mid
k\equiv{l+2}\Mod{e},k>i\}\\
= &\ \#\{i\mid
i\equiv{l+2}\Mod{e}\}
-\#\{i\mid
i\equiv{l+2}\Mod{e},k=i\}\\
&+\#\{i\mid
i\equiv{0}\Mod{e}\}\\
&+\#\{i\mid
k\equiv{l+1}\Mod{e}\}\\
&-\#\{i\mid
k\equiv{l+2}\Mod{e}\}
+\#\{i\mid
k\equiv{l+2}\Mod{e},k=i\}\\
= &\ \#\{i\mid
i\equiv{l+2}\Mod{e}\}
+\#\{i\mid
i\equiv{0}\Mod{e}\}\\
& +\#\{i\mid
k\equiv{l+1}\Mod{e}\}
-\#\{i\mid
k\equiv{l+2}\Mod{e}\}\\
= &\ \left\lfloor\tfrac{m+e-l-2}{e}\right\rfloor
+\lfloor\tfrac{l+1}{e}\rfloor
+\left\lfloor\tfrac{m}{e}\right\rfloor
+\#\{i\mid
k\equiv{l+1}\Mod{e}\}
-\#\{i\mid
k\equiv{l+2}\Mod{e}\}.
\end{align*}
\qedhere
\end{proof}

For any non-empty subset $\mathscr{T}\subseteq\std((n-m),(1^m))$, we define the set $A_{\mathscr{T}}:=\left\{a_{\ttt}\ | \ \ttt\in\mathscr{T}\right\}$. We now define the \emph{maximum degree of $\mathscr{T}$} to be $\maxdeg(\mathscr{T}):=\mx\{\deg(\ttt)\:|\:\ttt\in\mathscr{T}\}$ and the \emph{minimum degree of $\mathscr{T}$} to be $\mindeg(\mathscr{T}):=\min\{\deg(\ttt)\:|\:\ttt\in\mathscr{T}\}$. By \cref{lem:deg}, we have
\begin{itemize}
	\item $\maxdeg(\mathscr{T})=
	\left\lfloor\tfrac{m+e-l-2}{e}\right\rfloor
	+\left\lfloor\tfrac{l+1}{e}\right\rfloor
	+\left\lfloor\tfrac{m}{e}\right\rfloor
	+\mx(A_{\mathscr{T}})$,
	\item $\mindeg(\mathscr{T})=
	\left\lfloor\tfrac{m+e-l-2}{e}\right\rfloor
	+\left\lfloor\tfrac{l+1}{e}\right\rfloor
	+\left\lfloor\tfrac{m}{e}\right\rfloor
	+\min(A_{\mathscr{T}})$.
\end{itemize}

We now set
\begin{align*}
a_n:=&\ \#\{i\mid{1\leqslant{i}\leqslant{n},
	i\equiv{l+1}\Mod{e}}\},\\
b_n:=&\ \#\{i\mid{1\leqslant{i}\leqslant{n},
	i\equiv{l+2}\Mod{e}}\},\\
c_n:=&\ \#\{i\mid{1\leqslant{i}\leqslant{n},
	i-l\not\equiv{1,2}\Mod{e}}\}.
\end{align*}

\begin{rmk}\label{rmk:abc}
	The values of $a_n$, $b_n$ and $c_n$ in each of the cases given in \Cref{sec:UDN} are as follows.
	\begin{itemize}
		\item Case I: $a_n=b_n=\lfloor\frac{n-l-1}{e}\rfloor+1$ and $c_n=n-2\lfloor\frac{n-l-1}{e}\rfloor-2$,
		\item Case II: $a_n=\lfloor\frac{n-l-1}{e}\rfloor + 1$, $b_n=\lfloor\frac{n-l-1}{e}\rfloor$, $c_n=n-2\lfloor\frac{n-l-1}{e}\rfloor -1$,
		\item Case III: $a_n=\lfloor\frac{n}{e}\rfloor$, $b_n=\lfloor\frac{n}{e}\rfloor + 1$, $c_n=n-2\lfloor\frac{n}{e}\rfloor -1$,
		\item Cases IV: $a_n=b_n=\lfloor\frac{n}{e}\rfloor$ and $c_n=n-2\lfloor\frac{n}{e}\rfloor$.
	\end{itemize}
\end{rmk}

\begin{lem}\label{lem:maxmin}
Let $\mathscr{T}=\std((n-m),(1^m))$ and $1\leqslant m<n$.
\begin{enumerate}
\item{If $1\leqslant{m}\leqslant{\tfrac{n}{e}}$, then $\max(A_{\mathscr{T}})=m$ and $\min(A_{\mathscr{T}})=-m$.}
\item{If $\tfrac{n}{e}<m<n-\tfrac{n}{e}$, then $\max(A_{\mathscr{T}})=a_n$ and $\min(A_{\mathscr{T}})=-b_n$.}
\item{If $n-\tfrac{n}{e}\leqslant{m}<n$, then $\max(A_{\mathscr{T}})=n-m+a_n-b_n$ and $\min(A_{\mathscr{T}})=m-n+a_n-b_n$.}
\end{enumerate}
\end{lem}

\begin{proof}
Let $\tts,\ttt\in\mathscr{T}$ be such that $\deg(\ttt)=\maxdeg(\mathscr{T})$ and $\deg(\tts)=\mindeg(\mathscr{T})$.
It follows from \cref{lem:deg} that $\ttt$ (respectively, $\tts$) is a standard $((n-m),(1^m))$-tableau with the maximum (resp., minimum) number of entries congruent to $l+1$ modulo $e$, say $i_{\ttt}$ (resp., $i_{\tts}$), together with the minimum (resp., maximum) number of entries congruent to $l+2$ modulo $e$, say $j_{\ttt}$ (resp., $j_{\tts}$), which lie in the leg of $\ttt$ (resp., $\tts$). We then compute $\deg(\tts)=i_{\tts}-j_{\tts}$ and $\deg(\ttt)=i_{\ttt}-j_{\ttt}$.
\end{proof}

\begin{prop}\label{prop:grdimS}
Let $1\leqslant m\leqslant n$ and $\mathscr{T}=\std((n-m),(1^m))$. Then $\grdim \left(S_{((n-m),(1^m))}\right)$ is
\[
\sum_{i=0}^{\mx\left(A_{\mathscr{T}}\right)-\min\left(A_{\mathscr{T}}\right)}
\left(
\sum_{j=0}^{\mx\left(A_{\mathscr{T}}\right)}
\left(
\binom{a_n}{m-i+j}
\binom{b_n}{j}
\binom{c_n}{i-2j}
\right)
v^{\left(
\mx\left(A_{\mathscr{T}}\right)
-i+\left\lfloor\frac{m}{e}\right\rfloor+
\left\lfloor\frac{m+e-l-2}{e}\right\rfloor
+\left\lfloor\frac{l+1}{e}\right\rfloor
\right)}
\right).
\]
\end{prop}

\begin{proof}
Let $\ttt\in\mathscr{T}$. By \cref{lem:deg}, there are at most $\max(A_{\mathscr{T}})$ entries in the leg of $\ttt$ congruent to $l+1$ modulo $e$, and at most $\min(A_{\mathscr{T}})$ entries congruent to $l+2$ modulo $e$. Thus, there exists a tableau with degree
\[d_i:=\max(A_{\mathscr{T}})-i+\left\lfloor\tfrac{m}{e}\right\rfloor+
\left\lfloor\tfrac{m+e-l-2}{e}\right\rfloor+\left\lfloor\tfrac{l+1}{e}\right\rfloor\]
for all $i\in\{0,\dots,\max(A_{\mathscr{T}})-\min(A_{\mathscr{T}})\}$, and hence $\grdim(S_{((n-m),(1^m))})$ has $\max(A_{\mathscr{T}})-\min(A_{\mathscr{T}})+1$ terms.

Suppose that $\ttt$ has degree $d_i$ for some $i$ and that there are $j$ entries congruent to $l+2$ modulo $e$ in the leg of $\ttt$. These $j$ entries contribute $-j$ to the degree of $\ttt$. Hence, there must be $m-i+j$ entries congruent to $l+1$ modulo $e$ in the leg of $\ttt$, and the remaining $i-2j$ nodes in the leg of $\ttt$ must contain entries congruent to neither $l+1$ modulo $e$ nor $l+2$ modulo $e$. Thus, there are
$\binom{a_n}{m-i+j}
\binom{b_n}{j}
\binom{c_n}{i-2j}$
standard $((n-m),(1^m))$-tableaux with this combination of entries in its leg for some $j\in\left\{0,\dots,\left\lfloor\frac{i}{2}\right\rfloor\right\}$, and summing over $j$ gives the number of standard $((n-m),(1^m))$-tableaux with degree $d_i$.
\end{proof}

Later on, we will require the explicit leading and trailing terms in the graded dimensions of Specht modules labelled by hook bipartitions as given below.

\begin{cor}\label{cor:grdim}
Let $1\leqslant m\leqslant n$ and $x=\lfloor\frac{m}{e}\rfloor+
\lfloor\frac{m+e-l-2}{e}\rfloor
+\lfloor\frac{l+1}{e}\rfloor$.
Then the first and last two terms in the graded dimension of $S_{((n-m),(1^m))}$ are displayed in the following table.
\begin{center}
\begin{tabular}{ c | | c | c | c }
&  $1\leqslant{m}\leqslant{\frac{n}{e}}$ &  $\frac{n}{e}<m<n-\frac{n}{e}$ &  $n-\frac{n}{e}\leqslant{m}<n$ \\
	\hline\hline 
\rule{0pt}{4ex} $\nth{1}$ term
&\rule{0pt}{4ex} $\binom{a_n}{m}v^{\left(	m+x	\right)}$
 &\rule{0pt}{4ex}  $\binom{c_n}{m-a}v^{\left(a_n+x\right)}$ 
 &\rule{0pt}{4ex}  $\binom{b_n}{n-m} v^{\left(n-m+a_n-b_n+x\right)}$ \\
\rule{0pt}{4ex} $\nth{2}$ term 
&\rule{0pt}{4ex} $c_n\binom{a_n}{m-1}v^{\left(m-1+x\right)}$
&\rule{0pt}{4ex}  $\left(a_n\binom{c_n}{m-a_n+1}+b_n\binom{c_n}{m-a_n-1}\right)v^{\left(a_n-1+x\right)}$ 
&\rule{0pt}{4ex}  $c_n\binom{b_n}{n-m-1} v^{\left(n-m+a_n-b_n-1+x\right)}$ \\  
\rule{0pt}{4ex} $\nth{2}$ last term
&\rule{0pt}{4ex} $c_n\binom{b_n}{m-1}v^{\left(1-m+x\right)}$
 &\rule{0pt}{4ex}  $\left( b_n\binom{c_n}{m-b_n+1}+ a_n\binom{c_n}{m-b_n-1}\right) v^{\left(1-b_n+x\right)}$
 &\rule{0pt}{4ex}  $c_n\binom{a_n}{n-m-1} v^{\left(1+m-n+a_n-b_n+x\right)}$\\
\rule{0pt}{4ex} last term
&\rule{0pt}{4ex} $\binom{b_n}{m}v^{\left(-m+x\right)}$
&\rule{0pt}{4ex}  $\binom{c_n}{m-b_n}v^{\left(-b_n+x\right)}$
 &\rule{0pt}{4ex}  $\binom{a_n}{n-m}v^{\left(m-n+a_n-b_n+x\right)}$    
\end{tabular}
\end{center}
\end{cor}

\section{Graded dimensions of the composition factors of $S_{((n-m),(1^m))}$}\label{sec:GDI}

We now study the \emph{graded} composition factors of Specht modules labelled by hook bipartitions, and determine the leading terms in their graded dimensions. Our results rely on the basis elements that span these irreducible $\mathscr{H}_n^{\Lambda}$-modules, which we deduce from the spanning sets of the images and the kernels of certain Specht module homomorphisms given in~\cite[Lemma 5.10]{Sutton17I}.

Recalling from \cref{thm:selfdual} that irreducible $\mathscr{H}_n^{\Lambda}$-modules are self-dual as graded modules, leads us to the following.

\begin{prop}\label{prop:irrsymmetric}
	Let $\lambda\in\mathscr{RP}_n^l$. Then $\grdim\left(D_{\lambda}\right)$ is symmetric in $v$ and $v^{-1}$.
\end{prop}

Thus, by the symmetry of the graded dimensions of irreducible $\mathscr{H}_n^{\Lambda}$-modules, we automatically recover their trailing terms if we know their leading terms. Together with \cref{defn:grdim}, the following result is an immediate consequence.

\begin{cor}\label{lem:irr}
	Let $\lambda\in\mathscr{P}_n^l$ and $\mathscr{T}\subseteq\std(\lambda)$.
	Suppose that $M$ is an irreducible $\mathscr{H}_n^{\Lambda}$-module with spanning set $\spn\{v_{\ttt}\mid\ttt\in\mathscr{T}\}$ such that $M\cong D_{\mu}$ as ungraded $\mathscr{H}_n^{\Lambda}$-modules for some $\mu\in\mathscr{RP}_n^l$. Then
	\[
	\grdim (D_{\lambda})=
	v^i
	\sum_{\ttt\in\mathscr{T}}v^{\degr (\ttt)}
	\in\mathbb{N}\cup\{0\}[v+v^{-1}],
	\]
	where $2i=-\maxdeg (\mathscr{T})-\mindeg (\mathscr{T})$. Moreover, the highest degree in the graded dimension of $D_{\mu}$ is
	$\frac{1}{2}(\maxdeg (\mathscr{T})-\mindeg (\mathscr{T}))$.
\end{cor}

\subsection[Case I: $\kappa_2\not\equiv{\kappa_1-1}\Mod{e}$ and $n\not\equiv{l+1}\Mod{e}$]{Case I: $\kappa_2\not\equiv{\kappa_1-1}\Mod{e}$ and $n\not\equiv{l+1}\Mod{e}$}

We recall from \cref{factors} that $S_{((n-m),(1^m))}$ is irreducible in this case, and moreover, we know that ${S_{((n-m),(1^m))}}\cong{D_{\mu_{n,m}}}\langle{i}\rangle$ as graded $\mathscr{H}_n^{\Lambda}$-modules for some $i\in\mathbb{Z}$.

\begin{prop}\label{prop:GDI1}
Suppose that  $\kappa_2\not\equiv{\kappa_1-1}\Mod{e}$ and $n\not\equiv{l+1}\Mod{e}$, and let $1\leqslant m <n$. Then the leading term of $\grdim\left(D_{\mu_{n,m}}\right)$ is
\begin{enumerate}
\item{ $\displaystyle{\binom{\left\lfloor\tfrac{n-l-1}{e}\right\rfloor+1}{m}v^m}$ if $1\leqslant{m}\leqslant{\tfrac{n}{e}+1}$,}
\item{ $\displaystyle{\binom{n-2\left(\left\lfloor\tfrac{n-l-1}{e}\right\rfloor+1\right)}{m-\left\lfloor\tfrac{n-l-1}{e}\right\rfloor-1}}
v^{\left\lfloor\tfrac{n}{e}\right\rfloor}$ if $\tfrac{n}{e}+1<m<n-\tfrac{n}{e}-1$, }
\item{ $\displaystyle{\binom{\left\lfloor\tfrac{n-l-1}{e}\right\rfloor+1}{n-m}v^{n-m}}$ if $n-\tfrac{n}{e}-1\leqslant{m} < n$. }
\end{enumerate}

Moreover, $D_{\mu_{n,m}}\left\langle{\left\lfloor\tfrac{m}{e}\right\rfloor
	+\left\lfloor\tfrac{m+e-l-2}{e}\right\rfloor}\right\rangle\cong{S_{((n-m),(1^m))}}$ as graded $\mathscr{H}_n^{\Lambda}$-modules.
\end{prop}

\begin{proof}
Since $S_{((n-m),(1^m))}$ is irreducible, the coefficients of the leading terms in $\grdim\left(D_{\mu_{n,m}}\right)$ and $\grdim\left(S_{((n-m),(1^m))}\right)$ are equal, which we know from \cref{cor:grdim}.

Let $\mathscr{T} = \std {((n-m),(1^m))}$. If $1\leqslant{m}\leqslant{\tfrac{n}{e}+1}$, then $\maxdeg (\mathscr{T})=m+\lfloor\tfrac{m}{e}\rfloor+\lfloor\tfrac{m+e-l-2}{e}\rfloor$ and $\mindeg (\mathscr{T})=-m+\lfloor\tfrac{m}{e}\rfloor+\lfloor\tfrac{m+e-l-2}{e}\rfloor$, by \cref{lem:maxmin}. It thus follows from \cref{lem:irr} that the highest degree in the graded dimension of $D_{((n-m),(1^m))}$ is $\tfrac{1}{2}\left(\maxdeg (\mathscr{T})-\mindeg (\mathscr{T})\right)=m$. Similarly, one can deduce the leading degrees in the other two cases.

We now determine $i\in\mathbb{Z}$ such that $D_{((n-m),(1^m))}\cong{S_{((n-m),(1^m))}\langle{i}\rangle}$ as graded $\mathscr{H}_n^{\Lambda}$-modules.
By above, we also know from \cref{lem:irr} that,   for all $m\in\{1,\dots,n-1\}$, $i=-\tfrac{1}{2}\maxdeg (\mathscr{T})-\tfrac{1}{2}\mindeg (\mathscr{T})=-\lfloor\tfrac{m}{e}\rfloor-\lfloor\tfrac{m+e-l-2}{e}\rfloor$, as required.
\end{proof}

\begin{ex}
	
Let $e=3$, $\kappa=(0,1)$. The following tableaux index the basis vectors of $S_{((2),(1^2))}$
\[
\ttt_1=\young(34,,1,2)\quad
\ttt_2=\young(24,,1,3)\quad
\ttt_3=\young(23,,1,4)\quad
\ttt_4=\young(14,,2,3)\quad
\ttt_5=\young(13,,2,4)\quad
\ttt_6=\young(12,,3,4)
\]
It is easy to check that $\deg(\ttt_1)=\deg(\ttt_5)=1$, $\deg(\ttt_2)=\deg(\ttt_6)=-1$ and $\deg(\ttt_3)=\deg(\ttt_4)=0$. Hence $\grdim\left(S_{((2),(1^2))}\right)=2v+2+2v^{-1}$ is symmetric in $v$ and $v^{-1}$, and thus $S_{((2),(1^2))}\cong D_{\mu_{4,2}} = D_{((2),(1^2))}$ as graded $\mathscr{H}_n^{\Lambda}$-modules.	
\end{ex}

\subsection[Case II: $\kappa_2\not\equiv{\kappa_1-1}\Mod{e}$ and $n\equiv{l+1}\Mod{e}$]{Case II: $\kappa_2\not\equiv{\kappa_1-1}\Mod{e}$ and $n\equiv{l+1}\Mod{e}$}

For $m\in\{1,\dots,n-1\}$, we recall from \cref{thm:labels case 2} that $S_{((n-m),(1^m))}$ has graded composition factors ${D_{\mu_{n,m-1}}}$ and ${D_{\mu_{n,m}}}$ such that ${D_{\mu_{n,m-1}}}\cong{\im (\gamma_{m-1})\langle{i}\rangle}$ and ${D_{\mu_{n,m}}}\cong{\im (\gamma_m)\langle{j}\rangle}$ for some $i,j\in\mathbb{Z}$.

\begin{prop}\label{grdimD1}
Suppose that $\kappa_2\not\equiv{\kappa_1-1}\Mod{e}$ and $n\equiv{l+1}\Mod{e}$, and let $1\leqslant m<n$. Then the leading term of $\grdim\left(D_{\mu_{n,m}}\right)$ is
\begin{enumerate}
\item{
$\displaystyle{\binom{\left\lfloor\frac{n-l-1}{e}\right\rfloor}
{m}v^m}$
if $1\leqslant{m}\leqslant\frac{n}{e}$,
}
\item{
$\displaystyle{
\binom{n-2\left\lfloor\frac{n-l-1}{e}\right\rfloor-1}{m-\left\lfloor\frac{n-l-1}{e}\right\rfloor}
v^{\left\lfloor\frac{n-l-1}{e}\right\rfloor}
}$
if $\frac{n}{e}<m<n-\frac{n}{e}$,
}
\item{
$\displaystyle{
\binom{\left\lfloor\frac{n-l-1}{e}\right\rfloor}{n-m-1}v^{n-m-1}
}$
if $\frac{n}{e}\leqslant{m} < n$.
}
\end{enumerate}

Moreover, $D_{\mu_{n,m}}\left\langle{\left\lfloor\tfrac{m+e-l-2}{e}\right\rfloor
	+\left\lfloor\tfrac{m}{e}\right\rfloor}\right\rangle\cong\im (\gamma_m)$ as graded $\mathscr{H}_n^{\Lambda}$-modules.
\end{prop}

\begin{proof}
Let $\mathscr{T}=\{\ttt\in\std((n-m),(1^m))\mid\ttt(1,n-m,1)=n\}$. Then we know from~\cite[Lemma 5.10]{Sutton17I} that the set of vectors $\{v_{\ttt}\mid\ttt\in\mathscr{T}\}$ spans $\im (\gamma_m)$. By \cref{lem:irr}, we have
\[
\grdim \left(D_{\mu_{n,m}}\right)
=v^i\grdim \left(\im (\gamma_m)\right)
=v^i\sum_{\ttt\in\mathscr{T}}v^{\deg (\ttt)},
\]
where $2i=-\maxdeg (\mathscr{T})-\mindeg (\mathscr{T})$, and moreover, we know from \cref{prop:irrsymmetric} that the coefficients in the leading and trailing terms of the graded dimension of $D_{\mu_{n,m}}$ are equal. We now recall from \cref{rmk:abc} that $a_n=\left\lfloor\frac{n-l-1}{e}\right\rfloor + 1$, $b_n=\left\lfloor\frac{n-l-1}{e}\right\rfloor$ and $c_n=n-2\left\lfloor\frac{n-l-1}{e}\right\rfloor - 1$, and suppose that $\tts,\ttt\in\mathscr{T}$ are such that $\degr(\tts)=\maxdeg(\mathscr{T})$ and $\degr(\ttt)=\mindeg(\mathscr{T})$. Then the proof follows similarly to that of \cref{lem:maxmin}, by applying \cref{lem:deg}. We note that $n$, which is congruent to $l+1$ modulo $e$, lies in the hand node of both $\tts$ and $\ttt$.

\begin{enumerate}
	\item Observe that each node in the leg of $\tts$ contains one of the $a-1$ entries congruent to $l+1$ modulo $e$ (excluding $n$), and hence there are $\binom{a_n-1}{m}$ standard $((n-m),(1^m))$-tableaux with degree $\degr(\tts)$. We now observe that each node in the leg of $\ttt$ contains one of the $b_n$ entries congruent to $l+2$ modulo $e$. Hence $\mx(A_{\mathscr{T}})=m$ and $\min(A_{\mathscr{T}})=-m$.
	
	\item Firstly, the leg of $\tts$ contains all of the remaining $a_n-1$ entries congruent to $l+1$ modulo $e$ and $m-a_n+1$ of the $c_n$ entries neither congruent to $l+1$ modulo $e$ nor congruent to $l+2$ modulo $e$. Secondly, the leg of $\ttt$ contains all of the $b_n$ entries congruent to $l+2$ modulo $e$, and $m-b_n$ of the $c_n$ entries congruent to neither $l+1$ modulo $e$ nor $l+2$ modulo $e$. Hence $\mx(A_{\mathscr{T}})=a_n-1$ and $\min(A_{\mathscr{T}})=-b_n$.
	
	\item Except for the hand nodes of $\tts$ and $\ttt$, we observe that each node in the arm of $\tts$ contains one of the $b_n$ entries congruent to $l+2$ modulo $e$, and that every node in the arm of $\ttt$ contains one of the remaining $a-1$ entries congruent to $l+1$ modulo $e$. Hence $\mx(A_{\mathscr{T}})=n-m-1$ and $\min(A_{\mathscr{T}})=m-n+1$.
	
\end{enumerate}

For all $m$, we notice that $\min (A_{\mathscr{T}})=-\max (A_{\mathscr{T}})$. Moreover, $i=-\frac{1}{2}\maxdeg (\mathscr{T}) -\frac{1}{2} \mindeg (\mathscr{T}) = -\lfloor\tfrac{m+e-l-2}{e}\rfloor-\lfloor\tfrac{m}{e}\rfloor$, as required.
\end{proof}

\begin{ex}
	Let $e=3$, $\kappa=(0,0)$, $n=7$ and $\mathscr{T}=\{\ttt\in\std((5),(1^2))\mid\ttt(2,1,2)=7\}$. We know from~\cite[Lemma 5.10]{Sutton17I} that $\im (\gamma_1)$ is spanned by $\{v_{\ttt}\mid\ttt\in\mathscr{T}\}$. The tableaux lying in $\mathscr{T}$ are
	\[
	\ttt_1=\gyoung(23456,,1,!\gr7)\quad
	\ttt_2=\gyoung(13456,,2,!\gr7)\quad
	\ttt_3	=\gyoung(12456,,3,!\gr7)
	\]
	\[
	\ttt_4=\gyoung(12356,,4,!\gr7)\quad
	\ttt_5=\gyoung(12346,,5,!\gr7)\quad
	\ttt_6=\gyoung(12345,,6,!\gr7)
	\]
	Let $\tts,\ttt\in\mathscr{T}$ be such that $\deg (\tts)=\maxdeg (\mathscr{T})$ and $\deg (\ttt)=\mindeg (\mathscr{T})$. Then, by \cref{lem:entries}, $\{1,4\}\in\tts(1,1,2)$ and $\{2,5\}\in\ttt(1,1,2)$. Hence $\deg (\ttt_1)=\deg \left(\ttt_4\right)>\deg \left(\ttt_3\right)=\deg \left(\ttt_6\right)>
	\deg \left(\ttt_2\right)=\deg \left(\ttt_5\right)$. One can check that
	$\deg (\ttt_1)=3$, $\deg \left(\ttt_2\right)=1$ and $\deg \left(\ttt_3\right)=2$, obtaining $\grdim \left(\im (\gamma_1)\right)	=2v^3+2v^2+2v$.	By \cref{thm:labels case 2}, $\im (\gamma_1)\cong{D_{\mu_{7,1}}}=D_{((6),(1))}$ as ungraded $\mathscr{H}_7^{\Lambda}$-modules, and thus by shifting the grading of $\im (\gamma_1)$, we have 
	\[\grdim \left(D_{((6),(1))}\right)
	=\grdim \left(\im (\gamma_1)\langle{-2}\rangle\right)
	=2v+2+2v^{-1}.\]
\end{ex}

\subsection[Case III: $\kappa_2\equiv{\kappa_1-1}\Mod{e}$ and $n\not\equiv{0}\Mod{e}$]{Case III: $\kappa_2\equiv{\kappa_1-1}\Mod{e}$ and $n\not\equiv{0}\Mod{e}$}

We recall from \cref{labels:case 3} that $S_{((n-m),(1^m))}$ has graded composition factors ${D_{\mu_{n,2m}}}$ and ${D_{\mu_{n,2m+1}}}$, for all $m\in\{1,\dots,n-1\}$, such that ${D_{\mu_{n,2m}}}\cong{\text{im}(\chi_m)\langle{i}\rangle}$ and ${D_{\mu_{n,2m+1}}}\cong{(S_{((n-m),(1^m))}/\im (\chi_m))\langle{j}\rangle}$ for some $i,j\in\mathbb{Z}$.

\begin{prop}\label{grdimD4}
Suppose that $\kappa_2\equiv{\kappa_1-1}\Mod{e}$ and $n\not\equiv{0}\Mod{e}$, and let $1\leqslant m<n$.
\begin{enumerate}
\item{
Then the leading term of $\grdim (D_{\mu_{n,2m}})$ is
	\begin{enumerate}
	\item{
	$\displaystyle{\binom{\left\lfloor\frac{n}{e}\right\rfloor}{m}v^m}$
	if $1\leqslant{m}\leqslant{\frac{n}{e}}$,	
	}
	\item{
	$\displaystyle{\binom{n-2\left\lfloor\frac{n}{e}\right\rfloor-1}{m-\left\lfloor\frac{n}{e}\right\rfloor}
		v^{\left\lfloor\frac{n}{e}\right\rfloor}}$
	if $\frac{n}{e}<m<n-\frac{n}{e}$,
	}
	\item{
	$\displaystyle{
	\binom{\left\lfloor\frac{n}{e}\right\rfloor}{n-m-1}
	v^{\left(n-m-1\right)}	
	}$	
	if $n-\frac{n}{e}\leqslant{m}<n$.
	}
	\end{enumerate}
	
Moreover, $D_{\mu_{n,2m}}\left\langle{\left\lfloor\tfrac{m}{e}\right\rfloor+\left\lfloor\tfrac{m-1}{e}\right\rfloor-1}\right\rangle\cong\im (\chi_m)$ as graded $\mathscr{H}_n^{\Lambda}$-modules.
}
\item{
Then the leading term of $\grdim (D_{\mu_{n,2m+1}})$ is
	\begin{enumerate}
	\item{
	$\displaystyle{
	\binom{\left\lfloor\frac{n}{e}\right\rfloor}{m-1}
	v^{m-1}
	}$	
	if $1\leqslant{m}\leqslant{\frac{n}{e}}$,
	}
	\item{
	$\displaystyle{
	\binom{n-2\left\lfloor\frac{n}{e}\right\rfloor-1}
	{m-1-\left\lfloor\frac{n}{e}\right\rfloor}
	v^{\left\lfloor\frac{n}{e}\right\rfloor}
	}$	
	if $\frac{n}{e}<m<n-\frac{n}{e}$,
	}
	\item{
	$\displaystyle{
	\binom{\left\lfloor\frac{n}{e}\right\rfloor}{n-m}
	v^{n-m}	
	}$	
	if $n-\frac{n}{e}\leqslant{m}<n$.
	}
	\end{enumerate}
	
Moreover, $D_{\mu_{n,2m+1}}\left\langle{\left\lfloor\tfrac{m-1}{e}\right\rfloor+\left\lfloor\tfrac{m}{e}\right\rfloor}\right\rangle\cong{S_{((n-m),(1^m))}/\im (\chi_m)}$ as graded $\mathscr{H}_n^{\Lambda}$-modules.
	}
\end{enumerate}
\end{prop}

\begin{proof}
The proof follows the same structure as that of \cref{grdimD1}, using the spanning sets of $\im (\chi_m)$ and $S_{((n-m),(1^m))}/\im (\chi_m)$ determined from~\cite[Lemma 5.10]{Sutton17I}. In particular, we apply \cref{lem:irr} with $\lambda=\mu_{n,2m}$ and $\mathscr{T}=\{\ttt\in\std((n-m),(1^m))\mid\ttt(1,1,1)=1\}$ for the first part, and with $\lambda=\mu_{n,2m+1}$ and $\mathscr{T}=\{\ttt\in\std((n-m),(1^m))\mid\ttt(1,1,2)=1\}$ for the second part.
\end{proof}

\begin{ex}
	Let $e=3$, $\kappa=(0,2)$, $n=5$ and $\mathscr{T}=\{\ttt\in\std ((3),(1^2))\mid\ttt(1,1,1)=1\}$. By~\cite[Lemma 5.10]{Sutton17I}, $\im (\chi_2)$ is spanned by $\{v_{\ttt}\mid\ttt\in\mathscr{T}\}$. There are six tableaux in $\mathscr{T}$, namely
	\[
	\ttt_1=\gyoung(!\gr1!\wh45,,2,3)\quad
	\ttt_2=\gyoung(!\gr1!\wh35,,2,4)\quad
	\ttt_3=\gyoung(!\gr1!\wh34,,2,5)\quad
	\ttt_4=\gyoung(!\gr1!\wh25,,3,4)\quad
	\ttt_5=\gyoung(!\gr1!\wh24,,3,5)\quad
	\ttt_6=\gyoung(!\gr1!\wh23,,4,5)
	\]
	One can check from \cref{lem:deg} that $\deg (\ttt_1)=\deg (\ttt_5)=2$, $\deg (\ttt_2)=\deg (\ttt_6)=0$ and $\deg (\ttt_3)=\deg (\ttt_4)=1$, and hence $\grdim \left(\im (\chi_2)\right)=2v^2+2v+2$. By \cref{labels:case 3}, $\im (\chi_2)\cong{D_{\mu_{5,4}}}=D_{((3,1^2),\varnothing)}$ as ungraded $\mathscr{H}_5^{\Lambda}$-modules, and by shifting the degree of $\im (\chi_2)$, we have
	\[
	\grdim \left(D_{((3,1^2),\varnothing)}\right)
	=\grdim \left(\im (\chi_2)\langle{-1}\rangle\right)
	=2v+2+2v^{-1}.\]
	
	Let $\mathscr{S}=\{\tts\in\std ((3),(1^2))\mid\tts(1,1,2)=1\}$. 
	It follows from above that $S_{((3),(1^2))}/\im (\chi_2)$ is spanned by $\{v_{\tts}\mid\tts\in\mathscr{S}\}$, and moreover, we know that $S_{((3),(1^2))}/\im (\chi_2)\cong{D_{\mu_{5,5}}}=D_{((3),(1^2))}$ as ungraded $\mathscr{H}_5^{\Lambda}$-modules by \cref{labels:case 3}. We see that $\mathscr{S}$ contains the following tableaux
	\[
	\tts_1=\young(345,,!\gr1,!\wh2)\quad
	\tts_2=\young(245,,!\gr1,!\wh3)\quad
	\tts_3=\young(235,,!\gr1,!\wh4)\quad
	\tts_4=\young(234,,!\gr1,!\wh5)
	\]
	One can easily check that $\deg (\tts_1)=\deg (\tts_4)=0$, $\deg (\tts_2)=1$ and $\deg (\tts_3)=-1$. Thus
	\[ \grdim \left(D_{((3),(1^2))}\right) = \grdim \left(S_{((3),(1^2))}/\im (\chi_2)\right) = v+2+v^{-1},\]
	and $S_{((3),(1^2))}/\im (\chi_2)\cong{D_{((3),(1^2))}}$ as graded $\mathscr{H}_5^{\Lambda}$-modules.
\end{ex}

\subsection[Case IV: $\kappa_2\equiv{\kappa_1-1}\Mod{e}$ and $n\equiv{0}\Mod{e}$]{Case IV: $\kappa_2\equiv{\kappa_1-1}\Mod{e}$ and $n\equiv{0}\Mod{e}$}

Let $1<m<n$. Then we know from \cref{labels:case 4} that $S_{((n-m),(1^m))}$ has graded composition factors $D_{\mu_{n,2m}}$ and $D_{\mu_{n,2m+1}}$ such that ${D_{\mu_{n,2m}}}\langle{i}\rangle\cong{\im (\phi_m)}$ and ${D_{\mu_{n,2m+1}}}\langle{j}\rangle\cong{\ker (\gamma_m)/\im (\phi_m)}$ for some $i,j\in\mathbb{Z}$. Except for the one-dimensional Specht modules, recall from \cref{labels:case 4} that $S_{((n-m),(1^m))}$ has either three or four composition factors. Hence we not only find the leading terms of $\grdim\left(D_{\mu_{n,2m}}\right)$ and $\grdim\left(D_{\mu_{n,2m+1}}\right)$, but the second leading terms too. It will become apparent to the reader in \cref{sec:GDNum} that these extra terms are, in fact, necessary in order to determine the corresponding graded decomposition numbers in this case.

\begin{prop}\label{grdimD2}
Suppose that $\kappa_2\equiv{\kappa_1-1}\Mod{e}$ and $n\equiv{0}\Mod{e}$, and let $1< m<n$.
\begin{enumerate}
\item{Then the first two leading terms of $\grdim (D_{\mu_{n,2m}})$ are
	\begin{enumerate}
	\item{
	$\displaystyle{\binom{\frac{n-e}{e}}{m-1}v^{m-1}}$ and 
	$\displaystyle{
	\frac{(e-2)n}{e}
	\binom{\frac{n-e}{e}}{m-2}v^{m-2}}$	
	 if $1<{m}\leqslant\frac{n}{e}$,
	}
	\item{
	$\displaystyle{\binom{\frac{(e-2)n}{e}}{\frac{em-n}{e}}v^{\frac{n-e}{e}}}$ and 
	$\displaystyle{
	\frac{n-e}{e}
	\left(
	\binom{\frac{(e-2)n}{e}}{m-\frac{n}{e}+1}+
	\binom{\frac{(e-2)n}{e}}{m-\frac{n}{e}-1}
	\right)}v^{\frac{n-2e}{e}}$
	if $\frac{n}{e}<m\leqslant\frac{n(e-1)}{e}$,	
	}
	\item{
	$\displaystyle{\binom{\frac{n-e}{e}}{n-m-1} v^{n-m-1}}$ and 
	$\displaystyle{\frac{(e-2)n}{e}\binom{\frac{n-e}{e}}{n-m-2} v^{n-m-2}}$
	if $\frac{n(e-1)+e}{e}\leqslant{m}< n$.	
	}
	\end{enumerate}
	
Moreover, $D_{\mu_{n,2m}}\left\langle{\left\lfloor\tfrac{m-1}{e}\right\rfloor+\left\lfloor\tfrac{m}{e}\right\rfloor+2}\right\rangle\cong\im (\phi_m)$ as graded $\mathscr{H}_n^{\Lambda}$-modules.
}
\item{
Then the first two leading terms of $\grdim (D_{\mu_{n,2m+1}})$ are
	\begin{enumerate}
	\item{
	$\displaystyle{\binom{\frac{n-e}{e}}{m-2}v^{m-2}}$ and 
	$\displaystyle{\frac{(e-2)n}{e}
	\binom{\frac{n-e}{e}}{m-3}v^{m-3}}$
	if $1<{m}\leqslant{\frac{n}{e}}$,	
	}
	\item{
	$\displaystyle{\binom{\frac{(e-2)n}{e}}{\frac{e(m-1)-n}{e}}v^{\frac{n-e}{e}}}$ and 
	$\displaystyle{
	\frac{n-e}{e}
	\left(
	\binom{\frac{(e-2)n}{e}}{\frac{em-n}{e}}
	+
	\binom{\frac{(e-2)n}{e}}{\frac{e(m-2)-n}{e}}
	\right)	
	v^{\frac{n-2e}{e}}
	}$
	if $\frac{n}{e}<m\leqslant\frac{n(e-1)}{e}$,	
	}
	\item{
	$\displaystyle{\binom{\frac{n-e}{e}}{n-m}v^{n-m}}$ and 
	$\displaystyle{\frac{(e-2)n}{e}
	\binom{\frac{n-e}{e}}{n-m-1}v^{n-m-1}}$
	if $\frac{n(e-1)+e}{e}\leqslant{m}< n$.	
	}
	\end{enumerate}
	
Moreover, $D_{\mu_{n,2m+1}}\left\langle{\left\lfloor\tfrac{m-1}{e}\right\rfloor+\left\lfloor\tfrac{m}{e}\right\rfloor+1}\right\rangle\cong\ker (\gamma_m)/\im (\phi_m)$ as graded $\mathscr{H}_n^{\Lambda}$-modules.
}
\end{enumerate}
\end{prop}

\begin{proof}
We follow the same structure as the proof of \cref{grdimD1}, using the spanning sets of $\im (\phi_m)$ and $\ker (\gamma_m)/\im (\phi_m)$ determined from \cite[Lemma 5.10]{Sutton17I}. In particular, we apply \cref{lem:irr} with $\lambda=\mu_{n,2m}$ and $\mathscr{T}=\{\ttt\in\std ((n-m),(1^m))\mid\ttt(1,1,1)=1,\ttt(m,1,2)=n\}$ for the first part, and with $\lambda=\mu_{n,2m+1}$ and $\mathscr{T}=\{\ttt\in\std ((n-m),(1^m))\mid\ttt(1,1,2)=1,\ttt(m,1,2)=n\}$ for the second part.
\end{proof}

\begin{ex} 
	Let $e=3$, $\kappa=(0,2)$, $n=6$ and $\mathscr{T}=\{\ttt\in\std ((3),(1^3))\mid\ttt(1,1,1)=1,\ttt(3,1,2)=6\}$. By \cite[Lemma 5.10]{Sutton17I}, $\im (\phi_3)$ is spanned by $\{v_{\ttt}\mid \ttt\in\mathscr{T}\}$. There are six tableaux in $\mathscr{T}$, namely
	\[
	\ttt_1=\gyoung(!\gr1!\wh45,,2,3,!\gr6)\quad
	\ttt_2=\gyoung(!\gr1!\wh35,,2,4,!\gr6)\quad
	\ttt_3=\gyoung(!\gr1!\wh34,,2,5,!\gr6)\quad
	\ttt_4=\gyoung(!\gr1!\wh25,,3,4,!\gr6)\quad
	\ttt_5=\gyoung(!\gr1!\wh24,,3,5,!\gr6)\quad
	\ttt_6=\gyoung(!\gr1!\wh23,,4,5,!\gr6)
	\]
	One can check that $\deg (\ttt_1)=\deg (\ttt_5)=4$, $\deg (\ttt_2)=\deg(\ttt_6)=2$ and $\deg (\ttt_3)=\deg (\ttt_4)=3$, and hence $\grdim \left(\im(\phi_3)\right)=2v^4+2v^3+2v^2$. We know from \cref{labels:case 4} that $\im (\phi_3)\cong{D_{\mu_{6,6}}}=D_{((4,1^2),\varnothing)}$ as ungraded $\mathscr{H}_6^{\Lambda}$-modules. Thus, by shifting the grading on $\im (\phi_3)$, we obtain
	\[
	\grdim\left(D_{((4,1^2),\varnothing)}\right)
	=\grdim \left(\im(\phi_3)\langle -3\rangle\right)=2v+2+2v^{-1}.
	\]
	
	Let $\mathscr{S}=\{\tts\in\std ((3),(1^3))\mid\tts(1,1,2)=1,\tts(3,1,2)=6\}$. By \cite[Lemma 5.10]{Sutton17I}, $\ker (\gamma_3)/\im (\phi_3)$ is spanned by $\{v_{\tts}\mid\tts\in\mathscr{S}\}$. There are four tableaux in $\mathscr{S}$, namely
	\[
	\tts_1=\gyoung(345,,!\gr1,!\wh2,!\gr6)\quad
	\tts_2=\gyoung(245,,!\gr1,!\wh3,!\gr6)\quad
	\tts_3=\gyoung(235,,!\gr1,!\wh4,!\gr6)\quad
	\tts_4=\gyoung(234,,!\gr1,!\wh5,!\gr6)
	\]
	One can check that $\deg (\tts_1)=\deg (\tts_4)=2$, $\deg (\tts_2)=3$ and $\deg (\tts_3)=1$, and hence $\grdim \left(\ker (\gamma_3)/\im (\phi_3)\right)=v^3+2v^2+v$. We know from \cref{labels:case 4} that $\ker (\gamma_3)/\im (\phi_3)\cong{D_{\mu_{6,7}}}=D_{((4),(1^2))}$ as ungraded $\mathscr{H}_6^{\Lambda}$-modules. By shifting the grading on $\ker (\gamma_3)/\im (\phi_3)$, we obtain
	\[
	\grdim\left( D_{((4),(1^2))}\right)
	=\grdim \left(\ker (\gamma_3)/\im (\phi_3)\langle -2 \rangle\right)=v+2+v^{-1}.
	\]
\end{ex}

\section{Graded decomposition numbers corresponding to $S_{((n-m),(1^m))}$}\label{sec:GDNum}

Recall that we determined the ungraded decomposition numbers for $\mathscr{H}_n^{\Lambda}$ corresponding to Specht modules labelled by hook bipartitions in \Cref{sec:UDN}, and then in \Cref{sec:GDS} and \Cref{sec:GDI}, we determined the graded dimensions of Specht modules labelled by hook bipartitions and of their composition factors, respectively. These findings are equivalent to solving part of the Decomposition Number Problem, corresponding to hook bipartitions, which we now provide an answer to.

Recall from \Cref{back:decomp} that the \emph{graded decomposition numbers} are defined to be the Laurent polynomials $[S_{\lambda}:D_{\mu}]_v=\sum_{i\in\mathbb{Z}}[S_{\lambda}:D_{\mu}\langle{i}\rangle]v^i$ for all $\lambda\in\mathscr{P}_n^l$ and for all $\mu\in\mathscr{RP}_n^l$.

We first determine the grading shifts on the trivial and sign representations to obtain the analogous graded representations. The trivial representation $S_{((n),\varnothing)}$ is generated by $v_{\ttt_{((n),\varnothing))}}$ where $\deg (\ttt_{((n),\varnothing))})=0$, so that $S_{((n),\varnothing)}={D_{((n),\varnothing)}}$ as graded $\mathscr{H}_n^{\Lambda}$-modules.
Hence
\[
\left[S_{((n),\varnothing)}:D_{\mu}\right]_v
=\begin{cases}
1&\text{ if $\mu=((n),\varnothing)$,}\\
0&\text{ otherwise.}
\end{cases}
\]
Recall from \Cref{simpsign} that $S_{(\varnothing,(1^n))}\cong{D_{(\varnothing,(1^n))^R}}$ as ungraded $\mathscr{H}_n^{\Lambda}$-modules. We now find $i\in\mathbb{Z}$ such that $S_{(\varnothing,(1^n))}\cong{D_{(\varnothing,(1^n))^R}\langle{i}\rangle}$ as graded $\mathscr{H}_n^{\Lambda}$-modules.

\begin{lem}
Let $\lambda=(\varnothing,(1^n))^R$. Then
\[
\left[S_{(\varnothing,(1^n))}:D_{\lambda}\right]_v=
\begin{cases}
v^{2\lfloor\frac{n}{e}\rfloor}
&\text{if $\kappa_2\equiv{\kappa_1-1}\Mod{e}$,}\\
v^{\left(\lfloor\frac{n}{e}\rfloor+\lfloor\frac{n-l-1}{e}\rfloor+1\right)}
&\text{if $\kappa_2\not\equiv{\kappa_1-1}\Mod{e}$.}
\end{cases}
\]
Moreover, $[S_{(\varnothing,(1^n))}:D_{\mu}]_v=0$ for all other $\mu\in\mathscr{RP}_n^2$.
\end{lem}

\begin{proof}
We have $[S_{(\varnothing,(1^n))}:D_{\lambda}]_v=v^{\deg \left(\ttt_{(\varnothing,(1^n))}\right)}$ since $\grdim \left(D_{\lambda}\right)=1$ and $\grdim \left(S_{(\varnothing,(1^n))}\right)=\deg (\ttt_{(\varnothing,(1^n))})$. 
We now deduce from the proof of \cref{lem:deg} that
\begin{align*}
\degr\left(\ttt_{(\varnothing,(1^n))}\right)
&=\left\lfloor\tfrac{n}{e}\right\rfloor
+\#\left\{i\ |\ \ttt_{(\varnothing,(1^n))}(i,1,2)\equiv{l+2}\Mod{e}\right\},
\end{align*}
where $\res(i,1,2)\equiv\kappa_1\Mod{e}$. In the leg of $[(\varnothing,(1^n))]$, we notice that there are $\left\lfloor\tfrac{n}{e}\right\rfloor$ $\kappa_1$-nodes if $\kappa_2\equiv\kappa_1-1\Mod{e}$ and $\left\lfloor\tfrac{n-l-1}{e}\right\rfloor+1$ $\kappa_1$-nodes otherwise, and we are done.
\end{proof}

For all regular bipartitions $\lambda\in\mathscr{RP}_n^2$, we now establish the \emph{graded} composition multiplicities $[S_{((n-m),(1^m))}:D_{\lambda}]_v$ of irreducible $\mathscr{H}_n^{\Lambda}$-modules $D_{\lambda}$ arising as composition factors of $S_{((n-m),(1^m))}$, for all $m\in\{1,\dots,n-1\}$, depending on whether $\kappa_2\equiv{\kappa_1-1}\Mod{e}$ or not and whether $n\equiv{l+1}\Mod{e}$ or not.

\subsection[Case I: $\kappa_2\not\equiv{\kappa_1-1}\Mod{e}$ and $n\not\equiv{l+1}\Mod{e}$]{Case I: $\kappa_2\not\equiv{\kappa_1-1}\Mod{e}$ and $n\not\equiv{l+1}\Mod{e}$}

Let $\kappa_2\not\equiv{\kappa_1-1}\Mod{e}$ and $n\not\equiv{l+1}\Mod{e}$. We recall from \cref{factors} that $S_{((n-m),(1^m))}$ is irreducible and isomorphic to $D_{\mu_{n,m}}$ as an ungraded $\mathscr{H}_n^{\Lambda}$-module for all $m\in\{1,\dots,n\}$. To find the graded multiplicity of $D_{\mu_{n,m}}$ arising as a composition factor of $S_{((n-m),(1^m))}$, it suffices to find the grading shift on $D_{\mu_{n,m}}$ so that it is isomorphic to $S_{((n-m),(1^m))}$ as a graded $\mathscr{H}_n^{\Lambda}$-module.

\begin{thm}
Suppose that $\kappa_2\not\equiv{\kappa_1-1}\Mod{e}$ and $n\not\equiv{l+1}\Mod{e}$, and let $\mu\in\mathscr{RP}_n^2$.
Then, for all $m\in\{1,\dots,n-1\}$, we have
    \[
	\left[S_{((n-m),(1^m)}:D_{\mu_{n,m}}\right]_v=
	\begin{cases}
	v^{\left(\left\lfloor\frac{m}{e}\right\rfloor+
	\left\lfloor\frac{m+e-l-2}{e}\right\rfloor\right)}
	&\text{if $\mu=\mu_{n,m}$,}\\
	0&\text{otherwise.}
	\end{cases}
	\]
\end{thm}

\begin{proof}
We determine $i\in\mathbb{Z}$ where $\left[S_{((n-m),(1^m))}:D_{\mu_{n,m}}\right]_v=v^i$, which is equivalent to finding $i\in\mathbb{Z}$ such that $S_{((n-m),(1^m))}\cong{D_{\mu_{n,m}}\langle{i}\rangle}$ as graded $\mathscr{H}_n^{\Lambda}$-modules. Thus, the result follows from \cref{prop:GDI1}.
\end{proof}

\begin{ex}
	Let $e=3$ and $\kappa=(0,0)$. Then the decomposition submatrix of $\mathscr{H}_6^{\Lambda}$ with rows corresponding to Specht modules labelled by hook bipartitions can be written as
	\[
	\begin{blockarray}{ccccccccccc}
	\begin{block}{c(ccccccc|ccc)}
	S_{((6),\varnothing)}   &1& & &   &   &   &   &   &   									&							\\[0.2em]
	S_{((5),(1))}           & &1& &   &      &  \text{\huge${0}$} &   						&&&			\\[0.6em]
	S_{((4),(1^2))}          && & v   &   &   &   &   									&&&				\\
	S_{((3),(1^3))}          & && &   v^2   &   &   &   									&&\text{\Huge${0}$}&		\\[0.6em]
	S_{((2),(1^4))}          & && &   &   v^2   &   &   									&&&					\\[0.3em]
	S_{((1),(1^5))}          & & \text{\huge{$0$}}&&   &   &  v^3  &   						&&&		\\[0.6em]
	S _{(\varnothing,(1^6))} & & &    &   &   &   &v^4									&&&		\\
	\end{block}
	\end{blockarray}
	\]
\end{ex}

\begin{comment}
\begin{ex}
	Let $e=3$ and $\kappa=(0,0)$. Then the decomposition submatrix of $\mathscr{H}_8^{\Lambda}$ with rows corresponding to Specht modules labelled by hook bipartitions is
	\[
	\begin{blockarray}{ccccccccccccc}
	\begin{block}{c(ccccccccc|ccc)}
	S_{((8),\varnothing)}   &1& & &   &   &   &   &   &   									&&&							\\[0.2em]
	S_{((7),(1))}           & &1& &   &   &   &   &  \text{\huge${0}$} &   						&&&			\\[0.6em]
	S_{((6),(1^2))}          && &v&   &   &   &   &   &   									&&&	\\[0.6em]
	S_{((5),(1^3))}          && & &v^2&   &   &   &   &   									&&&				\\
	S_{((4),(1^4))}          & && &   &v^2&   &   &   &   									&&\text{\Huge${0}$}&		\\[0.6em]
	S_{((3),(1^5))}          & && &   &   &v^3&   &   &   									&&&					\\[0.6em]
	S_{((2),(1^6))}          & & &&   &   &   &v^4&   &   									&&&							\\[0.2em]
	S_{((1),(1^7))}          & & \text{\huge{$0$}}&&   &   &   &   &v^4&   						&&&		\\[0.6em]
	S _{(\varnothing,(1^8))} & & & &  &   &   &   &   &v^5									&&&		\\
	\end{block}
	\end{blockarray}.
	\]
\end{ex}
\end{comment}

\subsection[Case II: $\kappa_2\not\equiv{\kappa_1-1}\Mod{e}$ and $n\equiv{l+1}\Mod{e}$]{Case II: $\kappa_2\not\equiv{\kappa_1-1}\Mod{e}$ and $n\equiv{l+1}\Mod{e}$}

Let $\kappa_2\not\equiv{\kappa_1-1}\Mod{e}$ and $n\equiv{l+1}\Mod{e}$. Recall from \cref{thm:labels case 2} that $S_{((n-m),(1^m))}$ has ungraded composition factors $D_{\mu_{n,m-1}}$ and $D_{\mu_{n,m}}$ for all $m\in\{1,\dots,n-1\}$. We now determine the grading shifts $i,j\in\mathbb{Z}$ so that $D_{\mu_{n,m-1}}\langle{i}\rangle$ and $D_{\mu_{n,m}}\langle{j}\rangle$ are graded composition factors of $S_{((n-m),(1^m))}$.

\begin{thm}\label{thm:grdecomp2}
Let $\kappa_2\not\equiv{\kappa_1-1}\Mod{e}$ and $n\equiv{l+1}\Mod{e}$. Then, for all $m\in\{1,\dots,n-1\}$,
	\begin{itemize}
	\item{
	$
	\left[S_{((n-m),(1^m))}:D_{\mu_{n,m-1}}\right]_v=v^{\left(\lfloor\frac{m}{e}\rfloor+
	\lfloor\frac{m+e-2-l}{e}\rfloor+1\right)}$,
	}
	\item{
	$\left[S_{((n-m),(1^m))}:D_{\mu_{n,m}}\right]_v=v^{\left(\lfloor\frac{m}{e}\rfloor+
	\lfloor\frac{m+e-2-l}{e}\rfloor\right)}$.
	}
\end{itemize}
Moreover, $[S_{((n-m),(1^m))}:D_{\mu}]_v=0$ for all other $\mu\in\mathscr{RP}_n^2$.
\end{thm}

\begin{proof}
We determine $x,y\in\mathbb{Z}$ such that
$\grdim (S_{((n-m),(1^m))})=v^x\grdim (D_{\mu_{n,m-1}})+v^y\grdim (D_{\mu_{n,m}})$.
\begin{enumerate}
\item{
Let $0\leqslant{m}\leqslant{\lfloor\frac{n}{e}\rfloor}$.
By \cref{cor:grdim}, the leading and trailing terms, respectively, in the graded dimension of $S_{((n-m),(1^m))}$ are
\[
\binom{\lfloor\frac{n-l-1}{e}\rfloor+1}{m}
v^{\left(
m+\lfloor\frac{m}{e}\rfloor+
\lfloor\frac{m+e-2-l}{e}\rfloor
\right)}
\text{ and }
\binom{\lfloor\frac{n-l-1}{e}\rfloor}{m}
v^{\left(
-m+\lfloor\frac{m}{e}\rfloor+
\lfloor\frac{m+e-2-l}{e}\rfloor
\right)},
\]
and by \cref{grdimD1}, the leading terms in the graded dimensions of $\im (\gamma_{m-1})$ and $\im (\gamma_m)$, respectively, are
\[
\binom{\lfloor\frac{n-l-1}{e}\rfloor}{m-1}v^{m-1}
\text{ and }
\binom{\lfloor\frac{n-l-1}{e}\rfloor}{m}v^{m}.
\]
First observe that the graded dimensions of $D_{\mu_{n,m}}$ and $S_{((n-m),(1^m))}$ both have $2m+1$ terms, and hence $y=\lfloor\frac{m}{e}\rfloor+\lfloor\frac{m+e-2-l}{e}\rfloor$.
Thus $x-\lfloor\frac{m}{e}\rfloor-\lfloor\frac{m+e-2-l}{e}\rfloor$ equals $0$ or $1$ since the trailing coefficients in the graded dimensions of $D_{\mu_{n,m}}$ and $S_{((n-m),(1^m))}$ are equal. Now observe that the sum of the leading coefficients in the graded dimensions of $D_{\mu_{n,m-1}}$ and $D_{\mu_{n,m}}$ equals the leading coefficient in the graded dimension of $S_{((n-m),(1^m))}$. Hence $x=\lfloor\frac{m}{e}\rfloor+\lfloor\frac{m+e-2-l}{e}\rfloor+1$.
}
\item{
Let $\lfloor\frac{n}{e}\rfloor<m<n-\lfloor\frac{n}{e}\rfloor$.
By \cref{cor:grdim}, the leading and trailing terms in the graded dimension of $S_{((n-m),(1^m))}$, respectively, are
\[
\binom{n-2\lfloor\frac{n-l-1}{e}\rfloor-1}{m-\lfloor\frac{n-l-1}{e}\rfloor-1}
v^{\left(
\lfloor\frac{n-l-1}{e}\rfloor+1+\lfloor\frac{m}{e}\rfloor+\lfloor\frac{m+e-2-l}{e}\rfloor
\right)},\ 
\binom{n-2\lfloor\frac{n-l-1}{e}\rfloor-1}{m-\lfloor\frac{n-l-1}{e}\rfloor}
v^{\left(
-\lfloor\frac{n-l-1}{e}\rfloor+\lfloor\frac{m}{e}\rfloor+\lfloor\frac{m+e-2-l}{e}\rfloor
\right)}.
\]
By \cref{grdimD1}, the leading terms in the graded dimensions of $D_{\mu_{n,m-1}}$ and $D_{\mu_{n,m}}$, respectively, are
\[
\binom{n-2\lfloor\frac{n-l-1}{e}\rfloor-1}{m-\lfloor\frac{n-l-1}{e}\rfloor-1}
v^{
\lfloor\frac{n-l-1}{e}\rfloor}
\text{ and }
\binom{n-2\lfloor\frac{n-l-1}{e}\rfloor-1}{m-\lfloor\frac{n-l-1}{e}\rfloor}
v^{
\lfloor\frac{n-l-1}{e}\rfloor}.
\]
Observing that the leading coefficients in the graded dimensions of $S_{((n-m),(1^m))}$ and $D_{\mu_{n,m-1}}$ are equal, we deduce that $x=\lfloor\frac{m}{e}\rfloor+\lfloor\frac{m+e-2-l}{e}\rfloor+1$.
Similarly, observing that the trailing coefficients in the graded dimensions of $S_{((n-m),(1^m))}$ and $D_{\mu_{n,m}}$ are equal, we deduce that $y=\lfloor\frac{m}{e}\rfloor+\lfloor\frac{m+e-2-l}{e}\rfloor$.
}
\item{
Let $\lfloor\frac{n}{e}\rfloor\leqslant{m}\leqslant{n-1}$.
By \cref{cor:grdim}, the leading and trailing terms in the graded dimension of $S_{((n-m),(1^m))}$ are
\[
\binom{\lfloor\frac{n-l-1}{e}\rfloor}{n-m}
v^{\left(
n-m+1+\lfloor\frac{m}{e}\rfloor+\lfloor\frac{m+e-2-l}{e}\rfloor
\right)}
\text{ and }
\binom{\lfloor\frac{n-l-1}{e}\rfloor+1}{n-m}
v^{\left(
m-n+1+\lfloor\frac{m}{e}\rfloor+\lfloor\frac{m+e-2-l}{e}\rfloor
\right)},
\]
respectively, and by \cref{grdimD1}, the leading terms in the graded dimension of $D_{\mu_{n,m-1}}$ and $D_{\mu_{n,m}}$, respectively, are
\[
\binom{\lfloor\frac{n-l-1}{e}\rfloor}{n-m}
v^{\left(
n-m
\right)}
\text{ and }
\binom{\lfloor\frac{n-l-1}{e}\rfloor}{n-m-1}
v^{\left(
n-m-1
\right)}.
\]
First observe that the graded dimensions of $S_{((n-m),(1^m))}$ and $D_{\mu_{n,m-1}}$ both have $2n-2m+1$ terms, and hence $x=\lfloor\frac{m}{e}\rfloor+\lfloor\frac{m+e-2-l}{e}\rfloor+1$.
Thus $y-\lfloor\frac{m}{e}\rfloor-\lfloor\frac{m+e-2-l}{e}\rfloor$ equals $0$ or $1$ since the leading coefficients in the graded dimensions of $S_{((n-m),(1^m))}$ and $D_{\mu_{n,m-1}}$ are equal. Now observe that the sum of the trailing coefficients in the graded dimensions of $D_{\mu_{n,m-1}}$ and $D_{\mu_{n,m}}$ equals the trailing coefficient in the graded dimension of $S_{((n-m),(1^m))}$. Hence, $y=\lfloor\frac{m}{e}\rfloor+\lfloor\frac{m+e-2-l}{e}\rfloor$.
}
\end{enumerate}
\qedhere
\end{proof}

\begin{ex}
	Let $e=3$ and $\kappa=(0,0)$. Then the decomposition submatrix of $\mathscr{H}_7^{\Lambda}$ with rows corresponding to Specht modules labelled by hook bipartitions can be written as
	\[
	\begin{blockarray}{ccccccccccc}
	\begin{block}{c(ccccccc|ccc)}
	S_{((7),\varnothing)}
	&1&&&&&&						&&&\\[0.2em]
	S_{((6),(1))}
	&v&1&&&&\text{\huge${0}$}&						&&&	\\[0.6em]
	S_{((5),(1^2))}
	&&v^2&v&&&&						&&&\\
	S_{((4),(1^3))}
	&&&v^3&v^2&&&					&&\text{\Huge${0}$}&\\[0.6em]
	S_{((3),(1^4))}
	&&&&v^3&v^2&&					&&&\\[0.6em]
	S_{((2),(1^5))}
	&&&&&v^4&v^3&					&&&\\[0.2em]
	S_{((1),(1^6))}
	&&\text{\huge${0}$}&&&&v^5&v^4					&&&\\[0.6em]
	S_{(\varnothing,(1^7))}
	&&&&&&&v^5						&&&\\
	\end{block}
	\end{blockarray}
	\]
\end{ex}

\subsection[Case III: $\kappa_2\equiv{\kappa_1-1}\Mod{e}$ and $n\not\equiv{0}\Mod{e}$]{Case III: $\kappa_2\equiv{\kappa_1-1}\Mod{e}$ and $n\not\equiv{0}\Mod{e}$}

Let $\kappa_2\equiv{\kappa_1-1}\Mod{e}$ and $n\not\equiv{0}\Mod{e}$. Recall from \cref{labels:case 3} that the ungraded composition factors of $S_{((n-m),(1^m))}$ are $D_{\mu_{n,2m}}$ and $D_{\mu_{n,2m+1}}$ for all $m\in\{1,\dots,n-1\}$. Hence as graded $\mathscr{H}_n^{\Lambda}$-modules, the composition factors of $S_{((n-m),(1^m))}$ are $D_{\mu_{n,2m}}\langle{i}\rangle$ and $D_{\mu_{n,2m+1}}\langle{j}\rangle$ for some integers $i$ and $j$, which we now determine.

\begin{thm}
Let $\kappa_2\equiv{\kappa_1-1}\Mod{e}$ and $n\not\equiv{0}\Mod{e}$. Then, for all $m\in\{1,\dots,n-1\}$,
	\begin{itemize}
	\item{
	$\left[S_{((n-m),(1^m))}:D_{\mu_{n,2m}}\right]_v=
	v^{\left(\lfloor\frac{m}{e}\rfloor+
	\lfloor\frac{m+e-1}{e}\rfloor\right)}$,	
	}
	\item{
	$\left[S_{((n-m),(1^m))}:D_{\mu_{n,2m+1}}\right]_v=
	v^{\left(\lfloor\frac{m}{e}\rfloor+
	\lfloor\frac{m+e-1}{e}\rfloor-1\right)}$.
	}
\end{itemize}
Moreover, $[S_{((n-m),(1^m))}:D_{\mu}]_v=0$ for all other $\mu\in\mathscr{RP}_n^2$.
\end{thm}

\begin{proof}
Similar to the proof of \cref{thm:grdecomp2}: we determine $x,y\in\mathbb{Z}$ such that $\grdim(S_{((n-m),(1^m))})=v^x\grdim(D_{\mu_{n,2m}})+v^y\grdim(D_{\mu_{n,2m+1}})$ using \cref{cor:grdim} and \cref{grdimD4}, for each of the three cases $1\leqslant{m}\leqslant{\lfloor\frac{n}{e}\rfloor}$, $\lfloor\frac{n}{e}\rfloor<m<n-\lfloor\frac{n}{e}\rfloor$ and $n-\lfloor\frac{n}{e}\rfloor\leqslant m\leqslant n-1$.
\end{proof}

\begin{ex}
	Let $e=3$ and $\kappa=(0,2)$. Then the decomposition submatrix of $\mathscr{H}_7^{\Lambda}$ with rows corresponding to Specht modules labelled by hook bipartitions can be written as
	\[
	\begin{blockarray}{cccccccccccccccccc}
	\begin{block}{c(cccccccccccccc|ccc)}
	S_{((7),\varnothing)}
	&1&&&&&&&&&&&&&										&&&\\[0.2em]
	S_{((6),(1))}
	&&v&1&&&&&&&&&&\text{\huge${0}$}&					&&&\\[0.6em]
	S_{((5),(1^2))}
	&&&&v&1&&&&&&&&&									&&&\\
	S_{((4),(1^3))}
	&&&&&&v^2&v&&&&&&&									 &&\text{\Huge${0}$} &\\[0.6em]
	S_{((3),(1^4))}
	&&&&&&&&v^3&v^2&&&&&								&&&\\[0.6em]
	S_{((2),(1^5))}
	&&&&&&&&&&v^3&v^2&&&								&&&\\[0.2em]
	S_{((1),(1^6))}
	&&\text{\huge${0}$}&&&&&&&&&&v^4&v^3&				&&&\\[0.6em]
	S_{(\varnothing,(1^7))}
	&&&&&&&&&&&&&&v^4									&&&\\
	\end{block}
	\end{blockarray}
	\]
\end{ex}

\subsection[Case IV: $\kappa_2\equiv{\kappa_1-1}\Mod{e}$ and $n\equiv{0}\Mod{e}$]{Case IV: $\kappa_2\equiv{\kappa_1-1}\Mod{e}$ and $n\equiv{0}\Mod{e}$}

Let $\kappa_2\equiv{\kappa_1-1}\Mod{e}$ and $n\equiv{0}\Mod{e}$. Recall from \cref{labels:case 4} that $D_{\mu_{n,2m}}$, $D_{\mu_{n,2m+2}}$, $D_{\mu_{n,2m+1}}$ and $D_{\mu_{n,2m+3}}$ are the ungraded composition factors of $S_{((n-m),(1^m))}$ for all $m\in\{2,\dots,n-2\}$; $S_{((n-1),(1))}$ and $S_{((1),(1^n))}$ both have three composition factors. Hence as graded $\mathscr{H}_n^{\Lambda}$-modules, $S_{((n-m),(1^m))}$ has composition factors $D_{\mu_{n,2m}}\langle{i_1}\rangle$, $D_{\mu_{n,2m+2}}\langle{i_2}\rangle$, $D_{\mu_{n,2m+1}}\langle{i_3}\rangle$ and $D_{\mu_{n,2m+3}}\langle{i_4}\rangle$ for some $i_1,i_2,i_3,i_4\in\mathbb{Z}$, which we now determine.

Firstly, one observes that the graded dimension of $D_{\mu_{n,2m}}$ equals the graded dimension of $D_{\mu_{n,2m+3}}$, under a grading shift, which follows immediately from \cref{grdimD2}.

\begin{lem}\label{eqgrdim}
Let $\kappa_2\equiv{\kappa_1-1}\Mod{e}$ and $n\equiv{0}\Mod{e}$. Then, for all $m\in\{1,\dots,n-2\}$, 
\[
v^2
\left[S_{((n-m),(1^{m}))}:D_{\mu_{n,2m+3}}\right]_v
=
\left[S_{((n-m),(1^m))}:D_{\mu_{n,2m}}\right]_v.
\]
\end{lem}

\begin{thm}
Suppose that $\kappa_2\equiv{\kappa_1-1}\Mod{e}$ and $n\equiv{0}\Mod{e}$, and let $1\leqslant m<n$. Then
	\begin{itemize}
	\item{
	$\left[S_{((n-m),(1^m))}:D_{\mu_{n,2m}}\right]_v=
	v^{\left(\lfloor\frac{m}{e}\rfloor
	+\lfloor\frac{m+e-1}{e}\rfloor+1\right)}$
	for all $m\in\{1,\dots,n-1\}$,
	}
	\item{
	$\left[S_{((n-m),(1^m))}:D_{\mu_{n,2m+2}}\right]_v=
	v^{\left(\lfloor\frac{m}{e}\rfloor
	+\lfloor\frac{m+e-1}{e}\rfloor\right)}	
	$
	for all $m\in\{1,\dots,n-2\}$,
	}
	\item{
	$\left[S_{((n-m),(1^m))}:D_{\mu_{n,2m+1}}\right]_v=
	v^{\left(\lfloor\frac{m}{e}\rfloor
	+\lfloor\frac{m+e-1}{e}\rfloor\right)}	
	$
	for all $m\in\{2,\dots,n-2\}$,
	}
	\item{
	$\left[S_{((n-m),(1^m))}:D_{\mu_{n,2m+3}}\right]_v=
	v^{\left(\lfloor\frac{m}{e}\rfloor
	+\lfloor\frac{m+e-1}{e}\rfloor-1\right)}	
	$
	for all $m\in\{1,\dots,n-2\}$.
	}
\end{itemize}
Moreover, $[S_{((n-m),(1^m))}:D_{\mu}]_v=0$ for all other $\mu\in\mathscr{RP}_n^2$.
\end{thm}

\begin{proof}
\begin{enumerate}
\item{
Let $m=1$. We know from \cref{labels:case 4} that $S_{((n-1),(1))}$ has three composition factors, namely $D_{\mu_{n,2}}=S_{((n),\varnothing)}$, $D_{\mu_{n,4}}$ and $D_{\mu_{n,5}}$.
It follows thus from \cref{eqgrdim} that, for some $x,y\in\mathbb{Z}$, we have
\begin{align*}
&
\text{grdim}\left(S_{((n-1),(1))}\right)
=
v^x \grdim \left(S_{((n),\varnothing)}\right)
+ v^y \grdim \left(D_{\mu_{n,4}}\right)
+ v^{x-2} \grdim \left(D_{\mu_{n,5}}\right).
\end{align*}
Furthermore, one determines that
$S_{((n),\varnothing)}\cong \langle v(n)\rangle$ and
$D_{\mu_{n,5}}\cong\ker(\gamma_2)/\im(\phi_2)\cong\langle v(1)\rangle$ as ungraded $\mathscr{H}_n^{\Lambda}$-modules, and thus $\grdim\left(S_{((n),\varnothing)}\right)=\grdim\left(D_{\mu_{n,5}}\right)=1$.
Hence, by \cref{cor:grdim} and \cref{grdimD2}, we have
\begin{align*}
\grdim \left(S_{((n-1),(1))}\right)
&=
\tfrac{n}{e}v^2 + \tfrac{(e-2)n}{e}v + \tfrac{n}{e}
=
v^{2x-2} + v^y \left( \tfrac{n-e}{e}v + \tfrac{(e-2)n}{e} + \tfrac{n-e}{e}v^{-1} \right).
\end{align*}
Thus, by equating terms, $y=1=x-1$.
}
\item{
Let $1<m<n-1$. We know from \cref{labels:case 4} that $S_{((n-m),(1^m))}$ has four composition factors, $D_{\mu_{n,2m}}$, $D_{\mu_{n,2m+1}}$, $D_{\mu_{n,2m+2}}$ and $D_{\mu_{n,2m+3}}$.
Following \cref{eqgrdim}, we know that
\begin{align*}
\grdim \left(S_{((n-m),(1^m))}\right)=
&v^x\grdim \left(D_{\mu_{n,2m}}\right)
+v^y\grdim \left(D_{\mu_{n,2m+1}}\right)\\
&\quad+v^z\grdim \left(D_{\mu_{n,2m+2}}\right)
+v^{x-2}\grdim \left(D_{\mu_{n,2m+3}}\right)
\end{align*}
for some $x,y,z\in\mathbb{Z}$, which we now determine.

Firstly, let $2\leqslant{m}\leqslant{\frac{n}{e}}$.
Let $\eta=\lfloor\frac{m}{e}\rfloor+\lfloor\frac{m+e-1}{e}\rfloor$. Then it follows from \cref{cor:grdim} that the leading and trailing terms of $\grdim \left(S_{((n-m),(1^m))}\right)$ are as follows.

	\begin{center}
	\begin{tabular}{ c | c | c | c}
		$\nth{1}$ term & $\nth{2}$ term & $\nth{2}$ last term & last term \\
		\hline\hline
		\rule{0pt}{4ex}  $\binom{\frac{n}{e}}{m}
		v^{\left(m+\eta\right)}$
		& \rule{0pt}{4ex} $\tfrac{(e-2)n}{e} \tbinom{\frac{n}{e}}{m-1}
		v^{\left(m-1+\eta\right)}$
		& \rule{0pt}{4ex} $\tfrac{(e-2)n}{e}\binom{\frac{n}{e}}{m-1}
		v^{\left(1-m+\eta\right)}$ 
		& \rule{0pt}{4ex} $\binom{\frac{n}{e}}{m}
		v^{\left(-m+\eta\right)}$
	\end{tabular}
\end{center}

\begin{comment}
By \cref{cor:grdim}, the first two leading terms of $\grdim \left(S_{((n-m),(1^m))}\right)$ are 
\[
\tbinom{\frac{n}{e}}{m}
v^{\left(m+\lfloor\frac{m}{e}\rfloor+\lfloor\frac{m+e-1}{e}\rfloor\right)}
\text{ and }
\tfrac{(e-2)n}{e}
\tbinom{\frac{n}{e}}{m-1}
v^{\left(m-1+\lfloor\frac{m}{e}\rfloor+\lfloor\frac{m+e-1}{e}\rfloor\right)},
\]
respectively, 
and the last two trailing terms of $\grdim \left(S_{((n-m),(1^m))}\right)$ are 
\[
\tfrac{(e-2)n}{e}
\tbinom{\frac{n}{e}}{m-1}
v^{\left(1-m+\lfloor\frac{m}{e}\rfloor+\lfloor\frac{m+e-1}{e}\rfloor\right)}
\text{ and }
\binom{\frac{n}{e}}{m}
v^{\left(-m+\lfloor\frac{m}{e}\rfloor+\lfloor\frac{m+e-1}{e}\rfloor\right)},
\]
respectively.
\end{comment}

By \cref{grdimD2}, the first two leading terms in the graded dimensions of $D_{\mu_{n,2m}}$, $D_{\mu_{n,2m+1}}$, $D_{\mu_{n,2m+2}}$ and $D_{\mu_{n,2m+3}}$ are presented in the following table.

	\begin{center}
		\begin{tabular}{ c | | c | c | c | c}
			& $D_{\mu_{n,2m}}$ & $D_{\mu_{n,2m+1}}$ & $D_{\mu_{n,2m+2}}$ & $D_{\mu_{n,2m+3}}$ \\
			\hline\hline
			\rule{0pt}{4ex} $\nth{1}$ term
			& \rule{0pt}{4ex}  $\binom{\frac{n-e}{e}}{m-1}v^{m-1}$
			& \rule{0pt}{4ex} $\binom{\frac{n-e}{e}}{m-2}v^{m-2}$
			& \rule{0pt}{4ex} $\binom{\frac{n-e}{e}}{m}v^m$ 
			& \rule{0pt}{4ex} $\binom{\frac{n-e}{e}}{m-1}v^{m-1}$ \\
			\rule{0pt}{4ex} $\nth{2}$ term
			& \rule{0pt}{4ex} $\frac{(e-2)n}{e}	\binom{\frac{n-e}{e}}{m-2}v^{m-2}$
			& \rule{0pt}{4ex} $\frac{(e-2)n}{e}	\binom{\frac{n-e}{e}}{m-3}v^{m-3}$
			& \rule{0pt}{4ex} $\frac{(e-2)n}{e}	\binom{\frac{n-e}{e}}{m-1}v^{m-1}$
			& \rule{0pt}{4ex} $\frac{(e-2)n}{e}	\binom{\frac{n-e}{e}}{m-2}v^{m-2}$
		\end{tabular}
	\end{center}

\begin{comment}
By \cref{grdimD2}, the first two leading terms of the graded dimensions of
$D_{\mu_{n,2m}}$ and 
$D_{\mu_{n,2m+2}}$
are
\[
\binom{\frac{n-e}{e}}{m-1}v^{m-1},\ 
\frac{(e-2)n}{e}
\binom{\frac{n-e}{e}}{m-2}v^{m-2}\]
and
\[
\binom{\frac{n-e}{e}}{m}v^m,\
\frac{(e-2)n}{e}
\binom{\frac{n-e}{e}}{m-1}v^{m-1},
\]
respectively, and the first two leading terms of 
$D_{\mu_{n,2m+1}}$
and
$D_{\mu_{n,2m+3}}$
are
\[
\binom{\frac{n-e}{e}}{m-2}v^{m-2},\
\frac{(e-2)n}{e}
\binom{\frac{n-e}{e}}{m-3}v^{m-3}\]
and
\[
\binom{\frac{n-e}{e}}{m-1}v^{m-1},
\frac{(e-2)n}{e}
\binom{\frac{n-e}{e}}{m-2}v^{m-2},
\]
respectively.
\end{comment}

The graded dimensions of $S_{((n-m),(1^m))}$ and $D_{\mu_{n,2m+2}}$ both have $2m+1$ terms, and hence $z=\lfloor\frac{m}{e}\rfloor+\lfloor\frac{m+e-1}{e}\rfloor$. Now observe that the graded dimensions of $D_{\mu_{n,2m}}$ and $D_{\mu_{n,2m+3}}$ both have $2m-1$ terms, so together with \cref{eqgrdim}, $x=\lfloor\frac{m}{e}\rfloor+\lfloor\frac{m+e-1}{e}\rfloor+1$.

We thus have $-2\leqslant{y-\lfloor\frac{m}{e}\rfloor-\lfloor\frac{m+e-1}{e}\rfloor}\leqslant{2}$, and observe that the sum of the second leading (trailing, respectively) coefficients in the graded dimensions of $D_{\mu_{n,2m}}$ and $D_{\mu_{n,2m+2}}$ ($D_{\mu_{n,2m+2}}$ and $D_{\mu_{n,2m+3}}$, respectively) form the second leading (trailing, resp.) coefficient in the graded dimension of $S_{((n-m),(1^m))}$. Hence $y=\lfloor\frac{m}{e}\rfloor+\lfloor\frac{m+e-1}{e}\rfloor$.

We similarly find $x,y,z$ using \cref{cor:grdim} and \cref{grdimD2} for $\frac{n}{e}<m<\frac{n(e-1)}{e}$ and $\frac{n(e-1)}{e}\leqslant{m}\leqslant{n-2}$, respectively.
}
\item{
Let $m=n-1$. We know from \cref{labels:case 4} that $S_{((1),(1^{n-1}))}$ has three composition factors, namely $D_{(\varnothing,(1^n))^R}$, $D_{\mu_{n,2n-1}}$ and $D_{\mu_{n,2n-2}}$.
One determines that $D_{(\varnothing,(1^n))^R}\cong \langle v(1,2,\dots,n-1)\rangle$ and $D_{\mu_{n,2n-2}}\cong \im(\phi_{n-1})=\langle v(2,3,\dots,n)\rangle$ as ungraded $\mathscr{H}_n^{\Lambda}$-modules, and hence $\grdim\left(D_{(\varnothing,(1^n))^R}\right)=\grdim\left(D_{\mu_{n,2n-2}}\right)=1$. Moreover, we find that
\[
\deg \left( v(2,3,\dots,n) \right) =
2\lfloor\tfrac{m}{e}\rfloor+2
=\deg \left( v(1,2,\dots,n-1) \right) +2,
\]
so that
\[
v^2 \left[S_{((1),(1^{n-1}))}:D_{(\varnothing,(1^n))^R}\right]_v = \left[S_{((1),(1^{n-1}))}:D_{\mu_{n,2m-2}}\right]_v.
\]
It thus follows that
\begin{align*}
&\grdim \left(S_{((1),(1^{n-1}))}\right)
=v^x \grdim \left(D_{\mu_{n,2n-2}}\right)
+ v^y \grdim \left(D_{\mu_{n,2n-1}}\right)
+ v^{x-2} \grdim \left(D_{(\varnothing,(1^n))^R}\right)
\end{align*}
for some $x,y\in\mathbb{Z}$, which we now determine.
Applying \cref{cor:grdim} and \cref{grdimD2},
\begin{align*}
\grdim \left(S_{((1),(1^{n-1}))} \right)
= &
\tfrac{n}{e}
v^{\left( 1 + \lfloor\frac{m}{e}\rfloor + \lfloor\frac{m+e-1}{e}\rfloor \right)}
+ \tfrac{(e-2)n}{e}
v^{\left( \lfloor\frac{m}{e}\rfloor + \lfloor\frac{m+e-1}{e}\rfloor \right)}
+ \tfrac{n}{e}
v^{\left( \lfloor\frac{m}{e}\rfloor + \lfloor\frac{m+e-1}{e}\rfloor - 1 \right)}\\
= &
v^{2x-2} + v^y \left( \tfrac{n-e}{e}v + \tfrac{(e-2)n}{e} + \tfrac{n-e}{e}v^{-1} \right).
\end{align*}
Equating terms, we deduce that $y=\lfloor\frac{m}{e}\rfloor+\lfloor\frac{m+e-1}{e}\rfloor=x-1$, as required.
}
\end{enumerate}
\qedhere
\end{proof}

\begin{ex}
	Let $e=3$ and $\kappa=(0,2)$. Then the decomposition submatrix of $\mathscr{H}_6^{\Lambda}$ with rows corresponding to Specht modules labelled by hook bipartitions can be written as
	\[
	\begin{blockarray}{cccccccccccccc}
	\begin{block}{c(cccccccccc|ccc)}
	S_{((6),\varnothing)}
	&1&&&&&&&&&								&&&\\[0.2em]
	S_{((5),(1))}
	&v^2&v&1&&&&&&\text{\huge${0}$}&			&&&\\[0.6em]
	S_{((4),(1^2))}
	&&v^2&v&v&1&&&&&						&&&\\
	S_{((3),(1^3))}
	&&&&v^3&v^2&v^2&v&&&					&&\text{\Huge${0}$}&\\[0.6em]
	S_{((2),(1^4))}
	&&&&&&v^4&v^3&v^3&v^2&					&&&\\[0.2em]
	S_{((1),(1^5))}
	&&\text{\huge${0}$}&&&&&&v^4&v^3&v^2		&&&\\[0.6em]
	S_{(\varnothing,(1^6))}
	&&&&&&&&&&v^4							&&&\\
	\end{block}
	\end{blockarray}
	\]
\end{ex}

\renewcommand{\abstractname}{Acknowledgements}
\begin{abstract}
	This paper was written under the guidance of the author's PhD supervisor, Matthew Fayers, at Queen Mary University of London, and forms part of her PhD thesis. The author would like to thank Dr Fayers for his many helpful comments and ongoing support, as well as Chris Bowman and Liron Speyer for their useful remarks and guidance. The author is also thankful to the referee for their careful reading of the manuscript.
\end{abstract}

\bibliographystyle{plain}
\bibliography{master}

\end{document}